\input ./style/arxiv-general.cfg
\documentclass[aop,MSNbibl,seceqn,dvips]{arximspdf}
\makeatletter
   \@ifpackageloaded{graphicx}{}{\usepackage{graphicx}}
\makeatother
\usepackage{mathrsfs,mathbh}
\doi{10.1214/15-AOP1018}% Updated by VTEXPTS2LaTeX.exe, 07.04.2015
\volume{44}
\issue{3}
\pubyear{2016}
\firstpage{2064}
\lastpage{2116}
\docsubty{FLA}

\makeatletter
\newcommand{\rrVert}{\Vert}
\newcommand{\llVert}{\Vert}
\newtheorem{theorem}{Theorem}[section]
\newtheorem{lemma}{Lemma}[section]
\newproclaim{remark}{Remark}[section]
\newproclaim{definition}{Definition}[section]
\newtheorem{proposition}{Proposition}[section]
\newproclaim{Examples}{Examples}[section]
\newproclaim{example}{Example}
\newtheorem{corollary}{Corollary}[section]
\def\e{\varepsilon}
\def\s{\sigma}
\def\a{\alpha}
\def\Om{\Omega}
\def\vph{\varphi}
\def\p{\partial}
\def\d{{\mathord{{\mathrm d}}}}
\def\tr{\operatorname{tr}}
\def\cC{{\mathcal C}}
\def\cF{{\mathcal F}}
\def\cH{{\mathcal H}}
\def\cN{{\mathcal N}}
\def\cS{{\mathcal S}}

\newcommand{\mL}{{\mathbb L}}
\newcommand{\mN}{{\mathbb N}}
\newcommand{\mR}{{\mathbb R}}
\newcommand{\mS}{{\mathbb S}}
\newcommand{\mW}{{\mathbb W}}
\newcommand{\mX}{{\mathbb X}}
\def\bP{{\mathbf P}}

\def\sF{{\mathscr F}}
\def\sP{{\mathscr P}}
\def\sS{{\mathscr S}}

\def\bE{{\mathbf E}}
\makeatother

\begin{document}
\begin{frontmatter}

\title{On approximate continuity and the support of reflected
stochastic differential equations\thanksref{T1}}
\runtitle{On approximate continuity and the support of RSDEs}
\thankstext{T1}{Supported by NSFC (Nos. 11171358,11301553, 11471340)
and the Fundamental Research Funds for the Central Universities (No. 13lgpy64).}

\begin{aug}
% Corresponding author: Jing Wu - wjjosie@hotmail.com% Updated by
%VTEXPTS2LaTeX.exe, 07.04.2015 13:48
\author[A]{\fnms{Jiagang}~\snm{Ren}\ead[label=e1]{renjg@mail.sysu.edu.cn}}
\and
\author[A]{\fnms{Jing}~\snm{Wu}\corref{}\ead[label=e2]{wjjosie@hotmail.com}}
%\author[A]{\fnms{}~\snm{}\corref{}\ead[label=e1]{}}%,
%\author[]{\fnms{}~\snm{}\ead[label=]{}}
% \and
%\author[]{\fnms{}~\snm{}\ead[label=]{}}
\runauthor{J. Ren and J. Wu}
\affiliation{Sun Yat-sen University}
%\dedicated{}
\address[A]{School of Mathematics and Computational Science\\
Sun Yat-sen University\\
510275 Guangzhou \\
P.R. China\\
\printead{e1}\\
\phantom{E-mail:\ }\printead*{e2}}
%\address[A]{\\\printead{e1}}
%\address[]{\\\printead{}}
\end{aug}

% HISTORY:
%
\received{\smonth{10} \syear{2014}}% Updated by VTEXPTS2LaTeX.exe,
%07.04.2015 13:48
%
\revised{\smonth{3} \syear{2015}}% Updated by VTEXPTS2LaTeX.exe,
%07.04.2015 13: ABSTRACT
%
\begin{abstract}
In this paper we prove an approximate continuity result for stochastic
differential equations with normal reflections in domains satisfying
Saisho's conditions, which together with the Wong--Zakai approximation
result completes the support theorem for such diffusions in the uniform
convergence topology. Also by adapting Millet and Sanz-Sol\'e's idea,
we characterize in H\"older norm the support of diffusions reflected in
domains satisfying the Lions--Sznitman conditions by proving limit
theorems of adapted interpolations. Finally we apply the support
theorem to establish a boundary-interior maximum principle for
subharmonic functions.
\end{abstract}

% KEYWORDS

%
%\begin{keyword}[class=AMS]
%\kwd[Primary ]{}
%\kwd{}
%\kwd[; secondary ]{}
%\end{keyword}
%\begin{keyword}
%\kwd{}
%\end{keyword}
\end{frontmatter}

%s1 #&#
\section{Introduction}\label{sec1}
The support theorem for diffusion processes defined by stochastic
differential equations
has been a much studied topic for probabilists and analysts since the
seminal work
of Stroock and Varadhan \cite{str}. The typical approach to a support
theorem in the norm of uniform convergence consists of two steps.
One step is to establish a limit theorem for SDEs, meaning that the
solution of an SDE can be approximated
by a sequence of solutions of ODEs, obtained by regularizing the
Brownian paths \cite{wong2}; the other
is to prove a Denjoy-type approximate continuity theorem, stating that
the solution of an SDE
is approximately continuous at points in a dense set of the
Cameron--Martin space. Millet and Sanz-Sol\'e \cite{millet2,millet}
proposed a simple
approach to characterizing in H\"older spaces the support of diffusions
described by general SDEs, obtained by approximating Brownian motions
with linear adapted interpolations, and proved the two inclusions
through approximation results.

In this work we are concerned with the support problem of diffusions
constrained in a domain
$D$ with normal reflection boundary. Such diffusions have been
constructed by Anderson and Orey \cite{orey} if
$D$ has smooth boundary and by Tanaka \cite{tanaka} if $D$ is convex.
Correspondingly the support theorem has been established by Doss and
Priouret \cite{doss}
if $D$ has smooth boundary, and a limit theorem has been proved by
Pettersson \cite{pettersson} when
$D$ is a convex domain and the diffusion coefficient is constant. Recently
in \cite{rx}, a support theorem was proved for stochastic variational
inequalities; this means, in particular, that the support theorem holds
true for diffusions normally reflected in convex domains.

However, normally reflected diffusions have been constructed for
domains much wider than convex domains
and smooth domains (see Lions and Sznitman~\cite{ls} and Saisho \cite
{s}), so a natural (and application-motivated) question is whether or
not the
support theorem continues to hold true for such diffusions. The first
step in this respect was taken by
Evans and Stroock \cite{stroock} who proved, under the set of
conditions given by Lions and Sznitman, that a weak limit theorem
holds. Very recently this result was improved by Aida and Sasaki \cite
{aida}, and independently by Zhang \cite{zhangtusheng}, who used an
adapted version of the Wong--Zakai approximations rather than the usual
one, by removing the admissibility condition from the
set of conditions and proving that the convergence takes place, in
fact, in $L^p$ (and they obtained the convergence speed). Roughly
speaking, they proved a strong limit theorem for the reflected
diffusions studied by Saisho in \cite{s}. To date, this was the widest,
well-studied situation.

On the other hand, however, approximate continuity has not yet been
touched in such situations. Our first result fills this gap, and it,
together with the Wong--Zakai convergence result in \cite{aida} and
\cite{zhangtusheng}, will yield the support theorem in the locally
uniform convergence topology for normally reflected SDEs in domains,
satisfying the conditions of Lions and Sznitman \cite{ls},
except the admissibility. The second contribution of this paper is to
present a characterization of the support for reflected diffusions in
H\"older spaces in domains satisfying the conditions in \cite{stroock},
by extending the idea of Millet and Sanz-Sol\'e \cite{millet} to SDEs
with normal reflections.

We recall the Skorohod problem here. Let $D$ be a domain in $\mR^d$ and
$w_\cdot\in\mathcal{C}([0,+\infty);\mR^d)$ such that $w_0\in\bar
{D}$. A
pair of continuous functions $(x,k)$ is a solution of the Skorohod
problem if:
\begin{itemize}
\item$x_t\in\bar{D}$ for all $t\geq0$ and $x_0=w_0$;
\item for all $t\geq0$, $x_t=w_t+k_t$;
\item$k(0)=0$, and $k$ is of bounded variation on each finite interval
and satisfies
\[
k_t=\int_0^t n_s
\,\d|k|_s, \qquad |k|_t=\int_0^t
1_{\partial D}(x_s)\,\d|k|_s,
\]
where $n_s\in\mathcal{N}_{x_s}$ and $\mathcal{N}_{x}$ is the set of
inward normal unit vectors at $x\in\partial D$ defined by
\begin{eqnarray*}
\mathcal{N}_x&=&\bigcup_{r>0}
\mathcal{N}_{x,r},
\\
\mathcal{N}_{x,r}&=&\bigl\{n\in\mR^d; |n|=1, B(x-rn, r)
\cap D=\varnothing\bigr\}.
\end{eqnarray*}
\end{itemize}
Here and in what follows $B(a,r)=\{y\in\mR^d; |y-a|<r\},  a\in\mR^d,
 r>0$ and $|k|_t$ denotes the total variation of
$k$ on $[0,t]$.

Let $\Om=\cC_0([0,\infty),\mR^{d_1})$ be the space consisting of
continuous functions
from $[0,\infty)$ to $\mR^{d_1}$ vanishing at $0$. Let $\cF$ be the
completion of the Borel $\sigma$-algebra
on $\Om$ associated with the locally uniform convergence topology and
$\bP$ the distribution
of an $d_1$-dimensional Brownian motion.
Then $(\Om,\cF,\bP)$ is a complete probability space, and the
coordinate process
\[
w_t(\omega):=\omega(t),\qquad t\geq0
\]
is a $d_1$-dimensional standard Brownian motion. The natural filtration
generated by $(w_t)_{t\geq0}$ is denoted by $(\sF_t)_{t\geq0}$.

We consider the following reflected SDE:
%
%e1.1 #&#
\begin{equation}\qquad
\label{msde1} \cases{ %
 \displaystyle X_t=X_0+
\int_0^t\sigma(X_s)\circ\d
w_s+\int_0^tb(X_s)\,\d
s+K_t,&\quad $X_0=x\in\bar{D},$
\vspace*{2pt}\cr
\displaystyle |K|_t=\int_0^t
1_{\partial D}(X_s)\,\d|K|_s,&\quad $\displaystyle K_t=\int
_0^t \xi_s\,\d|K|_s,$}\hspace*{-6pt}
\end{equation}
where $\xi_s\in\cN_{X_s}$. In It\^o's notation, it takes the
following form:
\begin{eqnarray*}
\cases{ %
\displaystyle  X_t=X_0+\int
_0^t\sigma(X_s) \,\d
w_s+\int_0^t\tilde{b}(X_s)
\,\d s+K_t, &\quad $X_0=x\in\bar{D},$
\vspace*{2pt}\cr
\displaystyle|K|_t=\int_0^t
1_{\partial D}(X_s)\,\d|K|_s,&\quad $\displaystyle K_t=\int
_0^t \xi_s\,\d|K|_s$}
\end{eqnarray*}
with
\[
\tilde{b}^i(x):=b^i(x)+\frac{1}2\sum
_{j=1}^d\sum_{k=1}^{d_1}
\bigl[\s _k^i(x)\bigr]_j\s^j_k(x).
\]
Throughout the paper we will assume that $\sigma\dvtx  \mR^d\mapsto\mR
^d\otimes\mR^{d_1}$ and $b\dvtx \mR^d\mapsto\mR^d$ are
$\cC^2_b$ and $\cC^1_b$ functions, respectively. Then by Saisho \cite{s}
this equation has a unique solution
$(X,K)$.

Let $\mW^d$ (resp., $\mW^{d_1}$) denote the space of all $\mR^d$ (resp.,
$\mR^{d_1}$)-valued continuous functions defined on $[0,\infty)$,
and for each $\a\in(0,\frac{1}2)$, $\mW_\a^d$ denote the subspace of
$\mW
^d$ consisting of locally $\a$-H\"older continuous functions. Then for
every $\a\in[0,\frac{1}2)$, $\mW_\a^d$ is a Fr\'echet space with the
topology defined by the system of seminorms $\{\|\cdot\|_{T,\a}, T>0\}
$, where for $x\in\mW^d$,
\begin{eqnarray*}
\|x\|_{T}:=\sup_{0\leq t\leq T}|x_t|,\qquad \|x
\|_{T,\a}:=\|x\|_{T}+\sup_{0\leq s,t\leq T,s\neq t}
\frac
{|x(t)-x(s)|}{|t-s|^\a}.
\end{eqnarray*}
Denote
\begin{eqnarray*}
\cH&:=&\bigl\{h\dvtx h\in\mW^{d_1}; h(0)=0, h(\cdot) \mbox{ is
absolutely continuous and}\\
&&\hspace*{122pt} \dot{h}\in\mL^2\bigl([0,\infty);\mR^{d_1}\bigr), \forall
T>0\bigr\},
\\
\cS&:=&\bigl\{h\in\mW^{d_1}; h(0)=0, t\rightarrow h(t) \mbox{ is smooth}
\bigr\} ,
\\
\cS_p&:=&\bigl\{h\in\mW^{d_1}; h(0)=0, t\rightarrow h(t)
\mbox{ is piecewise smooth}\bigr\}.
\end{eqnarray*}
$\cH$ will be endowed with the topology given by the family of
seminorms $\{\|h\|_{\cH_T}:=(\int_0^T|\dot{h}_t|^2\,\d t)^{1/2}, T>0\}$.
Given $h\in\cH$, denote by $(Z(h), \psi(h))$ the solution to the
following deterministic Skorohod problem:
%
%e1.2 #&#
\begin{equation}
\label{skorohodZ} Z_t=x+\int_0^t
\sigma(Z_s)\dot{h}_s\,\d s+\int_0^tb(Z_s)
\,\d s+\psi_t.
\end{equation}
Let
\[
\sS(\cH):=\bigl\{Z(h), h\in\cH\bigr\}; \qquad \sS:=\bigl\{Z(h), h\in\cS\bigr\};
\qquad\sS_p:=\bigl\{ Z(h), h\in\cS_p\bigr\}.
\]
Denote by $\overline{\sS(\cH)}^\alpha$ the closure of $
\sS(\cH)$ in $\mW_\alpha^d$, and $\overline{\sS}$, $\overline
{\sS_p}$
and $\overline{\sS(\cH)}$ the closures of $\sS$, $\sS_p$
and $\sS(\cH)$ in $\mW^d$, respectively.
We are going to prove in Section~\ref{sec2} the approximate continuity theorem,
which together with the result in \cite{aida} and \cite{zhangtusheng}
yields that the support of $\bP\circ X^{-1}$ in $\mW^d$ coincides with
$\overline{\sS}$. We also prove in Section~\ref{sec3} an enhanced version of the
support theorem by showing that for every $\a\in(0,\frac{1}2)$, the
support of $\bP\circ X^{-1}$ in $\mW_\a^d$ coincides with $\overline
{\sS
(\cH)}^\a$.

The paper is organized as follows: in Section~\ref{sec2} an approximate
continuity theorem for normally reflected diffusions is proved,
and this result combined with the main result in \cite{aida} and \cite
{zhangtusheng} implies, of course, the
support theorem for such diffusions. Next, we provide in Section~\ref{sec3} an
alternate approach to solving the support problem in H\"{o}lder spaces. Finally in Section~\ref{sec4}, we give a first application of our support theorem to maximum
principle for $L$-subharmonic functions in domains having nonsmooth
boundaries and with possibly degenerate $L$.

Throughout the paper we use $C$ to denote a generic constant which may
be different in different
places, and we use summation convention for repeated indices. Finally
$A\lesssim B$ means that
there exists a $C\geq0$ such that $A\leq CB$.

%s2 #&#
\section{Approximate continuity}\label{sec2}
In this section we will work in the setup of \cite{ls}. But, as in
\cite
{s}, we will not need the admissibility condition on the domain.
Precisely, we assume that we are given a domain $D\subset\mR^d$ satisfying:

(H$_1$) There exists $c_0>0$ such that for any $x\in
\partial
D$, $y\in\bar{D}$
and $\xi\in\cN_x$,
\[
(y-x,\xi)+c_0|x-y|^2\geq0,
\]
where $\cN_x$ denotes the set of unit inward normals at $x$;

(H$_2$) There exist a function $\varphi\in\cC_b^3(\mR^d;\mR)$ and
a constant $\alpha>0$ such that
\[
D\varphi(x)\cdot\xi\geq\alpha c_0\qquad\forall x\in\partial D,\xi \in
\cN_x.
\]

It is obvious that under the conditions(H$_1$)--(H$_2$), $\overline{\sS}=\overline{\sS_p}=\overline{\sS(\cH
)}$. To
see this,
we only need to show $\overline{\sS}\supset\overline{\sS(\cH)}$. In
fact, for any $h\in\cH$, we can take a sequence $h^n\in\cS$ such that
$h^n\rightarrow h$ in $\cH$. Denote by $(Z,\Psi)$ and $(Z^n,\Psi^n)$
the corresponding solutions of the Skorohod problem (\ref{skorohodZ}).
Set $\rho(t):=e^{-({2}/\alpha)(\varphi(Z_t)+\varphi(Z^n_t))}$.
Then for
any $t\geq0$, by (H$_2$) and the assumptions $b\in\cC_b^1$
and $\sigma\in\mathcal{C}_b^2$, we have
\begin{eqnarray*}
&&\bigl|Z^n_t-Z_t\bigr|^2e^{-({2}/\alpha)(\varphi(Z_t)+\varphi(Z^n_t))}
\\
&&\qquad \leq C\int_0^t\rho(s)\bigl|Z^n_s-Z_s\bigr|^2
\bigl(1+\bigl|\dot{h}^n_s\bigr|+|\dot{h}_s|\bigr)\,\d
s+C\int_0^t\bigl|\dot{h}^n_s-
\dot{h}_s\bigr|^2\,\d s,
\end{eqnarray*}
which implies by Gronwall's lemma that $\sup_{0\leq t\leq
T}|Z^n_t-Z_t|^2\rightarrow0$ as $n\rightarrow\infty$ and thus $Z\in
\overline{\sS}$, yielding that $\overline{\sS}\supset\overline
{\sS(\cH)}$.

Before we proceed, a few words about these conditions are in order. The
constant $c_0$ appearing in condition (H$_1$) is also allowed
to equal to zero in \cite{ls}.
Then the function $\varphi$ in condition (H$_2$) can be taken
to be identically zero, and it turns out that some arguments below will
break down, and different treatments will be needed. But in this case
$D$ is a convex domain, and thus the equation is a special case of
stochastic variational inequalities
already treated in \cite{rx}. Hence we simply assume $c_0>0$ here.

For convenience we record here some basic facts which will be used
below; see~\cite{iw}.
Set for $i,j=1,\ldots,d_1$,
\begin{eqnarray*}
\kappa^{ij}(t):=\frac{1}{2}\int_0^t
\bigl[w^i_s\,\d w^j_s-w^j_s
\,\d w^i_s\bigr],\qquad \zeta^{ij}(t):=\int
_0^tw^i_s\circ\,\d
w^j_s.
\end{eqnarray*}

Let $T>0$ be arbitrarily fixed.

%le2.1 #&#
\begin{lemma}\label{w}
\textup{(i)} There exist two positive constants $c_1$ and $c_2$ such that
\[
\bP\bigl(\|w\|_T<\delta\bigr)\sim c_1\exp\biggl(-\frac{c_2}{\delta^2}
\biggr)\qquad \mbox{as } \delta \downarrow0.\vspace*{-9pt}
\]
\begin{longlist}[(iii)]
\item[(ii)] For all $i,j=1,\ldots,d_1$,
\[
\lim_{M\uparrow\infty}\sup_{0<\delta\leq1}\bP\bigl(\bigl\|
\kappa^{ij}\bigr\| _T>M\delta | \|w\|_T<\delta
\bigr)=0.
\]

\item[(iii)]
For all $i,j=1,\ldots,d_1$, we have
\[
\lim_{M\uparrow\infty}\sup_{0<\delta\leq1}\bP\bigl(\bigl\|
\zeta^{ij}\bigr\| _T>M\delta | \|w\|_T<\delta
\bigr)=0.
\]
\end{longlist}
\end{lemma}

In particular, we deduce from this lemma that for every $\e>0$ and
$\alpha\in(0,1)$,
%
%e2.1 #&#
\begin{equation}
\label{zalpha} \bP\bigl(\bigl\|\zeta^{ij}\bigr\|_T>\e
\delta^\alpha | \|w\|_T<\delta\bigr)\rightarrow0\qquad \mbox{as }\delta
\downarrow0. \label{est2}
\end{equation}
In fact, for arbitrary $M>0$, take $\delta_0>0$ such that $\e\delta
_0^{\alpha-1}\geq M$. Then for any $0<\delta<\delta_0$,
\[
\bP\bigl(\bigl\|\zeta^{ij}\bigr\|_T>\e\delta^\alpha | \|w
\|_T<\delta\bigr) \leq\bP\bigl(\bigl\|\zeta^{ij}
\bigr\|_T>M\delta | \|w\|_T<\delta\bigr).
\]
Thus
\begin{eqnarray*}
&&\limsup_{\delta\downarrow0}\bP\bigl(\bigl\|\zeta^{ij}\bigr\|_T>
\e\delta^\alpha | \|w\| _T<\delta\bigr)
\\
&&\qquad\leq\sup
_{0<\delta<1}\bP\bigl(\bigl\|\zeta^{ij}\bigr\|_T>M\delta |
\|w\| _T<\delta\bigr).
\end{eqnarray*}
By letting $M\uparrow\infty$ we arrive at (\ref{zalpha}) according to
(iii) in the above lemma.
In the same way, we can also obtain
%
%e2.2 #&#
\begin{equation}
\label{zalpha2} \bP\bigl(\bigl\|\kappa^{ij}\bigr\|_T>\e
\delta^\alpha | \|w\|_T<\delta \bigr)\rightarrow0 \qquad\mbox{as }\delta
\downarrow0. %\label{est2}
\end{equation}

We have the following exponential integrability result.

%pr2.1 #&#
\begin{proposition}\label{exp1}
There exists $\beta>0$ such that
\[
\bE\bigl[e^{\beta(|K|_T)^2}\bigr]<\infty,\qquad \bE\bigl[e^{\beta\|X\|_T^2}\bigr]<\infty.
\]
\end{proposition}

\begin{pf}
By It\^o's formula and (H$_2$) we have
%
%e2.3 #&#
\begin{eqnarray}\qquad
\label{itophi} \a c_0|K|_t&\leq&
\vph(X_t)-\vph(X_0)-\int_0^t(D
\vph) (X_s)\sigma (X_s)\,\d w_s-\int
_0^t D\vph(X_s)
\tilde{b}(X_s)\,\d s
\nonumber
\\[-8pt]
\\[-8pt]
\nonumber
&&{}-\frac{1}2\int_0^t\tr
\bigl[D^2\vph(X_s) \bigl(\s\s^*\bigr) (X_s)
\bigr]\,\d s.
\end{eqnarray}
Since $\vph\in\cC_b^2$, there exists a $\beta'>0$ such that
\[
\bE\biggl[\exp\biggl\{\beta'\biggl\|\int_0^\cdot(D
\vph) (X_s)\sigma(X_s)\,\d w_s\biggr\|
_T^2\biggr\} \biggr]<\infty.
\]
From this the first inequality follows immediately, and the second
follows from the first together with equation (\ref{msde1}).
\end{pf}

%le2.2 #&#
\begin{lemma}\label{estfk}
$
\lim_{\delta\downarrow0}\bP (|K|_T\geq\e\delta^{-{1}/{2}} | \|
w\|_T< \delta )= 0$.
\end{lemma}

\begin{pf}
We have by Lemma~\ref{w} and Proposition~\ref{exp1} that
%
%e2.4 #&#
\begin{eqnarray}
\label{est0} \lim_{\delta\downarrow0}\bP \bigl(|K|_T\geq\e
\delta^{-{3}/{2}} | \| w\|_T< \delta \bigr) \lesssim\lim
_{\delta\downarrow0}\frac{\exp\{-\e^2\delta
^{-3}\beta\}
}{\exp\{-C\delta^{-2}\}}= 0.
\end{eqnarray}
Next we prove that for $f\in\cC^2_b(\mR^d;\mR)$ and $1\leq k\leq d_1$ we have
%
%e2.5 #&#
\begin{equation}
\label{est1} \lim_{\delta\downarrow0} \bP\biggl(\biggl\|\int
_0^\cdot f(X_s)\circ\d
w^k_s\biggr\|_T\geq\e\delta^{-{1}/2} \Big| \| w
\|_T<\delta\biggr)=0.
\end{equation}
Set $f_i(x):=\frac{\p f}{\p x_i}(x)$. By It\^o's formula we have
\begin{eqnarray*}
\int_0^t f(X_s)\circ\d
w^k_s&=&f(X_t)w^k_t-
\int_0^t\bigl[f_i\s
_j^i\bigr](X_s)w^k_s
\circ\d w^j_s
\\
&&{}-\int_0^t\bigl[f_ib^i
\bigr](X_s)w^k\,\d s-\int_0^t
f_i(X_s)w^k_s\,\d
K_s^i
\\
&=:&I_1(t)-I_2(t)-I_3(t)+I_4(t).
\end{eqnarray*}
We need to prove
\begin{eqnarray*}
\lim_{\delta\downarrow0} \bP\bigl(\|I_i\|_T\geq\e
\delta^{-{1}/2} | \|w\|_T< \delta\bigr)=0,\qquad i=1,2,3,4.
\end{eqnarray*}
This is obvious for $I_1$ and $I_3$. To show this for $I_2$ we notice that
\begin{eqnarray*}
I_2(t)&=&\int_0^t
\bigl[f_i\s_j^i\bigr](X_s)w^k_s\,\mathrm{d}w^j_s+
\frac{1}2\int_0^t\bigl[f_i
\s ^i_j\bigr](X_s)\delta^{kj}\,\mathrm{d}s
\\
&&{}+\frac{1}2\int_0^t
\bigl[f_i\s^i_j\bigr]_q
\s_l^q(X_s)w_s^k
\delta^{lj}\,\mathrm{d}s
\\
&:=& I_{21}(t)+I_{22}(t)+I_{23}(t).
\end{eqnarray*}
Noticing that $f$ and $\sigma$ are bounded, the sets $\{\|I_{2i}\|
_T>\e
\delta^{-{1}/2}\}\cap\break \{\|w\|_T<\delta\},  i=2,3$ will be empty for
small $\delta$ and thus
\begin{eqnarray*}
\lim_{\delta\downarrow0} \bigl\{\bP\bigl(\|I_{22}
\|_T\geq\e\delta^{-{1}/2} | \|w\|_T< \delta\bigr) +
\bP\bigl(\|I_{23}\|_T\geq\e\delta^{-{1}/2} | \|w
\|_T< \delta\bigr) \bigr\}=0.
\end{eqnarray*}
Since for $t\in[0,T]$,
\[
{\langle}I_{21},I_{21}{\rangle}(t)=\sum
_{j=1}^{d_1}\int_0^t
\bigl[f_i\s _j^i\bigr]^2(X_s)
\bigl(w^k_s\bigr)^2\,\mathrm{d}s\lesssim\|w
\|_t^2.
\]
By the exponential inequality (cf. \cite{yor}, Exercise IV.3.16) we have
\[
\lim_{\delta\downarrow0}\bP\bigl(\|I_{21}\|_T\geq\e
\delta^{-{1}/2} | \| w\|_T< \delta\bigr)\lesssim \lim
_{\delta\downarrow0}\frac{\exp\{-\e^2\delta^{-3}\}}{\exp\{
-C\delta
^{-2}\}}=0.
\]
Hence
\[
\lim_{\delta\downarrow0} \bP\bigl(\|I_{2}\|_T\geq\e
\delta^{-{1}/2} | \|w\|_T\leq\delta\bigr)=0.
\]
Finally, since
\[
\|I_4\|_T\lesssim\|w\|_T|K|_T,
\]
we have by using (\ref{est0}) that
\[
\lim_{\delta\downarrow0} \bP\bigl(\|I_{4}\|_T\geq\e
\delta^{-{1}/2} | \|w\|_T< \delta\bigr)=0.
\]
Thus (\ref{est1}) has been proved. Now the result follows from (\ref
{itophi}) and (\ref{est1}).
\end{pf}

%co2.1 #&#
\begin{corollary}\label{kzc}
For every $\e>0$,
%
%e2.6 #&#
\begin{eqnarray}
\label{kz1} \lim_{\delta\downarrow0}\bP\bigl(\bigl\|\zeta^{ij}
\bigr\|_T|K|_T>\e | \|w\| _T<\delta\bigr)&=&0,\\
\label{kz2} \lim_{\delta\downarrow0}\bP \biggl(\biggl\llVert \int
_0^\cdot\zeta ^{ij}(s)\,\d
K_s\biggr\rrVert _T>\e \Big| \|w\|_T<\delta
\biggr)&\rightarrow&0.
\end{eqnarray}
\end{corollary}
\begin{pf}
It suffices to prove (\ref{kz1}). Using (\ref{zalpha}) with $\a
=\frac{1}2$ and the above
lemma we have
\begin{eqnarray*}
&&\bP\bigl(\bigl\|\zeta^{ij}\bigr\|_T|K|_T>\e | \|w
\|_T<\delta\bigr)
\\
&&\qquad\leq\bP\bigl(\bigl\|\zeta^{ij}\bigr\|_T>\delta^{{1}/{2}} |
\|w\|_T<\delta \bigr)+\bP \bigl(|K|_T>\e
\delta^{-{1}/{2}} | \|w\|_T<\delta\bigr)
\\
&&\qquad \rightarrow0,\qquad \delta\downarrow0.
\end{eqnarray*}
\upqed\end{pf}

Now we can prove the following:

%le2.3 #&#
\begin{lemma}\label{lem}
Suppose $f\in\cC_b(\mR^d;\mR)$ is uniformly continuous.
Then for all $\e>0$ and $i,j=1,2,\ldots,d_1$,
%
%e2.8 #&#
\begin{equation}
\label{es} \lim_{\delta\downarrow0}\bP \biggl(\biggl\llVert \int
_0^\cdot f(X_s)\,\d \zeta
^{ij}(s)\biggr\rrVert _T>\e \Big| \|w\|_T<\delta
\biggr)\rightarrow0.
\end{equation}
\end{lemma}
\begin{pf}
First we assume that $f\in\cC_b^2(\mR^d;\mR)$.
It\^o's formula gives us
\begin{eqnarray*}
\int_0^tf(X_s)\,\d
\zeta^{ij}(s)&=&f(X_t)\zeta^{ij}(t) -\int
_0^t\zeta^{ij}(s)f_l(X_s)
\sigma_k^l(X_s)\,\d w^k_s
\\
&&{}-\int_0^t(Lf) (X_s)
\zeta^{ij}(s)\,\d s-\int_0^tf_l(X_s)
\sigma _j^l(X_s)w^i_s
\,\d s
\\
&&{}-\int_0^tf_l(X_s)
\zeta^{ij}(s)\,\d K^l_s
\\
&=&:\sum_{q=1}^5I_{2q},
\end{eqnarray*}
where $L:=\frac{1}2\sum_{i,j}a_{ij}\partial_i\,\partial_j+\sum_i\tilde
{b}_i\partial_i$.

It is easy to see that for $q=1,3, 4$,
%
%e2.9 #&#
\begin{eqnarray}
\lim_{\delta\downarrow0}\bP\bigl(\|I_{2q}\|_T>\e | \|w
\|_T<\delta\bigr)=0.
\end{eqnarray}
Since
\[
\biggl\|\int_0^\cdot f_l(X_s)
\zeta_s^{ij}\,\d K^l_s
\biggr\|_T\lesssim\|\zeta\|_T|K|_T,
\]
we have
\[
\bP\bigl(\|I_{25}\|_T>\e | \|w\|_T<\delta\bigr)\leq P\bigl(
\|\zeta\|_T|K|_T>\e | \|w\| _T<\delta\bigr).
\]

Consequently by (\ref{kz1}),
\[
\lim_{\delta\downarrow0}\bP\bigl(\|I_{25}\|_T>\e | \|w
\|_T<\delta\bigr)=0.
\]
Now we deal with $I_{22}$.
Set $g_k(x):=-f_l(x)\sigma_k^l(x)$, $g_{k,l}:=\frac{\p}{\p x^l}g_k(x)$.
We have by It\^o's formula,
\begin{eqnarray*}
I_{22}&=&\int_0^t
g_k(X_s)\zeta^{ij}(s)\,\d w^k_s
\\
&=&g_k(X_t)\zeta^{ij}(t)w^k_t-
\int_0^tg_{k,l}(X_s)
\sigma_q^l(X_s)\zeta^{ij}(s)w^k_s
\,\d w^q_s
\\
&&{}-\int_0^t(Lg_{k})
(X_s)\zeta^{ij}(s)w^k_s\,\d s -
\int_0^tg_k(X_s)w^k_s
\,\d\zeta^{ij}(s)
\\
&&{}-\int_0^tg_j(X_s)w^i_s
\,\d s-\int_0^t\zeta^{ij}(s)g_{k,l}(X_s)
\sigma _q^l(X_s)\delta^{kq}\,\d s
\\
&&{}-\int_0^tg_{k,l}(X_s)
\sigma_j^l(X_s)w^k(s)w^i_s
\,\d s -\int_0^tg_{k,l}(X_s)
\zeta^{ij}(s)w^k_s\,\d K^l_s
\\
&:=&\sum_{i=1}^8I_{22i}.
\end{eqnarray*}
Obviously
\[
\lim_{\delta\downarrow0}\bP\bigl(\|I_{22i}\|_T>\e | \|w\|<
\delta \bigr)=0,\qquad i=1,3,4,5,7,
\]
and it is clear from Corollary~\ref{kzc} that it holds also for $i=8$.
For $I_{222}$ we notice
\[
I_{222}(t)=M_t,
\]
where
\[
M_t=\int_0^tg_{k,l}(X_s)
\sigma_q^l(X_s)\zeta^{ij}(s)w^k_s
\,\d w^q_s.
\]
It suffices to prove
%
%e2.10 #&#
\begin{equation}
\label{MT} \lim_{\delta\downarrow0}\bP\bigl(\|M\|_T>\e | \|w
\|_T<\delta\bigr)=0.
\end{equation}
Since
\[
\langle M\rangle (t)\lesssim\int_0^t\bigl\|
\zeta^{ij}\bigr\|_s^2\|w\|^2_s
\,\d s,
\]
we have by exponential inequality
\begin{eqnarray*}
&&\bP\bigl(\|M\|_T>\e,\bigl \|\zeta^{ij}\bigr\|_T<A
\delta, \|w\|_T<\delta\bigr)
\\
&&\qquad\leq \bP\bigl(\|M\|_T>\e, \langle M \rangle (T)\leq
cA^2\delta^4\bigr) \leq c\exp\bigl\{-cA^{-2}
\delta^{-4}\bigr\}\to0,\qquad \delta\downarrow0.
\end{eqnarray*}
Since
\begin{eqnarray*}
\bP\bigl(\|M\|_T>\e | \|w\|_T<\delta\bigr) &=&\bP\bigl(\|M
\|_T>\e, \bigl\|\zeta^{ij}\bigr\|_T>A\delta | \|w
\|_T<\delta\bigr)
\\
&&{}+\bP\bigl(\|M\|_T>\e, \bigl\|\zeta^{ij}\bigr\|_T\leq A
\delta | \|w\|_T<\delta\bigr)
\\
&\leq& \sup_{0\leq\delta\leq1}\bP\bigl(\|M\|_T>\e, \bigl\|
\zeta^{ij}\bigr\| _T>A\delta | \|w\|_T<\delta\bigr)
\\
&&{}+\bP\bigl(\|M\|_T>\e, \bigl\|\zeta^{ij}\bigr\|_T\leq A
\delta | \|w\|<\delta\bigr),
\end{eqnarray*}
we have
\begin{eqnarray*}
\lim_{\delta\downarrow0}\bP\bigl(\|M\|_T>\e | \|w\|_T<
\delta\bigr)\leq \sup_{0\leq\delta\leq1}\bP\bigl(\|M\|_T>\e,\bigl \|
\zeta^{ij}\bigr\|_T>A\delta | \| w\|_T<\delta\bigr).
\end{eqnarray*}
Hence by letting $A\to\infty$ we have
\[
\lim_{\delta\downarrow0}\bP\bigl(\|M\|_T>\e | \|w\|_T<
\delta\bigr)=0.
\]

Now we extend the result to $f\in\cC_b$, which is uniformly
continuous. Let $\e>0$ be given. For
any $\e'>0$ choose an $\eta\in(0,\frac{\e}{2T})$ sufficiently small
such that
\[
c_2-\frac{\e^2}{32\eta^2T}<0, \qquad\frac{4}{c_1}\exp\biggl
\{c_2-\frac
{\e
^2}{32\eta^2T}\biggr\}<\e',
\]
where $c_1$ and $c_2$ are constants appearing in Lemma~\ref{w}. Then
choose a $g\in\cC_b^2$ such that
$
\|f-g\|_T<\eta$. Note that
\begin{eqnarray*}
&&\int_0^t f(X_s)\,\d
\zeta^{ij}(s)-\int_0^t
g(X_s)\,\d\zeta^{ij}(s)
\\
&&\qquad=\int_0^t (f-g) (X_s)w^i_s
\,\d w^j_s+\frac{\delta_{ij}}{2}\int_0^t
(f-g) (X_s)\,\d s
\\
&&\qquad=:Y_1(t)+Y_2(t).
\end{eqnarray*}
It is easy to see $\|Y_2\|_T<\frac{\e}{4}$. Moreover, since
$
\langle Y_1\rangle (T)\leq\eta^2\|w\|_T^2T$,
we have by exponential inequality and with arguments similar to the
proof of (\ref{MT}) that if $\delta\in(0,1]$,
\begin{eqnarray*}
&&\bP\biggl(\|Y_1\|_T\geq\frac{\e}{4} \Big| \|w
\|_T<\delta\biggr)
\\
&&\qquad\leq\bP\biggl(\|Y_1\|_T\geq\frac{\e}{4}, \langle
Y_1\rangle (T)\leq \eta ^2\delta ^2T\biggr)P\bigl(
\|w\|_T<\delta\bigr)^{-1}
\\
&&\qquad\leq\frac{4}{c_1}\exp\biggl\{\frac{1}{\delta^2}\biggl(c_2-
\frac{\e
^2}{32\eta
^2T}\biggr)\biggr\}\leq\e'.
\end{eqnarray*}
Thus for such $\delta$,
\begin{eqnarray*}
&&\bP\biggl(\biggl\|\int_0^\cdot f(X_s)
\,\d\zeta^{ij}(s)\biggr\|_T\geq\e \Big| \|w\| _T<\delta
\biggr)
\\
&&\qquad\leq\e'+\bP\biggl(\biggl\|\int_0^\cdot
g(X_s)\,\d\zeta^{ij}(s)\biggr\|_T\geq\e /2 \Big| \|w\| <
\delta\biggr).
\end{eqnarray*}
Now we conclude by letting $\delta\to0$ and by the arbitrariness of
$\e'$.
\end{pf}

%le2.4 #&#
\begin{lemma}\label{wehave} We have:
\textup{(i)} For all $f\in\cC_b^2(\mR^d;\mR)$, $\e>0$ and $1\leq k\leq d_1$,
\[
\lim_{\delta\downarrow0}\bP\biggl(\biggl\|\int_0^\cdot
f(X_s)\circ \mathrm{d}w_s^k\biggr\| _T\geq \e
\Big| \|w\|_T<\delta\biggr)=0.
\]

\textup{(ii)} There exists a constant $c_3>0$ such that
\[
\lim_{\delta\downarrow0}\bP\bigl(|K|_T>c_3 | \|w
\|_T<\delta\bigr)=0.
\]
\end{lemma}
\begin{pf}
It suffices to prove (i), since then (ii) follows from (i) and (\ref{itophi}).

We have
\begin{eqnarray*}
\int_0^t f(X_s)\circ
\mathrm{d}w_s^k&=&f(X_t)w_t^k-
\int_0^t\bigl[f_i\s
_j^i\bigr](X_s)\,\mathrm{d}\zeta ^{kj}
\\
&&{}-\int_0^t\bigl[f_ib^i
\bigr](X_s)w_s^k\,\mathrm{d}s-\int
_0^tf_i(X_s)w_s^k\,\mathrm{d}K_s^i
\\
&:=&I_1(t)-I_2(t)-I_3(t)-I_4(t).
\end{eqnarray*}
Since
\[
\|I_4\|_T\lesssim\|w\|_T|K|_T,
\]
by Lemma~\ref{estfk} we have
\[
\lim_{\delta\downarrow0}\bP\bigl(\|I_4\|_T\geq\e | \|w
\|_T<\delta\bigr) \leq\lim_{\delta\downarrow0}\bP\bigl(|K|_T
\geq c\e\delta^{-1} | \|w\| _T<\delta\bigr)=0,
\]
while by Lemma~\ref{lem} we have
\[
\lim_{\delta\downarrow0}\bP\bigl(\|I_2\|_T\geq\e | \|w
\|_T<\delta\bigr)=0.
\]
Finally, it is trivial that
\[
\lim_{\delta\downarrow0}\bP\bigl(\|I_i\|_T\geq\e | \|w
\|_T<\delta\bigr)=0,\qquad i=1,3.
\]
This completes the proof.
\end{pf}

Now are ready to state our main result. Let $(Y,l)$ denote the
solution of the following
deterministic Skorohod problem:
%
%e2.11 #&#
\begin{equation}\qquad
\label{msdes}\cases{ %
\displaystyle Y_t=Y_0+
\int_0^t\sigma(Y_s)\,\d
h_s+\int_0^t b(Y_s)
\,\d s+l_t,&\quad $Y_0=x,$
\vspace*{2pt}\cr
\displaystyle |l|_t=\int_0^t
1_{\partial D}(Y_s)\,\d|l|_s, &\quad $\displaystyle l_t=\int
_0^t \eta(s)\,\d|l|_s,$}
\end{equation}
where $\eta(s)\in\cN_{Y_s}$.

%th2.1 #&#
\begin{theorem}\label{main}
For any $h\in\cS$ and $\e>0$,
\[
\bP\bigl(\|X-Y\|_T+\|K-l\|_T<\e|\|w-h\|_T<\delta\bigr)
\rightarrow1 \qquad \mbox{as }\delta\downarrow0.
\]
\end{theorem}

\begin{pf}
We first assume $h\equiv0$. Since $(X,K)$ and $(Y,l)$ are solutions to
equations (\ref{msde1}) and~(\ref{msdes}), respectively, we have
\[
X_t-Y_t=\int_0^t
\sigma(X_s)\circ\d w_s +\int_0^t
\bigl(b(X_s)-b(Y_s)\bigr)\,\d s+\int_0^t(
\d K_s-\d l_s).
\]
Set
\[
\Psi(x):=1-e^{-{|x|^2}/{2}};
\]
then
\begin{eqnarray*}
\Psi_i(x)&:=&\frac{\p}{\p x_i}\Psi(x)=e^{-{|x|^2}/{2}}x_i;
\\
\Psi_{i,j}(x)&:=&\frac{\p^2}{\p x_i\,\p x_j}\Psi(x)=-e^{-{|x|^2}/{2}}[x_ix_j+
\delta_{ij}];
\\
G(t)&:=&X_t-Y_t,\qquad \vph_i(x):=\frac{\partial}{\partial x_i}
\vph(x).
\end{eqnarray*}
By It\^o's formula we have
\begin{eqnarray*}
&&\exp\biggl\{\frac{2}{\a}\bigl(\vph(X_t)+
\vph(Y_t)\bigr)\biggr\}\d \biggl[\exp\biggl\{-\frac{2}{\a}\bigl(
\vph(X_t)+\vph(Y_t)\bigr)\biggr\}\Psi\bigl(G(t)\bigr)
\biggr]
\\
&&\qquad=\Psi_i\bigl(G(t)\bigr)\sigma_k^i(X_t)
\circ\d w^k_t+\Psi _i\bigl(G(t)\bigr)
\bigl(b^i(X_t)-b^i(Y_t)\bigr)\,\d
t\\
&&\qquad\quad{}+\Psi_i\bigl(G(t)\bigr) \bigl(\d K^i_t-\d
l^i_t\bigr)
\\
&&\qquad\quad{}-\frac{2}{\a}\Psi\bigl(G(t)\bigr) \bigl[\vph_i(X_t)
\sigma^i_k(X_t)\circ\d w^k_t+
\vph_i(X_t)b^i(X_t)\,\d t+
\vph_i(X_t)\,\d K^i_t
\\
&&\hspace*{167pt}\qquad\quad{}+\vph_i(Y_t)b^i(Y_t)\,\d t+
\vph_i(Y_t)\,\d l^i_t \bigr]\\
&&\qquad\quad{}-\frac{2}{\a}\Psi\bigl(G(t)\bigr)\Psi_i\bigl(G(t)\bigr)\vph_i(X_t)
\sigma^i_k\sigma^i_k(X_t)\d t.
\end{eqnarray*}
Using the elementary inequality $1-e^{-t}\geq te^{-t}$ for $t\geq0$
and conditions (H$_1$)--(H$_2$), we have
\begin{eqnarray*}
&&\Psi_i\bigl(G(t)\bigr)\,\d K^i_t-
\frac{2}{\a}\Psi\bigl(G(t)\bigr)\vph_i(X_t)\,\d
K^i_t
\\
&&\qquad= \biggl[e^{-{|X_t-Y_t|^2}/{2}}(X_t-Y_t)^*\xi_t-
\frac{2}{\a
}\bigl(1-e^{-{|X_t-Y_t|^2}/{2}}\bigr) \vph_i(X_t)
\xi^i_t \biggr]\,\d|K|_t
\\
&&\qquad\leq e^{-{|X_t-Y_t|^2}/{2}} \bigl[(X_t-Y_t)^*\xi
_t-c_0|X_t-Y_t|^2
\bigr]\,\d|K|_t\leq0,
\\
&& -\Psi_i\bigl(G(t)\bigr)\,\d l^i_t-
\frac{2}{\a}\Psi\bigl(G(t)\bigr)\vph_i(Y_t)\,\d
l^i_t\leq0.
\end{eqnarray*}
Combining these with the fact $|\Psi_i(x)x_i|\lesssim\Psi(x)$, we have
\begin{eqnarray*}
\exp\biggl\{-\frac{2}{\a}\bigl(\vph(X_t)+\vph(Y_t)
\bigr)\biggr\}\Psi\bigl(G(t)\bigr) \leq\int_0^t
\rho_k(s)\circ\,\d w_s^k+C\int
_0^t\Psi\bigl(G(s)\bigr)\,\d s,
\end{eqnarray*}
where
\begin{eqnarray*}
\rho_k(s)&:=&\exp\biggl\{-\frac{2}{\a}\bigl(
\vph(X_s)+\vph(Y_s)\bigr)\biggr\}\\
&&{}\times \biggl[
\Psi_i\bigl(G(s)\bigr)\sigma_k^i(X_s)-
\frac{2}{\a}\Psi\bigl(G(s)\bigr)\vph _i(X_s)\sigma
^i_k(X_s) \biggr].
\end{eqnarray*}
By It\^o's formula
\begin{eqnarray*}
&&\exp\biggl\{\frac{2}{\a}\bigl(\vph(X_s)+
\vph(Y_s)\bigr)\biggr\}\circ\d\rho_k(s)
\\
&& \qquad=\biggl[\Psi_i\bigl(G(s)\bigr)\s^i_{kj}(X_s)+
\s_k^i(X_s)\Psi_{ij}\bigl(G(s)
\bigr) -\frac{2}{\a} \bigl(\Psi\bigl(G(s)\bigr)\vph_i(X_s)
\s^i_{kj}(X_s)
\\
&&\hspace*{65pt}\qquad\quad{}+\Psi\bigl(G(s)\bigr)\vph_{ij}(X_s)
\s^i_{k}(X_s) +\Psi_j\bigl(G(s)
\bigr)\vph_i(X_s)\s^i_{k}(X_s)
\bigr) \biggr]
\\
&&\qquad\quad{}\times \bigl[\s_l^j(X_s)\circ\d
w^l_s+b^j(X_s)\,\d s+\d
K_s^j \bigr]
\\
&&\qquad\quad{}+\biggl[\frac{2}{\a}\Psi_j\bigl(G(s)\bigr)
\vph_i(X_s)\s^i_{k}(X_s)-
\s _k^i(X_s)\Psi _{ij}\bigl(G(s)
\bigr)\biggr] \bigl[b^j(Y_s)\,\d s+\d l_s^j
\bigr]
\\
&&\qquad\quad{}-\frac{2}{\a}\rho_k(s)\exp\biggl\{\frac{2}{\a}\bigl(
\vph(X_s)+\vph (Y_s)\bigr)\biggr\}
\\
&&\qquad\quad{}\times \bigl[\vph_j(X_s) \bigl(\s_l^j(X_s)
\circ \mathrm{d}w^l_s+b^j(X_s)\,\d s+\d
K_s^j\bigr)\\
&&\hspace*{125pt}{} +\vph_j(Y_s)
\bigl(b^j(Y_s)\,\d s+\d l_s^j
\bigr) \bigr],
\end{eqnarray*}
where $\sigma_{kl}^i(x):=\frac{\p}{\p x_l}\sigma_k^i(x)$. Rearranging,
we write
\begin{eqnarray*}
\d\rho_k(s)&=&F_{kl}(X_s,Y_s)
\circ\d w_s^l+G_{kj}(X_s,Y_s)b^j(X_s)
\,\d s+G_{kj}(X_s,Y_s)\,\d K_s^j
\\
&&{}+H_{kj}(X_s,Y_s)b^j(Y_s)
\,\d s+H_{kj}(X_s,Y_s)\,\d l_s^j,
\end{eqnarray*}
where
\begin{eqnarray*}
G_{kj}(x,y):&=&\exp\biggl\{-\frac{2}{\a}\bigl(\vph(x)+\vph(y)
\bigr)\biggr\}\\
&&{}\times \biggl[\Psi _i(x-y)\s ^i_{kj}(x)+
\s_k^i(x)\Psi_{ij}(x-y)
\\
&&\hspace*{18pt}{}-\frac{2}{\a} \bigl(\Psi(x-y)\vph_i(x)\s^i_{kj}(x)+
\Psi(x-y)\vph _{ij}(x)\s^i_{k}(x)\\
&&\hspace*{140pt}{} +
\Psi_j(x-y)\vph_i(x)\s^i_{k}(x)
\bigr) \biggr]
\\
&&{}-\frac{2}{\a}\exp\biggl\{-\frac{2}{\a}\bigl(\vph(x)
+\vph(y)\bigr)
\biggr\} \\
&&{}\times \biggl[\Psi_i(x-y)-\frac{2}{\a}\Psi(x-y)
\vph_i(x)\biggr]\s_k^i(x)
\vph_j(x),
\\
F_{kl}(x,y):&=&G_{kj}(x,y)\s_l^j(x),
\\
H_{kj}(x,y):&=&\biggl[\frac{2}{\a}\Psi_j(x-y)
\vph_i(x)\s^i_{k}(x)-\s _k^i(x)
\Psi_{ij}(x-y)\biggr]\
\\
&&{}\times\exp\biggl\{-\frac{2}{\a}\bigl(\vph(x)+
\vph(y)\bigr)\biggr\}
\\
&&{}-\frac{2}{\a}\exp\biggl\{-\frac{2}{\a}\bigl(\vph(x)+\vph(y)\bigr)
\biggr\}\\
&&{}\times \biggl[\Psi_i(x-y)-\frac{2}{\a}\Psi(x-y)
\vph_i(x)\biggr]\s_k^i(x)
\vph_j(y).
\end{eqnarray*}
Thus we have by It\^o's formula,
\begin{eqnarray*}
\int_0^t\rho_k(s)\circ\d
w_s^k&=&\rho_k(t)w^k_t-
\int_0^t F_{kl}(X_s,Y_s)w_s^k
\circ\d w^l_s
\\
&&{}-\int_0^t G_{kj}(X_s,Y_s)b^j(X_s)w_s^k
\,\d s -\int_0^t G_{kj}(X_s,Y_s)w_s^k
\,\d K^j_s
\\
&&{}-\int_0^t H_{kj}(X_s,Y_s)b^j(Y_s)w^k_s
\,\d s-\int_0^t H_{kj}(X_s,Y_s)w_s^k
\,\d l^j_s
\\
&=:&I_1(t)-I_2(t)-I_3(t)-I_4(t)-I_5(t)-I_6(t).
\end{eqnarray*}
Obviously,
\begin{eqnarray*}
\sum_{i\neq2}\|I_i\|_T
\lesssim\bigl(1+|K|_T\bigr)\|w\|_T.
\end{eqnarray*}
Thus
\begin{eqnarray*}
\lim_{\delta\downarrow0}\bP\biggl(\sum_{i\neq2}
\|I_i\|_T\geq\e \Big| \|w\| _T<\delta\biggr) \leq
\lim_{\delta\downarrow0}\bP \bigl(\bigl(1+|K|_T\bigr)\|w\|_T
\gtrsim\e | \|w\| _T<\delta \bigr)=0.
\end{eqnarray*}
As for $I_2$ we have
\begin{eqnarray*}
I_2(t)&=&\int_0^t
F_{kl}(X_s,Y_s)\,\d\zeta^{kl}_s+
\frac{1}2\int_0^t \frac{\partial}{\partial x_j}F_{kl}(X_s,Y_s)
\s_p^j(X_s)w_s^k
\delta^{p
l}\,\d s
\\
&=:&I_{21}+I_{22}.
\end{eqnarray*}
It is easily seen that
\[
\lim_{\delta\downarrow0}\bP\bigl(\|I_{22}\|_T\geq\e | \|w
\|_T<\delta\bigr)=0,
\]
and by applying Lemma~\ref{lem} to the functions $F_{kl}$ (in place of
$f$ there) and
the system satisfied by $(X,Y)$ (in place of $X$ there), we have that
\[
\lim_{\delta\downarrow0}\bP\bigl(\|I_{21}\|_T\geq\e | \|w
\|_T<\delta\bigr)=0.
\]
Consequently
\[
\lim_{\delta\downarrow0}\bP\bigl(\|I_{2}\|_T\geq\e | \|w
\|_T<\delta\bigr)=0.
\]
Combining all the above and the fact that $\vph$ is bounded, we have
\[
\Psi\bigl(G(t)\bigr) \leq C\int_0^t\Psi
\bigl(G(s)\bigr)\,\d s+A(t),
\]
where $A(t)$ satisfies that for every $\e>0$,
\[
\lim_{\delta\downarrow0}\bP\bigl(\|A\|_T>\e |\|w\|_T<
\delta\bigr)=0.
\]
On the set $\{\omega;\|A\|_T <\e\}$, we have
\[
\Psi\bigl(G(t)\bigr)\leq\e e^{C}\leq C \e,
\]
that is,
\[
\|X-Y\|_T\leq\sqrt{-2\ln(1-C\e)}.
\]
Since $\e$ is arbitrarily small,
\[
\bP\bigl(\|X-Y\|_T>\e | \|w\|_T<\delta\bigr)\rightarrow0 \qquad\mbox{as } \delta
\downarrow0.
\]
Finally, to see
\[
\bP\bigl(\|K-l\|_T<\e | \|w\|_T<\delta\bigr)\rightarrow1\qquad \mbox{as } \delta
\downarrow0,
\]
it suffices to notice that
\[
K_t-l_t=X_t-Y_t-\int
_0^t\sigma(X_s)\circ\d
w_s-\int_0^t\bigl(b(X_s)-b(Y_s)
\bigr)\,\d s
\]
and use Lemma~\ref{wehave}.

For general $h\in\cS$, just as in the proof of \cite{iw}, Theorem~8.2,
pages~527--528, we set
\begin{eqnarray*}
M_1(w):=\exp\biggl\{\int_0^T
\dot{h}_s\,\d w_s-\frac{1}2\int
_0^T|\dot {h}_s|^2\,\d s
\biggr\} ,\qquad \d\bP'=M_1\,\d\bP.
\end{eqnarray*}
Then $w'_t:=w_t-h_t$ is a Brownian motion under $\bP'$, and $(X,K)$,
$(Y,l)$ satisfy the following equations, respectively:
\begin{eqnarray*}
X_t&=&x+\int_0^tb'(s,X_s)
\,\d s+\int_0^t\sigma(X_t)\circ\d
w_s+K_t,
\\
Y_t&=&x+\int_0^tb'(s,Y_s)
\,\d s+l_t,
\end{eqnarray*}
where $b'(s,x):=b(x)+\sigma(x)\dot{h}_s$.

Therefore according to the case of $h\equiv0$ we have for every $\e>0$,
\begin{eqnarray*}
\bP'\bigl(\|X-Y\|_T>\e |\bigl \|w'
\bigr\|_T<\delta\bigr)&\rightarrow&0 \qquad\mbox{as } \delta \downarrow0,
\\
\bP'\bigl(\|K-l\|_T>\e | \bigl\|w'
\bigr\|_T<\delta\bigr)&\rightarrow&0\qquad \mbox{as } \delta\downarrow0,
\end{eqnarray*}
which, together with the fact that $M_1$ is a continuous functional of
$w$, yields that
\begin{eqnarray*}
&&\lim_{\delta\downarrow0}\bP\bigl(\|X-Y\|_T<\e |
\|w'\|_T<\delta\bigr)
\\
&&\qquad=\lim_{\delta\downarrow0}\frac{\bP(\|X-Y\|_T<\e,\|w'\|_T<\delta
)}{\bE
(M_1\mathbh{1}_{\{\|X-Y\|_T<\e,\|w-h\|_T<\delta\}})} \times\frac{\bE(M_1\mathbh{1}_{\{\|w-h\|_T<\delta\}})}{\bP(\|w-h\|
_T<\delta)}
\rightarrow1.
\end{eqnarray*}
\upqed\end{pf}

%re2.1 #&#
\begin{remark}
In the last step of the proof above, we encounter the situation that
the drift $b'$ depends also on time $t$. But as in \cite{iw}, Theorem~8.2, everything still works with trivial modifications.
\end{remark}

%s3 #&#
\section{The support problem}\label{sec3}
%s3.1 #&#
\subsection{Conditions and useful estimates}\label{sec3.1}
The approximate continuity theorem proved in the above section together
with the Wong--Zakai
approximation theorem proved in \cite{aida} gives, in a similar way
paved in \cite{str},
the support theorem for reflected diffusions under the conditions
(H$_1$)--(H$_2$). In this section we will prove the support theorem
based upon the idea in \cite{millet} when the domain $D$ is supposed to
satisfy the following conditions:
\begin{longlist}[(A)]
\item[(A)] There exists a constant $r_0>0$ such that for any
$x\in
\partial D$,
\begin{eqnarray*}
\mathcal{N}_x=\mathcal{N}_{x,r_0}\neq\varnothing.
\end{eqnarray*}
\item[(B)] There exist constants $\delta>0$ and $\beta\geq1$
satisfying that for any $x\in\partial D$, there exists a unit vector
$l_x$ such that
\begin{eqnarray*}
{\langle}l_x,n{\rangle}\geq1/\beta\qquad \mbox{for any } n\in \bigcup
_{y\in
B(x,\delta)\cap\partial D}\mathcal{N}_y,
\end{eqnarray*}
where ${\langle}\cdot,\cdot{\rangle}$ denotes the usual inner
product in $\mR^d$.
\item[(C)] There exists a function $\vph\in\mathcal
{C}_b^2(\mR^d)$
and a positive constant $\gamma$ such that for any $x\in\partial D$,
$y\in\bar{D}$ and $n\in\mathcal{N}_x$,
\begin{eqnarray*}
{\langle}y-x,n{\rangle}+\frac{1}{\gamma}{\bigl\langle}D\vph(x), n{\bigr
\rangle}|y-x|^2\geq0.
\end{eqnarray*}
\item[(D)] There exist $m\geq1$, $\lambda>0$, $R>0$,
$a_1,\ldots,
a_m\in\mS^{d-1}$ and $x_1,\ldots, x_m\in\partial D$ such that
$\partial D\subset\bigcup_{i=1}^m B(x_i,R)$ and $x\in\partial D\cap
B(x_i,2R)\Rightarrow n\cdot a_i\geq\lambda$, $\forall n\in\mathcal{N}_x$.
\end{longlist}

We will need some results from \cite{aida}.

%le3.1 #&#
\begin{lemma}[(\cite{aida}, Lemma~2.3)]\label{fv}
Assume \textup{(A)--(B)} hold, and $(x,k)$ is the solution to the
Skorohod problem associated with a continuous function $w$ such that
$x_0=w_0\in\bar{D}$. Then for $\theta\in(0,1]$, there exist constants
$c_1,  c_2,  C$ dependent on $\theta, \delta, \beta, \gamma_0$ such
that for all $0\leq s\leq t\leq T$,
\[
|k|_t^s\leq C\bigl(1+\|w\|_{[s,t],\theta}^{c_1}(t-s)
\bigr)e^{c_2\|w\|_{[s,t]}}\| w\|_{[s,t]},
\]
where (and throughout) $|k|_t^s$ denotes the total variation of $k$ on
$[s,t]$ and
\[
\|w\|_{[s,t],\theta}:=\sup_{u,v\in[s,t]}\frac
{|w_u-w_v|}{|u-v|^\theta
}, \qquad\|w
\|_{[s,t]}:=\sup_{u,v\in[s,t]}|w_u-w_v|.
\]
\end{lemma}

%le3.2 #&#
\begin{lemma}[(\cite{aida}, Lemma~2.4)]\label{bdde1}
Assume \textup{(A)} holds, and $(x,k)$ is the solution to the Skorohod
problem associated with a function $w$ having continuous bounded
variation path.
Then
\begin{eqnarray*}
|x|_t^s\leq2(\sqrt{2}+1)|w|_t^s.
\end{eqnarray*}
\end{lemma}

%le3.3 #&#
\begin{lemma}[(\cite{aida}, Lemma~2.8)]\label{solution}
Assume $D$ satisfies conditions \textup{(A)--(B)}, and $b, \sigma$
are bounded, Lipschitz continuous functions. Then there exists a unique
solution $(X,K)$ to equation (\ref{msde1}). Moreover, for all $0\leq
s<t<\infty$,
\begin{eqnarray*}
\label{finitemoment} \bE \bigl(\|X\|_{[s,t]} \bigr)^{2p}\leq
C_p|t-s|^p, \qquad\bE\bigl(|K|_t^s
\bigr)^{2p}\leq C_p|t-s|^p.
\end{eqnarray*}
\end{lemma}

Let $n\in\mN$ and $t_i=iT2^{-n}$ (here we should have used $t_i^n$
instead of $t_i$ to indicate the dependence on $n$,
but in order to not surcharge the notation, we omit the superscript
$n$), $\Delta=2^{-n}T$, and for $t\in[t_{i}, t_{i+1})$ set
\begin{eqnarray*}
\bar{t}_n&:=&t_{i-1}\vee0,\qquad \hat{t}_n:=t_i,\qquad
\Delta w_i:=w_{t_i}-w_{t_{i-1}\vee0},
\\
w^n_t&:=&w_{\bar{t}_n}+\frac{w_{\hat{t}_n}-w_{\bar{t}_n}}{\Delta
}(t-
\hat{t}_n).
\end{eqnarray*}
Consider the following reflected equation:
\[
X^n(t)=x+\int_0^t b
\bigl(X^n(s)\bigr)\,\d s+\int_0^t
\sigma\bigl(X^n(s)\bigr)\,\d w^n_s+K^n(t).
\]
Denote the solution by $(X^n,K^n)$.

%s3.2 #&#
\subsection{Support theorem}\label{sec3.2}
We first state our main theorem.

%th3.1 #&#
\begin{theorem}\label{support}
Suppose conditions \textup{(A)--(D)} hold and $\sigma\in\cC_b^2,
 b\in\cC_b^1$. Then for the solution $X$ to equation (\ref{msde1})
we have
\begin{eqnarray*}
\mbox{the support of } \bigl(\bP\circ X^{-1}\bigr) \mbox{ in }
\mW_\a =\overline {\sS(\cH)}^\a\qquad \forall\a\in\bigl[0,
\tfrac{1}2\bigr).
\end{eqnarray*}
\end{theorem}

To prove the theorem, we will apply the following results; cf. \cite{millet}.

%pr3.1 #&#
\begin{proposition}\label{milet}
Let $F$ be a measurable map from $\Omega$ to a Banach space $(\mX,\|
\cdot\|)$:
\begin{longlist}[(1)]
\item[(1)] Let $Z^\mX_1\dvtx \cH\rightarrow\mX$ be measurable and
$H_n\dvtx \Omega
\rightarrow\cH$ be a sequence of random variables such that for any
$\varepsilon>0$,
\[
\lim_n\bP \bigl(\bigl\|Z^\mX_1
\bigl(H_n(\omega)\bigr)-F(\omega)\bigr\|>\varepsilon \bigr)=0.
\]
Then $\operatorname{supp}(\bP\circ F^{-1})\subset\overline{Z^\mX_1(\cH)}$.
\item[(2)] Let $Z^\mX_2\dvtx \cH\rightarrow\mX$ be measurable and for fixed
$h$, $T_n^h\dvtx \Omega\rightarrow\Omega$ be a sequence of measurable
transformations such that $\bP\circ(T_n^h)^{-1}\ll\bP$, and for any
$\varepsilon>0$,
\[
\limsup_n\bP \bigl(\bigl\|F\bigl(T_n^h(
\omega)\bigr)-Z^\mX_2(h)\bigr\|<\varepsilon \bigr)>0.
\]
Then $\operatorname{supp}(\bP\circ F^{-1})\supset Z^\mX_2(\cH)$.
\end{longlist}
\end{proposition}

%pr3.2 #&#
\begin{proposition}\label{millet2}
Suppose $\{X^n_t\}$ is a sequence of finite dimensional processes
satisfying that for every $p\geq1$ and $s, t\in[0,T]$, there exists a
constant $C>0$,
%
%e3.1 #&#
\begin{equation}
\label{A1} \sup_n\bE\bigl|X^n_t-X^n_s\bigr|^{2p}
\leq C|t-s|^p.
\end{equation}
Then for any $\varepsilon>0$ and $\theta<\frac{1}{2}-\frac{1}{2p}$,
there exists a constant $C>0$ such that
\[
\sup_n\bP \bigl(\bigl\|X^n\bigr\|_{T,\theta}>
\varepsilon \bigr)\leq C\varepsilon^{-2p}.
\]
Moreover, besides (\ref{A1}), if for any $\varepsilon>0$,
\[
\lim_n\bP\Bigl(\sup_{1\leq i\leq2^n}\bigl|X_{t_i}^n\bigr|>
\varepsilon\Bigr)=0
\]
holds as well, then for any $\theta\in[0,1/2)$,
\[
\lim_n \bP\bigl(\bigl\|X^n\bigr\|_{T,\theta}>
\varepsilon\bigr)=0,
\]
where $\|\cdot\|_{T,\theta}$ is defined in the \hyperref[sec1]{Introduction}.
\end{proposition}

Following the idea in \cite{millet}, take
\[
Z^\mX_1=Z^\mX_2=Z(\cdot),\qquad
H_n(\omega)=w^n(\omega), \qquad T_n^h(
\omega)=w-w^n+h.
\]
Then by Girsanov's theorem, $\bP\circ(T_n^h)^{-1}\ll\bP$.

To prove Theorem~\ref{support}, by Proposition~\ref{milet}, it suffices
to prove that for every $\varepsilon>0$,
%
%e3.2 #&#
\begin{equation}
\label{wzk} \lim_n \bP \bigl(\bigl\|X-X^n
\bigr\|_{T,\theta}>\varepsilon \bigr)=0
\end{equation}
and
%
%e3.3 #&#
\begin{equation}
\label{convergence2} \lim_n\bP \bigl(\bigl\|X\bigl(w-w^n+h
\bigr)-Z(h)\bigr\|_{T,\theta}>\varepsilon \bigr)=0,
\end{equation}
where $Z(h)$ solves the following deterministic Skorohod problem:
\[
Z(h)_t=x+\int_0^t\sigma
\bigl(Z(h)_s\bigr)\dot{h}_s\,\d s+\int
_0^tb\bigl(Z(h)_s\bigr)\,\d s+
\psi_t.
\]
In what follows we will use $Z$ instead of $Z(h)$ if no confusion is possible.
(\ref{wzk}) is proved in \cite{zhangtusheng}, so we only need to prove
(\ref{convergence2}).

Using the Riemannian sum approximation of stochastic integrals, it is
easy to see that $Y^n:=X(w-w^n+h)$ solves the following RSDE:
\begin{eqnarray*}
Y^n_t=x+\int_0^t
\sigma\bigl(Y^n_s\bigr)\,\d w_s-\int
_0^t\sigma\bigl(Y^n_s
\bigr) \dot{w}^n_s \,\d s+\int_0^t
\sigma\bigl(Y^n_s\bigr) \dot{h}_s \,\d s+\int
_0^t\tilde{b}\bigl(Y^n_s
\bigr) \,\d s+\phi_t^n,
\end{eqnarray*}
where $\tilde{b}:=b+\frac{1}{2} (\nabla\sigma)\sigma$ and $\phi
^n(w)=K(w-w^n+h)$.

We first prepare some auxiliary results.

%le3.4 #&#
\begin{lemma}\label{holderZ}
For $0\leq s\leq t\leq T$,  $|Z_t-Z_s|^{2p}\leq C_p|t-s|^p$.
\end{lemma}

\begin{pf}
By Lemma~\ref{bdde1},
\begin{eqnarray*}
|Z_t-Z_s|^{2p}&\leq& \bigl[2(\sqrt{2}+1)
\bigr]^{2p}\biggl(\int_s^t\bigl|
\sigma(Z_u) \dot {h}_u +\tilde{b}(Z_u)\bigr|
\,\d u\biggr)^{2p}
\\
&\leq& C_p |t-s|^p.
\end{eqnarray*}
\upqed\end{pf}

%pr3.3 #&#
\begin{proposition}\label{holderC}
Let $p\geq1$.
Then there exists a constant $C_p>0$ independent of $n$ such that for
all $0\leq s\leq t\leq T$,
%
%e3.4 #&#
\begin{equation}
\label{3.8.2} \bE\bigl(\bigl|Y^n_t-Y^n_s\bigr|^{4p}
\bigr)\leq C_p|t-s|^p, \qquad\bE\bigl(\bigl|\phi^n_t-
\phi ^n_s\bigr|^{4p}\bigr)\leq C_p|t-s|^p.
\end{equation}
Moreover, for all $0\leq s\leq t\leq T$ and for any $\theta\in
(0,\frac{1}4)$,
%
%e3.5 #&#
\begin{equation}
\label{3.8.2theta} \bE\bigl(\bigl\|Y^n\bigr\|_{[s,t],\theta}^p\bigr)
\leq C_{p,\theta}, \qquad \bE\bigl(\bigl\|\phi^n\bigr\| _{[s,t],\theta}^p
\bigr)\leq C_{p,\theta}.
\end{equation}
\end{proposition}

To prove this proposition, we need some lemmas, and without loss of
generality we take $T=1$.

%le3.5 #&#
\begin{lemma}\label{exBM}
Let $\lambda, t>0$. Then there exists a constant $C>0$ independent of
$\lambda$ and $t$ such that
\[
\bE \bigl(e^{\lambda\|w\|_t} \bigr)\leq(1+C\lambda\sqrt{t})^{d_1}
Ce^{\lambda^2d_1t/2}.
\]
\end{lemma}

\begin{pf}
Set $\xi=\max_{0\leq s\leq t}|w_s|$. Note that $\bP(|w^i_t|\in\d
x)=\sqrt{\frac{2}{\pi t}}e^{-{x^2}/{(2t)}}\,\d x$, $i=1,\ldots, d_1$
and thus
\begin{eqnarray*}
\bE\bigl(e^{\lambda\xi}\bigr)&=&\int_0^\infty
e^s\bP(\lambda\xi>s)\,\d s+1 \leq2\int_0^\infty
e^s\bP\bigl(\lambda|w_t|>s\bigr)\,\d s+1
\\
&=&2\bE\bigl(e^{\lambda|w_t|}\bigr)-1 \leq2\prod_{i=1}^{d_1}
\sqrt{\frac{2}{\pi t}}\int_0^\infty
e^{\lambda
x}e^{-{x^2}/{(2t)}}\,\d x
\\
&\leq&(1+C\lambda\sqrt{t})^{d_1} Ce^{\lambda^2d_1t/2}.
\end{eqnarray*}
\upqed\end{pf}

%le3.6 #&#
\begin{lemma}\label{estimatemart}
Let $M_t:=\int_0^t f_s\,\d w_s$ and $|f_s|\leq c$ for some constant $c$.
Then there exists a constant $C>0$ such that for any integer $m$,
\[
\bE\|M\|_{[s,t]}^m \leq C^m
(m/2)^{{m}/2}(t-s)^{{m}/2}.
\]
\end{lemma}

\begin{pf}
It suffices to prove the result for $s=0$. Then $M_t=B_{{\langle
}M{\rangle}_t}$ where
$B$ is the DDS-Brownian motion of $M$. Note that
\[
{\langle}M{\rangle}_t\leq c^2t.
\]
The result follows from Doob's maximal inequality and that
%
%e3.6 #&#
\begin{equation}
\label{exBM2} \bE\bigl(|B_t|^{2m}\bigr)
\leq(2d_1)^mm^mt^m.
\end{equation}
\upqed\end{pf}

Set
\begin{eqnarray*}
L_t^n:=x+\int_0^t
\sigma\bigl(Y^n_s\bigr)\,\d w_s-\int
_0^t\sigma\bigl(Y^n_s
\bigr) \dot{w}^n_s \,\d s+\int_0^t
\sigma\bigl(Y^n_s\bigr) \dot{h}_s \,\d s+\int
_0^t\tilde{b}\bigl(Y^n_s
\bigr) \,\d s.
\end{eqnarray*}

%le3.7 #&#
\begin{lemma}\label{estimatevarphi}
There exists a constant $C_p$ such that for any $t\in[0,1]$ and any
$p\geq1$,
\begin{eqnarray*}
\bE\bigl(\bigl\|Y^n\bigr\|_{[\bar{t}_n,t]}^{2p}\bigr)\leq
C_p\Delta^p, \qquad\bE\bigl(\bigl|\phi ^n\bigr|^{\bar{t}_n}_{t}
\bigr)^{2p}\leq C_p\Delta^p.
\end{eqnarray*}
\end{lemma}

\begin{pf}
By Lemma~\ref{fv}, for any $\theta\in(0,1]$,
\begin{eqnarray*}
\bigl|\phi^n\bigr|^{\bar{t}_n}_{t}\leq C\bigl(1+
\bigl\|L^n\bigr\|_{[\bar{t}_n,t],\theta
}^{c_1}(t-\bar{t}_n)
\bigr)e^{c_2\|L^n\|_{[\bar{t}_n,t]}}\bigl\|L^n\bigr\|_{[\bar{t}_n,t]}.
\end{eqnarray*}

Note that for any $p\geq1$,
\begin{eqnarray*}
\bE\bigl(\bigl\|L^n\bigr\|_{[\bar{t}_n,t]}^{2p}\bigr) &\leq&
C_p \biggl[\bE\biggl(\int_{\bar{t}_n}^t\bigl\|
\sigma\bigl(Y^n_r\bigr)\bigr\|^2\,\d r
\biggr)^p+\bE\biggl(\int_{\bar
{t}_n}^t\bigl\|
\sigma\bigl(Y^n_r\bigr)\bigr\|\bigl|\dot{w}_r^n\bigr|
\,\d r\biggr)^{2p}
\\
&&{}+\bE\biggl(\int_{\bar{t}_n}^t\bigl\|\sigma
\bigl(Y^n_r\bigr)\bigr\||\dot{h}_r|\,\d r
\biggr)^{2p}+(t-\bar {t}_n)^{2p} \biggr]\leq
C_p(t-\bar{t}_n)^{p}.
\end{eqnarray*}
For any $c$, by Lemmas \ref{exBM} and~\ref{estimatemart},
\begin{eqnarray*}
&& \bE\bigl(e^{cp\|L^n\|_{[\bar{t}_n,t]}}\bigr)
\\
&&\qquad\leq\bE \bigl(e^{cp\max_{u,v\in[\bar{t}_n,t]}|\int_{u}^{v}\sigma
(Y^n_r)\,\d w_r+cp\int_{u}^{v}\sigma(Y^n_r)\dot{w}^n_r\,\d r|}
\\
&&\hspace*{42pt}\qquad\quad{}\times e^{cp|\int_{\bar{t}_n}^{t}\|\sigma(Y^n_r)\||\dot{h}_r|\,\d
r+cp\int
_{\bar{t}_n}^{t}\tilde{b}(Y^n_r)\,\d r|} \bigr)
\\
&&\qquad\leq \bigl(1+Cp\Delta^{1/2}\bigr)^{d_1}e^{Cd_1p^2\Delta+Cp\int_{\bar
{t}_n}^{t}(1+|\dot{h}_r|)\,\d r}\leq
C_p<\infty.
\end{eqnarray*}
Now combining these two estimates gives
\[
\bE\bigl(\bigl|\phi^n\bigr|^{\bar{t}_n}_{t}
\bigr)^{2p}\leq C_p\Delta^p.
\]
The other result follows from $Y^n_t=L^n_t+\phi^n_t$ and the above estimate.
\end{pf}

%le3.8 #&#
\begin{lemma}\label{sigmabar} For any $s, t\in[0,1]$,
\begin{eqnarray*}
\bE\sup_{u,v\in[s,t]}\biggl|\int_u^v
\sigma\bigl(Y^n_r\bigr)\,\d w^n_r\biggr|^{2p}
\leq C_p|t-s|^{p},\qquad \bE\bigl(\bigl\|L^n
\bigr\|_{[s,t]}^{2p}\bigr)\leq C_p|t-s|^{p}.
\end{eqnarray*}
\end{lemma}

\begin{pf} When $t_{i-1}\leq s\leq t_i\leq t\leq t_{i+1}$ for some
$1\leq i\leq2^n$, the result is trivial.
For general $s, t$, choose $1\leq l<m-1<m\leq2^n$ such that
$t_{l-1}\leq s\leq t_l<t_{m-1}\leq t\leq t_m$. Note that
\[
\int_s^t\sigma\bigl(Y^n_r
\bigr)\,\d w^n_r=\int_s^t
\bigl(\sigma\bigl(Y^n_r\bigr)-\sigma
\bigl(Y^n_{\bar
{r}_n}\bigr)\bigr)\,\d w^n_r+
\int_s^t\sigma\bigl(Y^n_{\bar{r}_n}
\bigr)\,\d w^n_r
\]
and
\begin{eqnarray*}
&&\int_s^t\sigma\bigl(Y^n_{\bar{r}_n}
\bigr)\,\d w^n_r
\\
&&\qquad=\int_s^{t_l}\sigma\bigl(Y^n_{\bar{r}_n}
\bigr)\,\d w^n_r+\sum_{j=l+1}^{m-1}
\int_{t_{j-1}}^{t_j}\sigma\bigl(Y^n_{\bar{r}_n}
\bigr)\,\d w^n_r+\int_{t_{m-1}}^{t}
\sigma \bigl(Y^n_{\bar{r}_n}\bigr)\,\d w^n_r
\\
&&\qquad=\sigma\bigl(Y^n_{t_{l-2}\vee0}\bigr)\frac{w_{t_{l-1}}-w_{t_{l-2}\vee
0}}{\Delta
}(t_l-s)
\\
&&\qquad\quad{}+\sum_{j=l+1}^{m-1}\sigma
\bigl(Y^n_{t_{j-2}\vee
0}\bigr) (w_{t_{j-1}}-w_{t_{j-2}\vee0})
\\
&&\qquad\quad{}+\sigma\bigl(Y^n_{t_{m-2}\vee0}\bigr)\frac{w_{t_{m-1}}-w_{t_{m-2}\vee
0}}{\Delta
}(t-t_{m-1}),
\end{eqnarray*}
$\int_0^\cdot\sigma(Y^n_{\bar{r}_n})\,\d w^n_r$ is the piecewise linear
interpolation of
\[
M^n_\cdot:=\int_0^{\cdot-\Delta}\sigma\bigl(Y^n\bigl(\pi
_n(r)\bigr)\bigr)\,\d w_r
\]
with
$\pi_n(r):=\max\{t_k;t_k\leq r\}$,
at $\{t_k\}_{k=0,1,\ldots, 2^n-1}$. Thus
\begin{eqnarray*}
\sup_{u,v\in[s,t]}\biggl|\int_u^v
\sigma\bigl(Y^n_{\bar{r}_n}\bigr)\,\d w^n_r\biggr|&
\leq& \sup_{l-2\leq k,k'\leq m-1}\bigl|M^n_{t_k}-M^n_{t_{k'}}\bigr|
\\
&\leq&2\sup_{t_{l-2}\leq r\leq t_{m-1}}\bigl|M^n_r-M^n_{t_l}\bigr|.
\end{eqnarray*}
Using Doob's inequality we get
\begin{eqnarray*}
\bE\sup_{u,v\in[s,t]}\biggl|\int_u^v
\sigma\bigl(Y^n_{\bar{r}_n}\bigr)\,\d w^n_r\biggr|^{2p}&
\leq& C_p\bE\sup_{t_{l-2}\leq r\leq t_{m-1}}\bigl|M^n_r-M^n_{t_l}\bigr|^{2p}
\\
&\leq&C_p\bE \biggl(\int_{t_{l-2}}^{t_{m-1}}\bigl\|
\sigma\bigl(Y^n_{\bar
{r}_n}\bigr)\bigr\| ^2\,\d r
\biggr)^p
\\
&\leq&C_p|t_{m-1}-t_{l-2}|^p
\\
&\leq&C_p|t-s|^p.
\end{eqnarray*}
By H\"older's inequalities and Lemma~\ref{estimatevarphi},
\begin{eqnarray*}
&&\bE\sup_{u,v\in[s,t]}\biggl|\int_u^v
\bigl(\sigma\bigl(Y^n_r\bigr)-\sigma
\bigl(Y^n_{\bar
{r}_n}\bigr)\bigr)\,\d w^n_r\biggr|^{2p}
\\
&&\qquad\leq\bE\int_s^t\bigl\|\sigma
\bigl(Y^n_r\bigr)-\sigma\bigl(Y^n_{\bar{r}_n}
\bigr)\bigr\| ^{2p}\bigl|\dot {w}^n_r\bigr|^{2p}\,\d
r(t-s)^{2p-1}
\\
&&\qquad\leq(t-s)^{2p-1}\int_s^t \bigl(\bE
\bigl\|\sigma\bigl(Y^n_r\bigr)-\sigma \bigl(Y^n_{\bar
{r}_n}
\bigr)\bigr\|^{4p} \bigr)^{1/2} \bigl(\bE\bigl|\dot{w}^n_r\bigr|^{4p}
\bigr)^{1/2}\,\d r
\\
&&\qquad\leq C_p(t-s)^{2p}.
\end{eqnarray*}

Now note that
\begin{eqnarray*}
L_t^n-L^n_s=\int
_s^t\sigma\bigl(Y^n_r
\bigr) \bigl(\d w_r-\d w^n_r\bigr)+\int
_s^t \sigma \bigl(Y^n_r
\bigr)\dot{h}_r\,\d r+\int_s^t
\tilde{b}\bigl(Y_r^n\bigr)\,\d r.
\end{eqnarray*}
Trivially by the Burkholder and H\"older inequalities we have
\begin{eqnarray*}
\bE\sup_{u,v\in[s,t]}\biggl|\int_u^v
\sigma\bigl(Y^n_r\bigr)\,\d w_r+\int
_u^v \sigma \bigl(Y^n_r
\bigr)\dot{h}_r\,\d r+\int_u^v
\tilde{b}\bigl(Y_r^n\bigr)\,\d r\biggr|^{2p}\leq
C|t-s|^p.
\end{eqnarray*}
From the estimates above we deduce
%
%e3.7 #&#
\begin{equation}
\label{ln} \bE\sup_{u,v\in[s,t]}\bigl|L^n_u-L^n_v\bigr|^{2p}
\leq C|t-s|^p.
\end{equation}
\upqed\end{pf}

Now we are ready to prove Proposition~\ref{holderC}.

\begin{pf*}{Proof of Proposition~\ref{holderC}}
%
%\begin{pf}
For cases of $s, t\in[t_{i-1},t_i]$ and $t_{i-1}\leq s\leq t_i<t\leq
t_{i+1}$ for some $1\leq i\leq2^n$, it follows from Lemmas \ref
{estimatevarphi}--\ref{sigmabar} that
%
%e3.8 #&#
\begin{equation}
\label{holderestimatevarphi} \bE\bigl\|Y^n\bigr\|_{[s,t]}^{2p}\leq
C_p|t-s|^p, \qquad\bE \bigl[\bigl(\bigl|\phi ^n\bigr|_t^s
\bigr)^{2p} \bigr]\leq C_p|t-s|^p.
\end{equation}

For general cases, choose $1\leq l<m-1<m\leq2^n$ such that
$t_{l-1}\leq s\leq t_l<t_{m-1}\leq t\leq t_m$. We get by It\^o's formula,
\begin{eqnarray*}
&&\d \bigl(e^{-({2}/\gamma)\vph(Y^n_t)}\bigl|Y^n_t-Y^n_s\bigr|^2
\bigr)
\\
&&\qquad=:U^n_s(t)\,\d w_t+U^n_s(t)
\,\d w^n_t+V^n_s(t)\,\d
t+Z^n_s(t)\,\d t+A^5_t,
\end{eqnarray*}
where according to (C),
\begin{eqnarray*}
A^5_t:=e^{-({2}/\gamma)\vph(Y^n_t)} \biggl[2{\bigl\langle
}Y^n_t-Y^n_s,\d\phi
^n_t{\bigr\rangle} -\frac{2}\gamma\bigl|Y^n_t-Y^n_s\bigr|^2{
\bigl\langle}D\vph\bigl(Y^n_t\bigr),\d\phi
^n_t{\bigr\rangle} \biggr]\leq0
\end{eqnarray*}
and
\begin{eqnarray*}
U^n_s(t)&:=&e^{-({2}/\gamma)\vph(Y^n_t)} \biggl(2
\bigl(Y^n_t-Y^n_s\bigr)-
\frac{2}\gamma \bigl|Y^n_t-Y^n_s\bigr|^2D
\vph\bigl(Y^n_t\bigr) \biggr)\sigma\bigl(Y^n_t
\bigr),
\\
V^n_s(t)&:=&e^{-({2}/\gamma)\vph(Y^n_t)} \biggl(2
\bigl(Y^n_t-Y^n_s\bigr)-
\frac{2}\gamma\bigl |Y^n_t-Y^n_s\bigr|^2D
\vph\bigl(Y^n_t\bigr) \biggr) \bigl(\sigma
\bigl(Y^n_t\bigr)\dot{h}_t+\tilde {b}
\bigl(Y^n_t\bigr) \bigr),
\\
Z^n_s(t)&:=&e^{-({2}/\gamma)\vph(Y^n_t)} \biggl[\operatorname{tr}\bigl(
\sigma \sigma ^*\bigr) \bigl(Y^n_t\bigr)-
\frac{1}\gamma\bigl|Y^n_t-Y^n_s\bigr|^2
\operatorname{tr}\bigl(D\vph\sigma \sigma ^*\bigr) \bigl(Y^n_t
\bigr)
\\
&&\hspace*{59pt}{}-\frac{4}\gamma\bigl(Y^n_t-Y^n_s
\bigr)\sigma\bigl(Y^n_t\bigr)D\vph\bigl(Y^n_t
\bigr)\sigma \bigl(Y^n_t\bigr)
\\
&&\hspace*{107pt}{}+\frac{2}{\gamma^2}\bigl|Y^n_t-Y^n_s\bigr|^2\bigl|D
\vph\bigl(Y^n_t\bigr)\sigma\bigl(Y^n_t
\bigr)\bigr|^2 \biggr].
\end{eqnarray*}
By the conditions on $\sigma, b, \vph$,
\begin{eqnarray*}
\bigl|U^n_s(t)\bigr|&\leq& C\bigl(\bigl|Y^n_t-Y^n_s\bigr|+\bigl|Y^n_t-Y^n_s\bigr|^2
\bigr),
\\
\bigl|U^n_s(t)-U^n_s
\bigl(t'\bigr)\bigr|&\leq& C\bigl|Y^n_t-Y^n_{t'}\bigr|
\bigl(1+\bigl|Y^n_t-Y^n_s\bigr|\bigr)
\\
&&{}+C\bigl|Y^n_t-Y^n_{t'}\bigr|
\bigl(\bigl|Y^n_t-Y^n_s\bigr|+\bigl|Y^n_{t'}-Y^n_s\bigr|+\bigl|Y^n_{t'}-Y^n_s\bigr|^2
\bigr)
\\
&\leq &C\bigl|Y^n_t-Y^n_{t'}\bigr|
\bigl(1+\bigl|Y^n_t-Y^n_s\bigr|
\bigr)+C\bigl|Y^n_t-Y^n_{t'}\bigr|^2
\bigl(1+\bigl|Y^n_{t'}-Y^n_{s}\bigr|^2\bigr),
\\
\bigl|V^n_s(t)\bigr|&\leq& C\bigl(\bigl|Y^n_t-Y^n_s\bigr|+\bigl|Y^n_t-Y^n_s\bigr|^2
\bigr) \bigl(1+|\dot{h}_t|\bigr),
\\
\bigl|Z^n_s(t)\bigr|&\leq& C\bigl(1+\bigl|Y^n_t-Y^n_s\bigr|+\bigl|Y^n_t-Y^n_s\bigr|^2
\bigr).
\end{eqnarray*}
Thus
\begin{eqnarray*}
&&\bE\bigl|Y^n_t-Y^n_s\bigr|^{4p}
\\
&&\qquad\leq C_p\bE \biggl(\biggl|\int_s^tU^n_s(r)
\,\d w_r\biggr|+\int_s^t\bigl|U^n_s(r)\bigr|
\bigl|\dot {w}^n_r\bigr|\,\d r+\biggl|\int_s^tV^n_s(r)
\,\d r\biggr|\\
&&\hspace*{227pt}{}+\int_s^t\bigl|Z^n_s(r)\bigr|
\,\d r \biggr)^{2p}.
\end{eqnarray*}

Using the BDG inequality we get
\begin{eqnarray*}
&&\bE \biggl(\int_s^tU^n_s(r)
\,\d w_r \biggr)^{2p}
\\
&&\qquad\leq C_p\bE \biggl(\int_s^t\bigl|U^n_s(r)\bigr|^2
\,\d r \biggr)^p
\\
&&\qquad\leq C_p(t-s)^{p-1}\bE \biggl(\int_s^t
\bigl(\bigl|Y^n_r-Y^n_s\bigr|^{2p}+\bigl|Y^n_r-Y^n_s\bigr|^{4p}
\bigr)\,\d r \biggr)
\\
&&\qquad\leq C_p(t-s)^p+C_p\bE \biggl(\int
_s^t\bigl|Y^n_r-Y^n_s\bigr|^{4p}
\,\d r \biggr),
\\
&&\bE \biggl(\int_s^tU^n_s(r)
\,\d w^n_r \biggr)^{2p}
\\
&&\qquad\leq C_p\bE \biggl[ \biggl(\int_s^t
\bigl(U^n_s(r)-U^n_s(
\bar{r}_n\vee s)\bigr)\,\d w^n_r+\int
_s^tU^n_s(
\bar{r}_n\vee s)\,\d w^n_r
\biggr)^{2p} \biggr].
\end{eqnarray*}
Note that $\int_0^\cdot U^n_s(\bar{r}_n)\,\d w^n_r$ is the piecewise
linear interpolation of $M^n_\cdot:= \break \int_0^{\cdot-\Delta}U^n_s(\pi
_n(r))\,\d w_r$ with
\[
\pi_n(r):=\max\{t_k;t_k\leq r\}.
\]
Thus by Doob's inequality and Lemma~\ref{estimatevarphi} we get
\begin{eqnarray*}
&&\bE\biggl|\int_s^tU^n_s(
\bar{r}_n\vee s)\,\d w^n_r\biggr|^{2p}
\\
&&\qquad\leq C_p\bE \biggl(\int_{s}^{t}\bigl|U^n_s(
\bar{r}_n\vee s)\bigr|^2\,\d r \biggr)^p
\\
&&\qquad\leq C_p|t-s|^{p-1}\bE \biggl(\int_s^t
\bigl(\bigl|Y^n_{\bar{r}_n\vee
s}-Y^n_s\bigr|^{2p}+\bigl|Y^n_{\bar{r}_n\vee s}-Y^n_s\bigr|^{4p}
\bigr)\,\d r \biggr)
\\
&&\qquad\leq C_p|t-s|^{p-1}\bE \biggl(\int_s^t
\bigl(\bigl|Y^n_r-Y^n_{\bar{r}_n\vee
s}\bigr|^{2p}+\bigl|Y^n_r-Y^n_s\bigr|^{2p}
\bigr)\,\d r \biggr)
\\
&&\qquad\quad{}+C_p|t-s|^{p-1}\bE \biggl(\int_s^t
\bigl(\bigl|Y^n_r-Y^n_{\bar{r}_n\vee
s}\bigr|^{4p}+\bigl|Y^n_r-Y^n_s\bigr|^{4p}
\bigr)\,\d r \biggr)
\\
&&\qquad\leq C_p|t-s|^p+C_p|t-s|^{p-1}\int
_s^t\bE\bigl|Y^n_r-Y^n_s\bigr|^{4p}
\,\d r
\\
&&\qquad\leq C_p|t-s|^p+C_p\int
_s^t\bE\bigl|Y^n_r-Y^n_s\bigr|^{4p}
\,\d r,
\\
&&\bE \biggl[\biggl|\int_s^t\bigl(U^n_s(r)-U^n_s(
\bar{r}_n\vee s)\bigr)\,\d w^n_r\biggr|^{2p}
\biggr]
\\
&&\qquad\leq C_p\bE \biggl[ \biggl(\int_s^t\bigl|Y^n_r-Y^n_{\bar{r}_n\vee
s}\bigr|
\bigl(1+\bigl|Y^n_r-Y^n_s\bigr|\bigr)
\frac{|w_{\hat{r}_n}-w_{\bar{r}_n\vee
s}|}{\Delta}\,\d r \biggr)^{2p} \biggr]
\\
&&\qquad\quad{}+C_p\bE \biggl[ \biggl(\int_s^t\bigl|Y^n_r-Y^n_{\bar{r}_n\vee
s}\bigr|^2
\bigl(1+\bigl|Y^n_r-Y^n_{s}\bigr|^2\bigr)
\frac{|w_{\hat{r}_n}-w_{\bar
{r}_n\vee s}|}{\Delta}\,\d r \biggr)^{2p} \biggr]
\\
&&\qquad\leq C_p\bE \biggl[ \biggl(\int_s^t\bigl|Y^n_r-Y^n_{\bar{r}_n\vee s}\bigr|
\frac
{|w_{\hat{r}_n}-w_{\bar{r}_n\vee s}|}{\Delta}\,\d r \biggr)^{2p} \biggr]
\\
&&\qquad\quad{}+C_p\bE \biggl[ \biggl(\int_s^t
\biggl(\bigl|Y^n_r-Y^n_s\bigr|^2+\bigl|Y^n_r-Y^n_{\bar
{r}_n\vee
s}\bigr|^2
\frac{|w_{\hat{r}_n}-w_{\bar{r}_n\vee s}|^2}{\Delta^2}\biggr)\d r \biggr)^{2p} \biggr]
\\
&&\qquad\quad{}+C_p\bE \biggl[ \biggl(\int_s^t\bigl|Y^n_r-Y^n_{\bar{r}_n\vee
s}\bigr|^2
\bigl(1+\bigl|Y^n_r-Y^n_{s}\bigr|^2\bigr)
\frac{|w_{\hat{r}_n}-w_{\bar
{r}_n\vee s}|}{\Delta}\,\mathrm{d} r \biggr)^{2p} \biggr]
\\
&&\qquad\leq C_p|t-s|^{2p}+C_p(t-s)^{2p-1}
\int_s^t\bE\bigl|Y^n_r-Y^n_s\bigr|^{4p}
\,\d r
\\
&&\qquad\leq C_p|t-s|^p+C_p\int
_s^t\bE\bigl|Y^n_r-Y^n_s\bigr|^{4p}
\,\d r.
\end{eqnarray*}
Moreover,
\begin{eqnarray*}
&&\bE \biggl(\biggl|\int_s^tV^n_s(r)
\,\d r\biggr|^{2p}+\biggl|\int_s^tZ^n_s(r)
\,\d r\biggr|^{2p} \biggr)
\\
&&\qquad\leq\bE \biggl[ \biggl(\int_s^t
\bigl(\bigl|Y^n_r-Y^n_s\bigr|+\bigl|Y^n_r-Y^n_s\bigr|^2
\bigr) \bigl(1+|\dot {h}_r|\bigr)\,\d r \biggr)^{2p}
\biggr]+(t-s)^{2p}
\\
&&\qquad\leq C_p|t-s|^{p}\biggl(1+\int_s^t|
\dot{h}_r|^2\d r\biggr)^p\\
&&\qquad\quad{}+C_p
\int_s^t\bE \bigl|Y^n_r-Y^n_s\bigr|^{4p}
\,\d r\biggl(1+\int_s^t|\dot{h}_r|^2
\,\d r\biggr)^p.
\end{eqnarray*}

Summing up we have
\[
\bE\bigl|Y^n_t-Y^n_s\bigr|^{4p}
\leq C_p|t-s|^{p}+C_p\int_s^t
\bE\bigl|Y^n_r-Y^n_s\bigr|^{4p}
\,\d r\biggl(1+\int_s^t|\dot {h}_r|^2
\,\d r\biggr)^p,
\]
which together with Gronwall's lemma yields
\[
\bE\bigl|Y^n_t-Y^n_s\bigr|^{4p}
\leq C_p|t-s|^{p}.
\]
It follows from this estimate and Lemma~\ref{sigmabar} that
\[
\bE\bigl|\phi^n_t-\phi^n_s\bigr|^{4p}
\leq C_p|t-s|^p.
\]

Now (\ref{3.8.2theta}) holds due to Kolmogorov's continuity criterion.
\end{pf*}

%pr3.4 #&#
\begin{proposition}\label{phiBV}
\[
\bE\sup_{t\in[0,1]}\bigl|Y^n_t\bigr|^{2p}<C_p
\bigl(1+|x|^{2p}\bigr),\qquad \sup_n\bE \bigl[\bigl(\bigl|
\phi^n\bigr|_1\bigr)^{2p} \bigr]<C_p.
\]
\end{proposition}

\begin{pf}
Using Proposition~\ref{holderC}, choose a $\theta\in[0,\frac{1}4)$, and
we get
\begin{eqnarray*}
\bE\sup_{t\in[0,1]}\bigl|Y^n_t\bigr|^{2p}&
\leq& 2^{2p-1}\bE\sup_{t\in[0,1]}\bigl|Y^n_t-x\bigr|^{2p}+2^{2p-1}|x|^{2p}
\\
&\leq&2^{2p-1}\bE \biggl(\sup_{t\in[0,1]}
\frac
{|Y^n_t-x|}{|t|^{\theta
}}|t|^{\theta} \biggr)^{2p}+2^{2p-1}|x|^{2p}
\\
&\leq&2^{2p-1}\bE\bigl\|Y^n\bigr\|_{[0,1],\theta}^{2p}+2^{2p-1}|x|^{2p}
\\
&\leq&C_p\bigl(1+|x|^{2p}\bigr).
\end{eqnarray*}

Similar to \cite{stroock}, Theorem~3.6, by (D) we get for all
$0\leq s<t\leq1$,
\begin{eqnarray*}
\bigl|\phi^n\bigr|_t^s\leq C \bigl(|t-s|R^{-4}
\bigl\|Y^n\bigr\|^4_{[s,t],\theta}+1 \bigr)\| \phi
^n\|_{[s,t]}.
\end{eqnarray*}
From this and Proposition~\ref{holderC},
\begin{eqnarray*}
\bE \bigl[\bigl(\bigl|\phi^n\bigr|_1\bigr)^{2p} \bigr]&
\leq& C_p\bE \bigl[ \bigl(R^{-4}\bigl\|Y^n
\bigr\|^4_{[0,1],\theta}+1 \bigr)^{2p}\bigl\|\phi ^n\bigr\|
_{[0,1]}^{2p} \bigr]
\\
&\leq&C_{R,p}<\infty.
\end{eqnarray*}
\upqed\end{pf}

%pr3.5 #&#
\begin{proposition}\label{convergence3}
\begin{eqnarray*}
\sup_{1\leq k\leq2^n}\bE\bigl(\bigl|Y^n_{t_k}-Z_{t_k}\bigr|^2
\bigr)\leq C\biggl[\Delta ^{\theta
/2}+\sup_{2\leq k\leq2^n}\biggl(\int
_{t_{k-2}}^{t_k}|\dot{h}_s|^2\,\d s
\biggr)^{1/2}\biggr],\qquad \theta\in(0,1).
\end{eqnarray*}
\end{proposition}

\begin{pf}
Set
\begin{eqnarray*}
\mu_n(t)&:=&e^{-({2}/{\gamma}) (\vph(Y^n_t)+\vph(Z_t)
)},\qquad m_n(t):=
\mu_n(t)\bigl|Y^n_t-Z_t\bigr|^2,\\
a(t)&:=&\bE\bigl(m_n(t)\bigr).
\end{eqnarray*}
Using the condition $\vph\in\mathcal{C}_b^2$, Lemma~\ref{holderZ} and
(\ref{holderestimatevarphi}) it is trivial to prove the following:

%le3.9 #&#
\begin{lemma}\label{mulemma}
\begin{eqnarray*}
\bE \Bigl(\sup_{t,t'\in[t_{k-2},t_k]}\bigl|\mu_n(t)-
\mu_n\bigl(t'\bigr)\bigr|^2 \Bigr)&\leq&C\Delta,
\\
 \bE \Bigl(\sup_{t,t'\in[t_{k-2},t_k]}\bigl|m_n(t)-m_n
\bigl(t'\bigr)\bigr| \Bigr)&\leq& C\Delta^{1/2}.
\end{eqnarray*}
\end{lemma}

For all $t_{k-1}\leq t\leq t_k$, $2\leq k\leq2^n$,
\begin{eqnarray*}
\d\mu_n(t)\bigl|Y^n_t-Z_t\bigr|^2
&=&\sum_{i=1}^{11}\,\d I_i(t)+2
\mu_n(t){\bigl\langle}Y^n_t-Z_t,
\d\phi _t^n-\d\psi _t{\bigr\rangle}
\\
&&{}-\frac{2}{\gamma} \mu_n(t)\bigl|Y^n_t-Z_t\bigr|^2
\bigl({\bigl\langle}D\vph \bigl(Y_t^n\bigr),\d \phi
_t^n{\bigr\rangle}+{\bigl\langle}D\vph(Z_t),
\d\psi_t{\bigr\rangle} \bigr),
\end{eqnarray*}
where
\begin{eqnarray*}
I_1(s)&:=&2\int_{t_{k-1}}^s
\mu_n(t){\bigl\langle}Y^n_t-Z_t,
\sigma \bigl(Y^n_t\bigr)-\sigma (Z_t){\bigr
\rangle}\dot{h}_t\,\d t,
\\
I_2(s)&:=&2\int_{t_{k-1}}^s
\mu_n(t){\bigl\langle}Y^n_t-Z_t,
\tilde {b}\bigl(Y^n_t\bigr)-b(Z_t){\bigr
\rangle} \,\d t,
\\
I_3(s)&:=&2\int_{t_{k-1}}^s
\mu_n(t){\bigl\langle}Y^n_t-Z_t,
\sigma \bigl(Y^n_t\bigr)-\sigma (Y_{\bar{t}_n}){\bigr
\rangle}\,\d w_t,
\\
I_4(s)&:=&2\int_{t_{k-1}}^s
\mu_n(t) \bigl(\operatorname{tr}\bigl(\sigma\sigma ^*\bigr)
\bigl(Y_t^n\bigr)-\operatorname{tr}\bigl(\sigma\sigma^*\bigr)
\bigl(Y_{\bar{t}_n}^n\bigr) \bigr)\,\d t
\\
&&{}+\int_{t_{k-1}}^s\bigl(\mu_n(t)-
\mu_n(\bar{t}_n)\bigr)\operatorname{tr}\bigl(\sigma \sigma ^*
\bigr) \bigl(Y_{\bar{t}_n}^n\bigr)\,\d t,
\\
I_5(s)&:=&2\int_{t_{k-1}}^s
\mu_n(t) \biggl({\bigl\langle}Y^n_t-Z_t,
\sigma \bigl(Y^n_{\bar
{t}_n}\bigr) \bigl(\d w_t-
\d{w}^n_t\bigr){\bigr\rangle}
\\
&&\hspace*{74pt}{}+\int_{t_{k-1}}^s\mu_n(
\bar{t}_n)\operatorname{tr}\bigl(\sigma\sigma ^*\bigr)
\bigl(Y_{\bar
{t}_n}^n\bigr) \biggr)\,\d t,
\\
I_6(s)&:=&-2\int_{t_{k-1}}^s
\mu_n(t){\bigl\langle}Y^n_t-Z_t,
\bigl(\sigma \bigl(Y^n_t\bigr)-\sigma
\bigl(Y^n_{\bar{t}_n}\bigr) \bigr)\dot{w}^n_t{
\bigr\rangle}\,\d t,
\\
I_7(s)&:=&-\frac{2}{\gamma} \int_{t_{k-1}}^s
\mu _n(t)\bigl|Y^n_t-Z_t\bigr|^2{
\bigl\langle} D\vph\bigl(Y_t^n\bigr),\sigma
\bigl(Y^n_t\bigr) \bigl(\d w_t-\d
w^n_t\bigr){\bigr\rangle},
\\
I_8(s)&:=&-\frac{2}{\gamma} \int_{t_{k-1}}^s
\mu_n(t)\bigl|Y^n_t-Z_t\bigr|^2
\\
&&\hspace*{40pt}{}\times \bigl({\bigl\langle}D\vph\bigl(Y_t^n\bigr),
\sigma\bigl(Y^n_t\bigr)\dot{h}_t{\bigr\rangle
}+{\bigl\langle}D\vph (Z_t),\sigma (Z_t)
\dot{h}_t{\bigr\rangle} \bigr)\,\d t,
\\
I_9(s)&:=&-\frac{2}{\gamma} \int_{t_{k-1}}^s
\mu_n(t)\bigl|Y^n_t-Z_t\bigr|^2
\\
&&\hspace*{40pt}{}\times \biggl({\bigl\langle}D\vph\bigl(Y_t^n\bigr),
\tilde{b}\bigl(Y^n_t\bigr){\bigr\rangle }+{\bigl\langle}D
\vph (Z_t),{b}(Z_t){\bigr\rangle}\\
&&\hspace*{99pt}{} +\frac{1}{2}
\operatorname{tr}\bigl(D^2\vph\bigl(Y^n_t\bigr)
\sigma\sigma^*\bigl(Y^n_t\bigr)\bigr) \biggr)\,\d t,
\\
I_{10}(s)&:=&-\frac{4}{\gamma} \int_{t_{k-1}}^s
\mu_n(t) \biggl(\sum_i{\bigl\langle}
D\vph\bigl(Y_t^n\bigr),\sigma\bigl(Y^n_t
\bigr)e_i{\bigr\rangle} {\bigl\langle}Y^n_t-Z_t,
\sigma \bigl(Y^n_t\bigr)e_i{\bigr\rangle}
\biggr)\,\d t,
\\
I_{11}(s)&:=&\frac{2}{\gamma^2}\int_{t_{k-1}}^s
\mu _n(t)\bigl|Y^n_t-Z_t\bigr|^2\bigl|D
\vph\bigl(Y_t^n\bigr)\sigma\bigl(Y^n_t
\bigr)\bigr|^2\,\d t.
\end{eqnarray*}

By (C),
\begin{eqnarray*}
{\bigl\langle}Y^n_t-Z_t,\,\d
\phi_t^n-\d\psi_t{\bigr\rangle}-
\frac{1}{\gamma} \bigl|Y^n_t-Z_t\bigr|^2
\bigl({\bigl\langle} D\vph\bigl(Y_t^n\bigr),\d
\phi_t^n{\bigr\rangle}+{\bigl\langle}D
\vph(Z_t),\d\psi _t{\bigr\rangle} \bigr)\leq0.
\end{eqnarray*}
Thus
%
%e3.9 #&#
\begin{equation}
\label{akminus} \mu_n(t_k)\bigl|Y^n_{t_k}-Z_{t_k}\bigr|^2
\leq\mu _n(t_{k-1})\bigl|Y^n_{t_{k-1}}-Z_{t_{k-1}}\bigr|^2+
\sum_{i=1}^{11}I_i(t_k).
\end{equation}

By the hypotheses $\sigma\in\cC_b^2,  \vph\in\cC_b^2$,
\begin{eqnarray*}
&&\bigl|I_1(t_k)+I_8(t_k)\bigr|
\\
&&\qquad\leq\biggl|\int_{t_{k-1}}^{t_k}2\mu_n(s){\bigl
\langle}Y_s^n-Z_s,\sigma
\bigl(Y^n_s\bigr)-\sigma (Z_s){\bigr
\rangle}\dot{h}_s\,\d s\biggr|
\\
&&\qquad\quad{}+\biggl|\frac{2}{\gamma}\int_{t_{k-1}}^{t_k}\mu
_n(s)\bigl|Y_s^n-Z_s\bigr|^2
\bigl({\bigl\langle}D\vph\bigl(Y^n_s\bigr),\sigma
\bigl(Y^n_s\bigr){\bigr\rangle }\\
&&\hspace*{159pt}{}+{\bigl\langle}D\vph
(Z_s),\sigma (Z_s){\bigr\rangle}\bigr)
\dot{h}_s\,\d s\biggr|
\\
&&\qquad\leq C\int_{t_{k-1}}^{t_k}m_n(s)|
\dot{h}_s|\,\d s.
\end{eqnarray*}

For $I_3$, $\bE(I_3(t_k))=\bE(\int_{t_{k-1}}^{t_k}2{\langle}
Y_t^n-Z_t,\sigma
(Y^n_t)-\sigma(Y^n_{\bar{t}_n}){\rangle}\,\d w_t)=0$.

Throughout the proof we need several lemmas which will be proved afterward.
Now we deal with the terms $I_2$ and $I_6$. Note that for $I_2$,
\begin{eqnarray*}
I_2&=&2\int_{t_{k-1}}^{t_k}
\mu_n(t){\bigl\langle}Y_t^n-Z_t,
\tilde {b}\bigl(Y^n_t\bigr)-b(Z_t){\bigr
\rangle} \,\d t
\\
&=&2\int_{t_{k-1}}^{t_k}\mu_n(t){\bigl
\langle}Y_t^n-Z_t-\bigl(Y^n_{\bar
{t}_n}-Z_{\bar
{t}_n}
\bigr),\tilde{b}\bigl(Y^n_t\bigr)-b(Z_t){
\bigr\rangle}\,\d t
\\
&&{}+2\int_{t_{k-1}}^{t_k}\bigl(\mu_n(t)-
\mu_n(\bar{t}_n)\bigr){\bigl\langle }Y^n_{\bar
{t}_n}-Z_{\bar{t}_n},
\tilde{b}\bigl(Y^n_t\bigr)-b(Z_t){\bigr
\rangle}\,\d t
\\
&&{}+2\int_{t_{k-1}}^{t_k}\mu_n(
\bar{t}_n){\bigl\langle}Y^n_{\bar
{t}_n}-Z_{\bar
{t}_n},
\tilde{b}\bigl(Y^n_t\bigr)-\tilde{b}\bigl(Y^n_{\bar{t}_n}
\bigr)-b(Z_t)+b(Z_{\bar
{t}_n}){\bigr\rangle}\,\d t
\\
&&{}+2\int_{t_{k-1}}^{t_k}\mu_n(
\bar{t}_n){\bigl\langle}Y^n_{\bar
{t}_n}-Z_{\bar
{t}_n},{b}
\bigl(Y^n_{\bar{t}_n}\bigr)-b(Z_{\bar{t}_n}){\bigr\rangle}\,\d
t
\\
&&{}+2\int_{t_{k-1}}^{t_k}\mu_n(t_{k-2}){
\bigl\langle} Y^n_{t_{k-2}}-Z_{t_{k-2}},\tilde {b}
\bigl(Y^n_{t_{k-2}}\bigr)-{b}\bigl(Y^n_{t_{k-2}}
\bigr){\bigr\rangle}\,\d t
\\
&=:&\sum_{i=1}^5I_{2,i}.
\end{eqnarray*}
Taking expectations and applying Lemmas \ref{holderZ}, \ref
{estimatevarphi} and \ref{mulemma}, we get
\begin{eqnarray*}
&&\Biggl|\bE\sum_{i=1}^4I_{2,i}\Biggr|
\\
&&\qquad\leq C\Delta^{3/2}+\biggl|\bE\int_{t_{k-1}}^{t_k}
\bigl(\mu_n(t)-\mu_n(\bar {t}_n)
\bigr)\bigl|Y^n_{\bar{t}_n}-Z_{\bar{t}_n}\bigr|\bigl|\tilde{b}
\bigl(Y^n_t\bigr)-b(Z_t)\bigr|\,\d t\biggr|
\\
&&\qquad\quad{}+C\biggl|\bE\int_{t_{k-1}}^{t_k}\mu_n(
\bar{t}_n)\bigl|Y^n_{\bar
{t}_n}-Z_{\bar
{t}_n}\bigr|\bigl|Y_t^n-Z_t-
\bigl(Y^n_{\bar{t}_n}-Z_{\bar{t}_n}\bigr)\bigr|\,\d
t\biggr|+Ca(t_{k-2})\Delta
\\
&&\qquad\leq C\Delta^{3/2}+Ca(t_{k-2})\Delta.
\end{eqnarray*}
Note that
\begin{eqnarray*}
I_6&=&-2\int_{t_{k-1}}^{t_k}
\mu_n(t){\bigl\langle}Y_t^n-Z_t-
\bigl(Y^n_{\bar
{t}_n}-Z_{\bar
{t}_n}\bigr),\sigma
\bigl(Y^n_t\bigr)-\sigma\bigl(Y^n_{\bar{t}_n}
\bigr){\bigr\rangle}\,\d w^n_t
\\
&&{}-2\int_{t_{k-1}}^{t_k}\bigl(\mu_n(t)-
\mu_n(\bar{t}_n)\bigr){\bigl\langle }Y^n_{\bar
{t}_n}-Z_{\bar{t}_n},
\sigma\bigl(Y^n_t\bigr)-\sigma\bigl(Y^n_{\bar{t}_n}
\bigr){\bigr\rangle }\,\d w^n_t
\\
&&{}-2\int_{t_{k-1}}^{t_k}\mu_n(
\bar{t}_n){\bigl\langle}Y^n_{\bar
{t}_n}-Z_{\bar
{t}_n},
\sigma\bigl(Y^n_t\bigr)-\sigma\bigl(Y^n_{\bar{t}_n}
\bigr){\bigr\rangle}\,\d w^n_t=:\sum
_i I_{6,i}
\end{eqnarray*}
and
\begin{eqnarray*}
\bigl|\bE(I_{6,1}+I_{6,2})\bigr|\leq C\Delta^{3/2}.
\end{eqnarray*}
As for $I_{6,3}+I_{2,5}$, note that $I_{6,3}+I_{2,5}=-2A_k^n$, where
\begin{equation}\qquad
\label{Ak} A_k^n:=\int_{t_{k-1}}^{t_k}
\mu_n(\bar {t}_n){\biggl\langle}Y^n_{\bar
{t}_n}-Z_{\bar{t}_n},
\bigl(\sigma\bigl(Y^n_t\bigr)-\sigma\bigl(Y^n_{\bar{t}_n}
\bigr)\bigr)\dot {w}^n_t-\frac{1}{2}(\nabla\sigma)
\sigma\bigl(Y^n_{\bar{t}_n}\bigr){\biggr\rangle }\,\d t,
\end{equation}
and by Lemma~\ref{Aklemma},
\begin{eqnarray*}
\Biggl|\bE \Biggl(\sum_{i=1}^{2^n}A_i^n
\Biggr)\Biggr|\leq C\biggl[\Delta^{1/2}+\sup_{2\leq
k\leq2^n}\biggl(
\int_{t_{k-2}}^{t_k}|\dot{h}_s|^2
\,\d s\biggr)^{1/2}\biggr].
\end{eqnarray*}

By Lemmas \ref{estimatevarphi} and~\ref{mulemma},
\begin{eqnarray*}
&&\bigl|\bE I_4(t_k)\bigr|\\
&&\qquad\leq\biggl|\bE\int_{t_{k-1}}^{t_k}
\mu_n(t) \bigl(\operatorname {tr}\bigl(\sigma \sigma^*\bigr)
\bigl(Y^n_t\bigr)-\operatorname{tr}\bigl(\sigma\sigma^*\bigr)
\bigl(Y^n_{\bar{t}_n}\bigr)\bigr)\,\d t\biggr|
\\
&&\qquad\quad{}+\biggl|\bE\int_{t_{k-1}}^{t_k}\bigl(\mu_n(t)-
\mu_n(\bar{t}_n)\bigr)\operatorname {tr}\bigl(\sigma \sigma^*
\bigr) \bigl(Y^n_{\bar{t}_n}\bigr)\,\d t\biggr|\leq C
\Delta^{3/2},
\\
&&\bigl|\bE\bigl[I_9(t_k)+I_{11}(t_k)
\bigr]\bigr|\\
&&\qquad\leq\frac{2}{\gamma}\biggl|\bE\int_{t_{k-1}}^{t_k}
\bigl(m_n(s)-m_n(\bar {s}_n)+m_n(
\bar{s}_n)\bigr)
\\
&&\hspace*{36pt}\qquad\quad{}\times \biggl(\bigl|D\vph\bigl(Y^n_s\bigr)\bigr|\tilde{b}
\bigl(Y^n_s\bigr)\bigl|+\bigl|D\vph (Z_s)\bigr|\bigl|{b}(Z_s)\bigr|\\
&&\hspace*{148pt}{}+
\frac{1}2\operatorname{tr}\bigl(D^2\vph\sigma\sigma^*
\bigl(Y^n_s\bigr)\bigr) \biggr)\,\d s\biggr|
\\
&&\qquad\quad{}+\frac{2}{\gamma^2}\biggl|\bE\int_{t_{k-1}}^{t_k}
\bigl(m_n(s)-m_n(\bar {s}_n)+m_n(
\bar{s}_n)\bigr)\bigl|D\vph\bigl(Y^n_s\bigr)
\sigma\bigl(Y^n_s\bigr)\bigr|^2\,\d s\biggr|
\\
&&\qquad\leq C\Delta^{3/2}+Ca(t_{k-2})\Delta.
\end{eqnarray*}

For $I_5$, we have
\begin{eqnarray*}
I_5(t_k) &=&2\int_{t_{k-1}}^{t_k}
\mu_n(\bar{t}_n){\bigl\langle}Y_t^n-Z_t,
\sigma \bigl(Y^n_{\bar
{t}_n}\bigr){\bigr\rangle}\bigl(\d
w_t-\d w^n_t\bigr)
\\
&&{}+\int_{t_{k-1}}^{t_k}\mu_n(
\bar{t}_n)\operatorname{tr}\bigl(\sigma\sigma ^*\bigr)
\bigl(Y^n_{\bar{t}_n}\bigr)\,\d t
\\
&&{}+2\int_{t_{k-1}}^{t_k}\bigl(\mu_n(t)-
\mu_n(\bar{t}_n)\bigr){\bigl\langle}
Y_t^n-Z_t,\sigma \bigl(Y^n_{\bar{t}_n}
\bigr){\bigr\rangle}\bigl(\d w_t-\d w^n_t
\bigr)
\\
&=:&I_{5,1}+I_{5,2}.
\end{eqnarray*}
However, by Lemma~\ref{lemmaI51},
\begin{eqnarray*}
\Biggl|\bE \Biggl(\sum_{i=0}^{2^n}{I}_{5,1}(t_i)
\Biggr)\Biggr|\leq C\Delta^{\theta
/2}\qquad \forall\theta\in(0,1).
\end{eqnarray*}
With respect to $I_{5,2}$,
\begin{eqnarray*}
I_{5,2}&=&-\frac{4}{\gamma}\int_{t_{k-1}}^{t_k}
\int_{\bar
{t}_n}^t\mu _n(s){\bigl\langle}D
\vph\bigl(Y^n_s\bigr),\sigma\bigl(Y^n_s
\bigr) \bigl(\d w_s-\d w^n_s\bigr){\bigr
\rangle}
\\
&&\hspace*{54pt}{}\times{\bigl\langle}Y_t^n-Z_t,\sigma
\bigl(Y^n_{\bar{t}_n}\bigr){\bigr\rangle}\bigl(\d
w_t-\d w^n_t\bigr)
\\
&&{}-\frac{4}{\gamma}\int_{t_{k-1}}^{t_k}\int
_{\bar{t}_n}^t\mu _n(s) \bigl({\bigl
\langle} D\vph\bigl(Y^n_s\bigr),\tilde{b}
\bigl(Y^n_s\bigr){\bigr\rangle}\,\d s+{\bigl\langle}D\vph
(Z_s),{b}(Z_s){\bigr\rangle}\,\d s\bigr)
\\
&&\hspace*{59pt}{}\times{\bigl\langle}Y_t^n-Z_t,\sigma
\bigl(Y^n_{\bar{t}_n}\bigr){\bigr\rangle}\bigl(\d
w_t-\d w^n_t\bigr)
\\
&&{}-\frac{4}{\gamma}\int_{t_{k-1}}^{t_k}\int
_{\bar{t}_n}^t\mu _n(s) \bigl({\bigl
\langle} D\vph\bigl(Y^n_s\bigr),\sigma
\bigl(Y^n_s\bigr){\bigr\rangle}\dot{h}_s\,\d
s+{\bigl\langle}D\vph (Z_s),\sigma (Z_s){\bigr\rangle}
\dot {h}_s\,\d s\bigr)
\\
&&\hspace*{55pt}{}\times{\bigl\langle}Y_t^n-Z_t,\sigma
\bigl(Y^n_{\bar{t}_n}\bigr){\bigr\rangle}\bigl(\d
w_t-\d w^n_t\bigr)
\\
&&{}-\frac{4}{\gamma}\int_{t_{k-1}}^{t_k}\int
_{\bar{t}_n}^t\mu _n(s) \bigl({\bigl
\langle} D\vph\bigl(Y^n_s\bigr),\d\phi_s^n{
\bigr\rangle}+{\bigl\langle}D\vph(Z_s),\d\psi _s{\bigr
\rangle}\bigr)
\\
&&\hspace*{55pt}{}\times{\bigl\langle}Y_t^n-Z_t,\sigma
\bigl(Y^n_{\bar{t}_n}\bigr){\bigr\rangle}\bigl(\d
w_t-\d w^n_t\bigr)
\\
&&{}+\frac{2}{\gamma^2}\int_{t_{k-1}}^{t_k}\int
_{\bar{t}_n}^t\mu _n(s)\bigl|D\vph
\bigl(Y^n_s\bigr)\sigma\bigl(Y^n_s
\bigr)\bigr|^2\,\d s
\\
&&\hspace*{64pt}{}\times{\bigl\langle}Y_t^n-Z_t,\sigma
\bigl(Y^n_{\bar{t}_n}\bigr){\bigr\rangle}\bigl(\d
w_t-\d w^n_t\bigr)
\\
&&{}-\frac{1}{\gamma}\int_{t_{k-1}}^{t_k}\int
_{\bar{t}_n}^t\mu _n(s)\operatorname {tr}
\bigl(D^2\vph\bigl(Y^n_s\bigr)\sigma
\bigl(Y^n_s\bigr)\bigr)\,\d s
\\
&&\hspace*{60pt}{}\times{\bigl\langle}Y_t^n-Z_t,\sigma
\bigl(Y^n_{\bar{t}_n}\bigr){\bigr\rangle}\bigl(\d
w_t-\d w^n_t\bigr)
\\
&=:&\sum_{i=1}^6 I_{5,2,i}.
\end{eqnarray*}
Applying the BDG inequality, the conditions $\sigma\in\cC_b^2,  b\in
\cC
_b^1,  \vph\in\cC_b^2$, Lemmas \ref{holderZ}, \ref{estimatevarphi}
and Proposition~\ref{phiBV}, we get
\begin{eqnarray*}
|\bE I_{5,2,2}|&\leq&C\bE \biggl(\int_{t_{k-1}}^{t_k}C
\Delta ^2\bigl|Y_t^n-Z_t\bigr|^2
\bigl\|\sigma\bigl(Y^n_{\bar{t}_n}\bigr)\bigr\|^2\,\d t
\biggr)^{1/2}
\\
&&{}+\biggl|\bE\int_{t_{k-1}}^{t_k}C\Delta\bigl|Y_t^n-Z_t\bigr|
\bigl\|\sigma\bigl(Y^n_{\bar
{t}_n}\bigr)\bigr\| \Big|\,\d w^n_t\biggr|
\\
&\leq& C\Delta^{3/2},
\\
|\bE I_{5,2,3}|&\leq&C\bE \biggl(\int_{t_{k-1}}^{t_k}\bigl|Y_t^n-Z_t\bigr|^2
\bigl\| \sigma \bigl(Y^n_{\bar{t}_n}\bigr)\bigr\|^2\biggl|\int
_{\bar{t}_n}^t\dot{h}_s\,\d s\biggr|^2
\,\d t \biggr)^{1/2}
\\
&&{}+C\bE\int_{t_{k-1}}^{t_k}\bigl|Y_t^n-Z_t\bigr|
\bigl\|\sigma\bigl(Y^n_{\bar{t}_n}\bigr)\bigr\| \biggl|\int_{\bar{t}_n}^t
\dot{h}_s\,\d s\biggr|\bigl|\d w^n_t\bigr|
\\
&\leq&C\Delta^{1/2}\int_{t_{k-2}}^{t_k}|
\dot{h}_s|\,\d s,
\\
\bigl|\bE(I_{5,2,5}+I_{5,2,6})\bigr|&\leq& C\Delta^{3/2},
\\
|\bE I_{5,2,4}|&\leq&\frac{4}{\gamma}\biggl|\bE \biggl(\int
_{t_{k-1}}^{t_k}\int_{\bar{t}_n}^t
\mu_n(s){\bigl\langle}D\vph\bigl(Y^n_s
\bigr)-D\vph\bigl(Y^n_{\bar
{t}_n}\bigr),\d \phi
_s^n{\bigr\rangle}
\\
&&\hspace*{57pt}{}\times{\bigl\langle}Y_t^n-Z_t,\sigma
\bigl(Y^n_{\bar{t}_n}\bigr){\bigr\rangle}\bigl(\d
w_t-\d w^n_t\bigr) \biggr)\biggr|
\\
&&{}+\biggl|\bE \biggl(\int_{t_{k-1}}^{t_k}\int
_{\bar{t}_n}^t\mu _n(s){\bigl\langle}D\vph
(Z_s)-D\vph(Z_{\bar{t}_n}),\d\psi_s{\bigr\rangle}
\\
&&\hspace*{64pt}{}\times{\bigl\langle}Y_t^n-Z_t,\sigma
\bigl(Y^n_{\bar{t}_n}\bigr){\bigr\rangle}\bigl(\d
w_t-\d w^n_t\bigr) \biggr)\biggr|
\\
&&{}+\biggl|\bE \biggl(\int_{t_{k-1}}^{t_k}\int
_{\bar{t}_n}^t\mu _n(s) \bigl({\bigl
\langle} D\vph \bigl(Y^n_{\bar{t}_n}\bigr),\d
\phi_s^n{\bigr\rangle}+{\bigl\langle}D
\vph(Z_{\bar
{t}_n}),\d\psi_s{\bigr\rangle} \bigr)
\\
&&\hspace*{100pt}{}\times{\bigl\langle}Y_t^n-Z_t,\sigma
\bigl(Y^n_{\bar{t}_n}\bigr){\bigr\rangle}\bigl(\d
w_t-\d w^n_t\bigr) \biggr)\biggr|
\\
&\leq&C\Delta^{3/2}+C\bE G_k,
\end{eqnarray*}
where
%
%e3.10 #&#
\begin{eqnarray}
\label{defGk} G_k:=\max_{t\in[t_{k-1},t_k]}\bigl|Y^n_t-Z_t\bigr||
\Delta w_{k-1}|\times \bigl(\bigl|\phi ^n\bigr|_{t_k}^{t_{k-2}}+|
\psi|_{t_k}^{t_{k-2}}\bigr).
\end{eqnarray}
Again by Lemma~\ref{holderZ} and Proposition~\ref{phiBV},
\begin{eqnarray*}
\sum_{k=1}^{2^n}\bE G_k&\leq&
\bE \Bigl(\max_{1\leq k\leq2^n}\max_{t\in
[t_{k-1},t_k]}\bigl|Y^n_t-Z_t\bigr||
\Delta w_{k-1}|\times\bigl(\bigl|\phi^n\bigr|_1+|\psi
|_1\bigr) \Bigr)
\\
&\leq& \Bigl[\bE\Bigl(\sup_{0\leq t\leq1}\bigl|Y_t^n-Z_t\bigr|^{2p}
\Bigr)\bE\Bigl(\sup_{1\leq
k\leq2^n}|\Delta w_k|^{2p}
\Bigr) \Bigr]^{1/2p}\\
&&{}\times \bigl[\bE\bigl(\bigl|\phi^n\bigr|_1+|
\psi |_1\bigr)^q \bigr]^{1/q}
\\
&\leq&C\Delta^{(p-1)/2p},\qquad p, q>1, 1/p+1/q=1
\end{eqnarray*}
and
\begin{eqnarray*}
I_{5,2,1}&=&-\frac{4}{\gamma}\int_{t_{k-1}}^{t_k}
\int_{{\bar
{t}_n}}^t\bigl(\mu_n(s)-
\mu_n(\bar{t}_n)\bigr){\bigl\langle}D\vph
\bigl(Y^n_s\bigr),\sigma \bigl(Y^n_s
\bigr) \bigl(\d w_s-\d w^n_s\bigr){\bigr
\rangle}
\\
&&\hspace*{54pt}{}\times{\bigl\langle}Y_t^n-Z_t,\sigma
\bigl(Y^n_{\bar{t}_n}\bigr){\bigr\rangle}\bigl(\d
w_t-\d w^n_t\bigr)
\\
&&{}-\frac{4}{\gamma}\int_{t_{k-1}}^{t_k}\int
_{\bar{t}_n}^t\mu _n(\bar
{t}_n){\bigl\langle}D\vph\bigl(Y^n_s
\bigr),\sigma\bigl(Y^n_s\bigr) \bigl(\d
w_s-\d w^n_s\bigr){\bigr\rangle }
\\
&&\hspace*{58pt}{}\times{\bigl\langle}Y_t^n-Z_t-
\bigl(Y_{\bar{t}_n}^n-Z_{\bar{t}_n}\bigr),\sigma
\bigl(Y^n_{\bar
{t}_n}\bigr){\bigr\rangle}\bigl(\d
w_t-\d w^n_t\bigr)
\\
&&{}-\frac{4}{\gamma}\int_{t_{k-1}}^{t_k}\int
_{\bar{t}_n}^t\mu _n(\bar
{t}_n){\bigl\langle}D\vph\bigl(Y^n_s
\bigr),\bigl(\sigma\bigl(Y^n_s\bigr)-\sigma
\bigl(Y^n_{\bar
{t}_n}\bigr)\bigr) \bigl(\d w_s-\d
w^n_s\bigr){\bigr\rangle}
\\
&&\hspace*{58pt}{}\times{\bigl\langle}Y_{\bar{t}_n}^n-Z_{\bar{t}_n},\sigma
\bigl(Y^n_{\bar
{t}_n}\bigr){\bigr\rangle} \bigl(\d
w_t-\d w^n_t\bigr)
\\
&&{}-\frac{4}{\gamma}\int_{t_{k-1}}^{t_k}\int
_{\bar{t}_n}^t\mu _n(\bar
{t}_n){\bigl\langle}D\vph\bigl(Y^n_s
\bigr)-D\vph\bigl(Y^n_{\bar{t}_n}\bigr),\sigma
\bigl(Y^n_{\bar
{t}_n}\bigr) \bigl(\d w_s-\d
w^n_s\bigr){\bigr\rangle}
\\
&&\hspace*{58pt}{}\times{\bigl\langle}Y_{\bar{t}_n}^n-Z_{\bar{t}_n},\sigma
\bigl(Y^n_{\bar
{t}_n}\bigr){\bigr\rangle} \bigl(\d
w_t-\d w^n_t\bigr)
\\
&&{}-\frac{4}{\gamma}\int_{t_{k-1}}^{t_k}
\mu_n(\bar{t}_n){\bigl\langle }D\vph
\bigl(Y^n_{\bar
{t}_n}\bigr),\sigma\bigl(Y^n_{\bar{t}_n}
\bigr) (w_t-w_{\bar{t}_n}- \bigl(w^n_t-w^n_{\bar
{t}_n}
\bigr) {\bigr\rangle}
\\
&&\hspace*{41pt}{}\times{\bigl\langle}Y_{\bar{t}_n}^n-Z_{\bar{t}_n},\sigma
\bigl(Y^n_{\bar
{t}_n}\bigr){\bigr\rangle} \bigl(\d
w_t-\d w^n_t\bigr)
\\
&=:&\sum_{j=1}^5I_{5,2,1}^j.
\end{eqnarray*}
Using the BDG inequality, the fact that $\sigma\in\cC_b^2, \vph\in
\cC
_b^2$, Lemmas \ref{holderZ}, \ref{estimatevarphi}, \ref{mulemma} and
Proposition~\ref{phiBV}, we get
\begin{eqnarray*}
\bigl|\bE I_{5,2,1}^1\bigr| &\leq&C\Delta \bigl(\bE|\Delta
w_{k-1}|^2 \bigr)^{1/2}\leq C\Delta
^{3/2},
\\
\bigl|\bE I_{5,2,1}^2\bigr|&\leq& C\Delta^{1/2} \biggl(\bE
\biggl(\int_{t_{k-1}}^{t_k}\bigl|Y_t^n-Z_t-
\bigl(Y_{\bar
{t}_n}^n-Z_{\bar{t}_n}\bigr)\bigr|^2|
\Delta w_{k-1}|^2\Delta^{-1}\,\d t \biggr)
\biggr)^{1/2}
\\
&\leq& C\Delta^{3/2},
\\
\bigl|\bE I_{5,2,1}^j\bigr|&\leq&C\Delta^{3/2}, \qquad j=3,4
\end{eqnarray*}
and
%
%e3.11 #&#
\begin{eqnarray}\label{I5215}\qquad
&&\bE\bigl(I_{5,2,1}^5|\sF_{t_{k-2}}
\bigr)
\nonumber
\\[-8pt]
\\[-8pt]
\nonumber
&&\qquad=\frac{2\Delta}{\gamma}\mu_n({t_{k-2}}) \sum
_i{\bigl\langle}D\vph\bigl(Y^n_{t_{k-2}}
\bigr),\sigma \bigl(Y^n_{t_{k-2}}\bigr)e_i{\bigr
\rangle} {\bigl\langle} Y_{t_{k-2}}^n-Z_{t_{k-2}},\sigma
\bigl(Y^n_{t_{k-2}}\bigr)e_i{\bigr\rangle}.
\nonumber
\end{eqnarray}
This estimate will be used in Lemma~\ref{lemmaI5710}.

As for the term $I_7$,
\begin{eqnarray*}
I_7&=&-\frac{2}{\gamma}\int_{t_{k-1}}^{t_k}
\mu _n(t)\bigl|Y^n_t-Z_t\bigr|^2{
\bigl\langle} D\vph \bigl(Y^n_s\bigr),\sigma
\bigl(Y^n_s\bigr) \bigl(\d w_s-
\d{w}_s^n\bigr){\bigr\rangle}
\\
&=&-\frac{2}{\gamma}\int_{t_{k-1}}^{t_k}
\mu_n(t)\bigl|Y^n_{\bar
{t}_n}-Z_{\bar
{t}_n}\bigr|^2{
\bigl\langle}D\vph\bigl(Y^n_s\bigr),\sigma
\bigl(Y^n_s\bigr) \bigl(\d w_s-\d
{w}_s^n\bigr){\bigr\rangle}
\\
&&{}-\frac{4}{\gamma}\int_{t_{k-1}}^{t_k}
\mu_n(t)\int_{\bar
{t}_n}^t{\bigl\langle}
Y^n_s-Z_s,\sigma\bigl(Y^n_s
\bigr) \bigl(\d w_s-\d{w}_s^n\bigr){\bigr
\rangle}
\\
&&\hspace*{84pt}{}\times{\bigl\langle}D\vph\bigl(Y^n_t\bigr),\sigma
\bigl(Y^n_t\bigr) \bigl(\d w_t-\d
{w}_t^n\bigr){\bigr\rangle}
\\
&&{}-\frac{4}{\gamma}\int_{t_{k-1}}^{t_k}
\mu_n(t)\int_{\bar
{t}_n}^t{\bigl\langle}
Y^n_s-Z_s, \bigl(\sigma
\bigl(Y^n_s\bigr)-\sigma(Z_s) \bigr)
\dot{h}_s{\bigr\rangle}\,\d s
\\
&&\hspace*{84pt}{}\times{\bigl\langle}D\vph\bigl(Y^n_t\bigr),\sigma
\bigl(Y^n_t\bigr) \bigl(\d w_t-\d
{w}_t^n\bigr){\bigr\rangle}
\\
&&{}-\frac{4}{\gamma}\int_{t_{k-1}}^{t_k}
\mu_n(t)\int_{\bar
{t}_n}^t{\bigl\langle}
Y^n_s-Z_s,\tilde{b}\bigl(Y^n_s
\bigr)-b(Z_s){\bigr\rangle}\,\d s
\\
&&\hspace*{84pt}{}\times{\bigl\langle}D\vph\bigl(Y^n_t\bigr),\sigma
\bigl(Y^n_t\bigr) \bigl(\d w_t-\d
{w}_t^n\bigr){\bigr\rangle}
\\
&&{}-\frac{4}{\gamma}\int_{t_{k-1}}^{t_k}
\mu_n(t)\int_{\bar
{t}_n}^t{\bigl\langle}
Y^n_s-Z_s,\d\phi_s^n-
\d\psi_s{\bigr\rangle}
\\
&&\hspace*{78pt}{}\times{\bigl\langle}D\vph\bigl(Y^n_t\bigr),\sigma
\bigl(Y^n_t\bigr) \bigl(\d w_t-\d
{w}_t^n\bigr){\bigr\rangle}
\\
&&{}-\frac{2}{\gamma}\int_{t_{k-1}}^{t_k}
\mu_n(t)\int_{\bar{t}_n}^t \operatorname{tr}
\bigl(\sigma\sigma^*\bigl(Y^n_s\bigr)\bigr)\,\d s\times{
\bigl\langle}D\vph \bigl(Y^n_t\bigr),\sigma
\bigl(Y^n_t\bigr) \bigl(\d w_t-
\d{w}_t^n\bigr){\bigr\rangle}
\\
&=:&\sum_{i=1}^6 I_{7,i}.
\end{eqnarray*}
Notice that
\begin{eqnarray*}
I_{7,1}
&=&-\frac{2}{\gamma}\int_{t_{k-1}}^{t_k}\bigl(
\mu_n(t)-\mu _n(\bar {t}_n)
\bigr)\bigl|Y^n_{\bar{t}_n}-Z_{\bar{t}_n}\bigr|^2 {\bigl
\langle}D\vph\bigl(Y^n_t\bigr),\sigma
\bigl(Y^n_t\bigr) \bigl(\d w_t-
\d{w}_t^n\bigr){\bigr\rangle}
\\
&&{}-\frac{2}{\gamma}\int_{t_{k-1}}^{t_k}m_n(
\bar{t}_n){\bigl\langle }D\vph \bigl(Y^n_t
\bigr),\sigma\bigl(Y^n_t\bigr) \bigl(\d w_t-
\d{w}_t^n\bigr){\bigr\rangle}
\\
&=&\frac{4}{\gamma^2}\int_{t_{k-1}}^{t_k}\bigl|Y^n_{\bar{t}_n}-Z_{\bar
{t}_n}\bigr|^2
\int_{\bar{t}_n}^t\mu_n(s){\bigl\langle}D
\vph\bigl(Y^n_s\bigr),\sigma \bigl(Y^n_s
\bigr) \bigl(\d w_s-\d w_s^n\bigr){\bigr
\rangle}
\\
&&\hspace*{102pt}{}\times{\bigl\langle}D\varphi\bigl(Y^n_t\bigr),\sigma
\bigl(Y^n_t\bigr) \bigl(\d w_t-\d
{w}_t^n\bigr){\bigr\rangle}
\\
&&{}+\frac{4}{\gamma^2}\int_{t_{k-1}}^{t_k}\bigl|Y^n_{\bar{t}_n}-Z_{\bar
{t}_n}\bigr|^2
\int_{\bar{t}_n}^t\mu_n(s) \bigl({\bigl
\langle}D\vph\bigl(Y^n_s\bigr),\sigma
\bigl(Y^n_s\bigr){\bigr\rangle}\\
&&\hspace*{148pt}{} +{\bigl\langle} D
\vph(Z_s),\sigma(Z_s){\bigr\rangle}\bigr)
\dot{h}_s\,\d s
\\
&&\hspace*{115pt}{}\times{\bigl\langle}D\varphi\bigl(Y^n_t\bigr),\sigma
\bigl(Y^n_t\bigr) \bigl(\d w_t-\d
{w}_t^n\bigr){\bigr\rangle}
\\
&&{}+\frac{4}{\gamma^2}\int_{t_{k-1}}^{t_k}\bigl|Y^n_{\bar{t}_n}-Z_{\bar
{t}_n}\bigr|^2
\int_{\bar{t}_n}^t\mu_n(s) \bigl({\bigl
\langle}D\vph\bigl(Y^n_s\bigr),\tilde {b}
\bigl(Y^n_s\bigr){\bigr\rangle} +{\bigl\langle}D
\vph(Z_s),b(Z_s){\bigr\rangle}\bigr)\,\d s
\\
&&\hspace*{115pt}{}\times{\bigl\langle}D\vph\bigl(Y^n_t\bigr),\sigma
\bigl(Y^n_t\bigr) \bigl(\d w_t-\d
{w}_t^n\bigr){\bigr\rangle}
\\
&&{}+\frac{4}{\gamma^2}\int_{t_{k-1}}^{t_k}\bigl|Y^n_{\bar{t}_n}-Z_{\bar
{t}_n}\bigr|^2
\int_{\bar{t}_n}^t\mu_n(s) \bigl({\bigl
\langle}D\vph\bigl(Y^n_s\bigr),\d\phi
^n_s{\bigr\rangle} +{\bigl\langle} D
\vph(Z_s),\d\psi_s{\bigr\rangle}\bigr)
\\
&&\hspace*{114pt}{}\times{\bigl\langle}D\vph\bigl(Y^n_t\bigr),\sigma
\bigl(Y^n_t\bigr) \bigl(\d w_t-\d
{w}_t^n\bigr){\bigr\rangle}
\\
&&\hspace*{-1pt}{}+\frac{4}{\gamma^2}\hspace*{-1pt}\int_{t_{k-1}}^{t_k}\bigl|Y^n_{\bar{t}_n}-Z_{\bar
{t}_n}\bigr|^2\hspace*{-1pt}
\int_{\bar{t}_n}^t\mu_n(s) \biggl(
\frac{1}{2}\operatorname {tr}\bigl(D^2\vph \sigma\sigma^*\bigr)
\bigl(Y^n_s\bigr)-\frac{1}{\gamma}\bigl|D\vph\sigma
\bigl(Y^n_s\bigr)\bigr|^2 \biggr)\,\d s
\\
&&\hspace*{116pt}{}\times{\bigl\langle}D\vph\bigl(Y^n_t\bigr),\sigma
\bigl(Y^n_t\bigr) \bigl(\d w_t-\d
{w}_t^n\bigr){\bigr\rangle}
\\
&&\hspace*{94pt}{}-\frac{2}{\gamma}\int_{t_{k-1}}^{t_k}m_n(
\bar{t}_n){\bigl\langle }D\vph \bigl(Y^n_t
\bigr),\sigma\bigl(Y^n_t\bigr) \bigl(\d w_t-
\d{w}_t^n\bigr){\bigr\rangle}\\
&=:&\sum
_{i=1}^6 I_{7,1}^i.
\end{eqnarray*}
For the first term $I_{7,1}^1$,
\begin{eqnarray*}
I_{7,1}^1
&=&\frac{4}{\gamma^2}\int_{t_{k-1}}^{t_k}\bigl|Y^n_{\bar{t}_n}-Z_{\bar
{t}_n}\bigr|^2
\int_{\bar{t}_n}^t\bigl(\mu_n(s)-
\mu_n(\bar{t}_n)\bigr){\bigl\langle }D\vph
\bigl(Y^n_s\bigr),\sigma\bigl(Y^n_s
\bigr) \bigl(\d w_s-\d w_s^n\bigr){\bigr
\rangle}
\\
&&\hspace*{101pt}{}\times{\bigl\langle}D\vph\bigl(Y^n_t\bigr),\sigma
\bigl(Y^n_t\bigr) \bigl(\d w_t-\d
{w}_t^n\bigr){\bigr\rangle}
\\
&&{}+\frac{4}{\gamma^2}\int_{t_{k-1}}^{t_k}m_n(
\bar{t}_n)\int_{\bar
{t}_n}^t{\bigl\langle}D
\vph\bigl(Y^n_s\bigr),\sigma\bigl(Y^n_s
\bigr) \bigl(\d w_s-\d w_s^n\bigr){\bigr
\rangle }
\\
&&\hspace*{94pt}{}\times{\bigl\langle}D\vph\bigl(Y^n_t\bigr),\sigma
\bigl(Y^n_t\bigr) \bigl(\d w_t-\d
{w}_t^n\bigr){\bigr\rangle}.
\end{eqnarray*}
However, note that $\bE\sup_{t\in[0,T]}|Y^n_t-Z_t|^{4}<\infty$
by Lemma~\ref{holderZ} and Proposition~\ref{holderC},
and by Lemma~\ref{mulemma},
\begin{eqnarray*}
&&\bE\biggl|\int_{t_{k-1}}^{t_k}\int_{\bar{t}_n}^t
\bigl(\mu_n(s)-\mu_n(\bar {t}_n)\bigr){
\bigl\langle}D\vph\bigl(Y^n_s\bigr),\sigma
\bigl(Y^n_s\bigr) \bigl(\d w_s-\d
w_s^n\bigr){\bigr\rangle }
\\
&&\hspace*{100pt}\quad{}\times{\bigl\langle}D\vph\bigl(Y^n_t\bigr),\sigma
\bigl(Y^n_t\bigr) \bigl(\d w_t-\d
{w}_t^n\bigr){\bigr\rangle}\biggr|^2
\\
&&\qquad\leq\bE\biggl|\int_{t_{k-1}}^{t_k}\int_{\bar{t}_n}^t
\bigl(\mu_n(s)-\mu _n(\bar {t}_n)\bigr){
\bigl\langle}D\vph\bigl(Y^n_s\bigr),\sigma
\bigl(Y^n_s\bigr) \bigl(\d w_s-\d
w_s^n\bigr){\bigr\rangle }
\\
&&\hspace*{150pt}\qquad\quad{}\times{\bigl\langle}D\vph\bigl(Y^n_t\bigr),\sigma
\bigl(Y^n_t\bigr)\,\d w_t{\bigr
\rangle}\biggr|^2
\\
&&\qquad\quad{}+\bE\biggl|\int_{t_{k-1}}^{t_k}\int_{\bar{t}_n}^t
\bigl(\mu_n(s)-\mu_n(\bar {t}_n)\bigr){
\bigl\langle}D\vph\bigl(Y^n_s\bigr),\sigma
\bigl(Y^n_s\bigr)\,\d w_s{\bigr\rangle}
\\
&&\hspace*{122pt}\qquad\quad{}\times{\bigl\langle}D\vph\bigl(Y^n_t\bigr),\sigma
\bigl(Y^n_t\bigr)\,\d w^n_t{\bigr
\rangle}\biggr|^2
\\
&&\qquad\quad{}+\bE\biggl|\int_{t_{k-1}}^{t_k}\int_{\bar{t}_n}^t
\bigl(\mu_n(s)-\mu_n(\bar {t}_n)\bigr){
\bigl\langle}D\vph\bigl(Y^n_s\bigr),\sigma
\bigl(Y^n_s\bigr)\,\d w_s^n{
\bigr\rangle}
\\
&&\hspace*{124pt}\qquad\quad{}\times{\bigl\langle}D\vph\bigl(Y^n_t\bigr),\sigma
\bigl(Y^n_t\bigr)\,\d w^n_t{\bigr
\rangle}\biggr|^2
\\
&&\qquad\leq C\Delta^3+ \biggl[\bE\max_{t\in[t_{k-1},t_k]} \biggl(
\int_{\bar
{t}_n}^t\bigl(\mu_n(s)-
\mu_n(\bar{t}_n)\bigr){\bigl\langle}D\vph
\bigl(Y^n_s\bigr),\sigma \bigl(Y^n_s
\bigr)\,\d w_s{\bigr\rangle} \biggr)^4
\biggr]^{1/2}
\\
&&\qquad\quad{}\times \biggl[\bE \biggl(\int_{t_{k-1}}^{t_k}{\bigl
\langle}D\vph \bigl(Y^n_t\bigr),\sigma
\bigl(Y^n_t\bigr)\,\d w^n_t{\bigr
\rangle} \biggr)^4 \biggr]^{1/2},
\\
&&\qquad\quad{}+ \biggl[\bE\max_{t\in[t_{k-1},t_k]} \biggl(\int_{\bar{t}_n}^t
\bigl(\mu _n(s)-\mu _n(\bar{t}_n)\bigr){
\bigl\langle}D\vph\bigl(Y^n_s\bigr),\sigma
\bigl(Y^n_s\bigr)\,\d w^n_s{
\bigr\rangle } \biggr)^4 \biggr]^{1/2}
\\
&&\qquad\quad{}\times \biggl[\bE \biggl(\int_{t_{k-1}}^{t_k}{\bigl
\langle}D\vph \bigl(Y^n_t\bigr),\sigma
\bigl(Y^n_t\bigr)\,\d w^n_t{\bigr
\rangle} \biggr)^4 \biggr]^{1/2}\leq C\Delta^3.
\end{eqnarray*}
Similar to the term $I_{5,2,1}$,
\begin{eqnarray*}
&&\frac{4}{\gamma^2}\int_{t_{k-1}}^{t_k}m_n(
\bar{t}_n)\int_{\bar
{t}_n}^t{\bigl\langle}D
\vph\bigl(Y^n_s\bigr),\sigma\bigl(Y^n_s
\bigr) \bigl(\d w_s-\d w_s^n\bigr){\bigr
\rangle }
\\
&&\hspace*{70pt}\quad{}\times{\bigl\langle}D\vph\bigl(Y^n_t\bigr),\sigma
\bigl(Y^n_t\bigr)\,\d w^n_t{\bigr
\rangle}
\\
&&\qquad=\frac{4}{\gamma^2}\int_{t_{k-1}}^{t_k}m_n(
\bar{t}_n)\int_{\bar
{t}_n}^t{\bigl\langle}D
\vph\bigl(Y^n_s\bigr)-D\vph\bigl(Y^n_{\bar{t}_n}
\bigr),\sigma \bigl(Y^n_s\bigr) \bigl(\d
w_s-\d w_s^n\bigr){\bigr\rangle}
\\
&&\hspace*{83pt}\qquad\quad{}\times{\bigl\langle}D\vph\bigl(Y^n_t\bigr),\sigma
\bigl(Y^n_t\bigr)\,\d w^n_t{\bigr
\rangle}
\\
&&\qquad\quad{}+\frac{4}{\gamma^2}\int_{t_{k-1}}^{t_k}m_n(
\bar{t}_n)\int_{\bar
{t}_n}^t{\bigl\langle}D
\vph\bigl(Y^n_{\bar{t}_n}\bigr),\bigl(\sigma\bigl(Y^n_s
\bigr)-\sigma \bigl(Y^n_{\bar
{t}_n}\bigr)\bigr) \bigl(\d
w_s-\d w_s^n\bigr){\bigr\rangle}
\\
&&\hspace*{94pt}\qquad\quad{}\times{\bigl\langle}D\vph\bigl(Y^n_t\bigr),\sigma
\bigl(Y^n_t\bigr)\,\d w^n_t{\bigr
\rangle}
\\
&&\qquad\quad{}+\frac{4}{\gamma^2}\int_{t_{k-1}}^{t_k}m_n(
\bar{t}_n)\int_{\bar
{t}_n}^t{\bigl\langle}D
\vph\bigl(Y^n_{\bar{t}_n}\bigr),\sigma\bigl(Y^n_{\bar{t}_n}
\bigr) \bigl(\d w_s-\d w_s^n\bigr){\bigr
\rangle}
\\
&&\hspace*{94pt}\qquad\quad{}\times{\bigl\langle}D\vph\bigl(Y^n_t\bigr)-D\vph
\bigl(Y^n_{\bar{t}_n}\bigr),\sigma \bigl(Y^n_t
\bigr)\,\d w^n_t{\bigr\rangle}
\\
&&\qquad\quad{}-\frac{4}{\gamma^2}\int_{t_{k-1}}^{t_k}m_n(
\bar{t}_n)\int_{\bar
{t}_n}^t{\bigl\langle}D
\vph\bigl(Y^n_{\bar{t}_n}\bigr),\sigma\bigl(Y^n_{\bar{t}_n}
\bigr) \bigl(\d w_s-\d w_s^n\bigr){\bigr
\rangle}
\\
&&\hspace*{94pt}\qquad\quad{}\times{\bigl\langle}D\vph\bigl(Y^n_{\bar{t}_n}\bigr),\sigma
\bigl(Y^n_t\bigr)-\sigma \bigl(Y^n_{\bar
{t}_n}
\bigr){\bigr\rangle} \,\d w^n_t
\\
&&\qquad\quad{}-\frac{4}{\gamma^2}\int_{t_{k-1}}^{t_k}m_n(
\bar{t}_n)\int_{\bar
{t}_n}^t{\bigl\langle}D
\vph\bigl(Y^n_{\bar{t}_n}\bigr),\sigma\bigl(Y^n_{\bar{t}_n}
\bigr) \bigl(\d w_s-\d w_s^n\bigr){\bigr
\rangle}
\\
&&\hspace*{94pt}\qquad\quad{}\times{\bigl\langle}D\vph\bigl(Y^n_{\bar{t}_n}\bigr),\sigma
\bigl(Y^n_{\bar
{t}_n}\bigr){\bigr\rangle}\,\d w^n_t
\end{eqnarray*}
and
\begin{eqnarray*}
&&\biggl|\bE \biggl[\int_{t_{k-1}}^{t_k}m_n(
\bar{t}_n)\int_{\bar
{t}_n}^t{\bigl\langle}
D\vph \bigl(Y^n_{\bar{t}_n}\bigr),\sigma\bigl(Y^n_{\bar{t}_n}
\bigr) \bigl(\d w_s-\d w_s^n\bigr){\bigr
\rangle }
\\
&&\hspace*{102pt}\quad{}\times{\bigl\langle}D\vph\bigl(Y^n_{\bar{t}_n}\bigr),\sigma
\bigl(Y^n_{\bar{t}_n}\bigr)\,\d w^n_t{
\bigr\rangle} \biggr]\biggr|
\\
&&\qquad=\biggl|\bE \biggl[\bE \biggl(\int_{t_{k-1}}^{t_k}m_n(
\bar{t}_n)\int_{\bar
{t}_n}^t{\bigl\langle}D
\vph\bigl(Y^n_{\bar{t}_n}\bigr),\sigma\bigl(Y^n_{\bar{t}_n}
\bigr) \bigl(\d w_s-\d w_s^n\bigr){\bigr
\rangle}
\\
&&\hspace*{98pt}\qquad\quad{}\times{\bigl\langle}D\vph\bigl(Y^n_{\bar{t}_n}\bigr),\sigma
\bigl(Y^n_{\bar{t}_n}\bigr)\,\d w^n_t{
\bigr\rangle}\Big |\sF _{t_{k-2}} \biggr) \biggr]\biggr|
\\
&&\qquad\leq Ca(t_{k-2})\Delta.
\end{eqnarray*}
Summing up we get
\[
\bigl|\bE I_{7,1}^1\bigr|\leq Ca(t_{k-2})\Delta+C \Delta
^{3/2}.
\]
Next, we have by Proposition~\ref{phiBV} and Lemma~\ref{holderZ},
\begin{eqnarray*}
\bigl|\bE I_{7,1}^2\bigr| &\leq&C\bE \bigl(\bigl|Y^n_{{t}_{k-2}}-Z_{{t}_{k-2}}\bigr|^2|
\Delta w_{k-1}| \bigr)\int_{t_{k-2}}^{t_k}|
\dot{h}_s|\,\d s
\\
&\leq&C\Delta^{1/2}\int_{t_{k-2}}^{t_k}|
\dot{h}_s|\,\d s,
\\
\bigl|\bE I_{7,1}^i\bigr| &\leq& C\Delta^{3/2},\qquad i=3,5
\end{eqnarray*}
and
\[
\bigl|\bE I_{7,1}^4\bigr|\leq C\bE\bigl(G_k^1
\bigr),
\]
where
%
%e3.12 #&#
\begin{eqnarray}
\label{Gk1} G_k^1&:=&\bigl|Y^n_{{t}_{k-2}}-Z_{{t}_{k-2}}\bigr|^2
\bigl(\bigl|\phi^n\bigr|_{t_k}^{t_{k-2}\vee
0}+|\psi|_{t_k}^{t_{k-2}\vee0}
\bigr)|\Delta w_{k-1}|,
\nonumber
\\[-8pt]
\\[-8pt]
\nonumber
\Biggl|\bE\Biggl(\sum_k^{2^n}
G_k^1\Biggr)\Biggr|&\leq &C\Delta^{\theta/2}\qquad \forall \theta
\in(0,1).
\end{eqnarray}
For the term $I_{7,1}^6$,
\begin{eqnarray*}
I_{7,1}^6&=&-\frac{2}{\gamma}m_n(t_{k-2})
\int_{t_{k-1}}^{t_k}{\bigl\langle}D\vph
\bigl(Y^n_t\bigr),\sigma\bigl(Y^n_t
\bigr) \bigl(\d w_t-\d{w}_t^n\bigr){\bigr
\rangle}
\\
&=&-\frac{2}{\gamma}m_n(t_{k-2})\int
_{t_{k-1}}^{t_k}{\bigl\langle}D\vph \bigl(Y^n_t
\bigr)-D\vph \bigl(Y^n_{\bar{t}_n}\bigr),\sigma
\bigl(Y^n_t\bigr) \bigl(\d w_t-
\d{w}_t^n\bigr){\bigr\rangle}
\\
&&{}-\frac{2}{\gamma}m_n(t_{k-2})\int
_{t_{k-1}}^{t_k}{\bigl\langle}D\vph \bigl(Y^n_{\bar
{t}_n}
\bigr),\bigl(\sigma\bigl(Y^n_t\bigr)-\sigma
\bigl(Y^n_{\bar{t}_n}\bigr)\bigr) \bigl(\d w_t-\d
{w}_t^n\bigr){\bigr\rangle}
\\
&&{}-\frac{2}{\gamma}m_n(t_{k-2})\int
_{t_{k-1}}^{t_k}{\bigl\langle}D\vph \bigl(Y^n_{\bar
{t}_n}
\bigr),\sigma\bigl(Y^n_{\bar{t}_n}\bigr) \bigl(\d
w_t-\d{w}_t^n\bigr){\bigr\rangle}.
\end{eqnarray*}
Moreover,
\begin{eqnarray*}
&&\biggl|\bE \biggl[m_n(t_{k-2})\int_{t_{k-1}}^{t_k}{
\bigl\langle}D\vph \bigl(Y^n_t\bigr)-D\vph
\bigl(Y^n_{\bar{t}_n}\bigr),\sigma\bigl(Y^n_t
\bigr)\,\d{w}_t^n{\bigr\rangle} \biggr]\biggr|
\\
&&\qquad\leq C\bE \biggl[m_n(t_{k-2})\int_{t_{k-1}}^{t_k}\bigl|Y^n_t-Y^n_{\bar
{t}_n}\bigr|\bigl|
\d{w}_t^n\bigr| \biggr]
\\
&&\qquad=C\bE \biggl[m_n(t_{k-2})\int_{t_{k-1}}^{t_k}\biggl|
\int_{\bar
{t}_n}^t\sigma \bigl(Y^n_s
\bigr) \bigl(\d w_s-\d w^n_s\bigr)
\\
&&\hspace*{86pt}\qquad\quad{}+\int_{\bar{t}_n}^t\sigma\bigl(Y^n_s
\bigr)\dot{h}_s\,\d s+\int_{\bar
{t}_n}^t
\tilde {b}\bigl(Y^n_s\bigr)\,\d s+\phi^n_t-
\phi^n_{\bar{t}_n}\biggr|\bigl|\d{w}_t^n\bigr| \biggr]
\\
&&\qquad\leq\bE \biggl[m_n(t_{k-2})\int_{t_{k-1}}^{t_k}\biggl|
\int_{\bar
{t}_n}^t\bigl(\sigma \bigl(Y^n_s
\bigr)-\sigma\bigl(Y^n_{\bar{s}_n}\bigr)+\sigma
\bigl(Y^n_{\bar{s}_n}\bigr)\bigr) \bigl(\d w_s-\d
w^n_s\bigr)\biggr|\bigl|\d{w}_t^n\bigr|
\biggr]
\\
&&\qquad\quad{}+\biggl|\bE \biggl[m_n(t_{k-2})\int_{t_{k-1}}^{t_k}\biggl|
\int_{\bar
{t}_n}^t\sigma \bigl(Y^n_s
\bigr)\dot{h}_s\,\d s+\int_{\bar{t}_n}^t
\tilde{b}\bigl(Y^n_s\bigr)\,\d s+\phi ^n_{\bar
{t}_n}-
\phi^n_t\biggr|\bigl|\d{w}_t^n\bigr| \biggr]
\\
&&\qquad\leq\bE \biggl[m_n(t_{k-2})\max_{t\in[t_{k-1},t_k]}\biggl|
\int_{\bar
{t}_n}^t\bigl(\sigma\bigl(Y^n_s
\bigr)-\sigma\bigl(Y^n_{\bar{s}_n}\bigr)\bigr) \bigl(\d
w_s-\d w^n_s\bigr)\biggr|\int
_{t_{k-1}}^{t_k}\bigl|\d{w}_t^n\bigr|
\biggr]
\\
&&\qquad\quad{}+\bE \biggl[m_n(t_{k-2})\int_{t_{k-1}}^{t_k}\biggl|
\int_{\bar
{t}_n}^t\sigma \bigl(Y^n_{\bar{s}_n}
\bigr) \bigl(\d w_s-\d w^n_s\bigr)\biggr|\bigl|
\d{w}_t^n\bigr| \biggr]
\\
&&\qquad\quad{}+\bE \biggl[m_n(t_{k-2})\int_{t_{k-1}}^{t_k}\biggl|
\int_{\bar
{t}_n}^t\sigma \bigl(Y^n_s
\bigr)\dot{h}_s\,\d s+\int_{\bar{t}_n}^t
\tilde{b}\bigl(Y^n_s\bigr)\,\d s\biggr|\bigl|\d {w}_t^n\bigr|
\biggr]
\\
&&\qquad\quad{}+\bE \biggl[m_n(t_{k-2})\int_{t_{k-1}}^{t_k}
\max_{t\in
[t_{k-1},t_k]}\bigl|\phi ^n_{\bar{t}_n}-
\phi^n_t\bigr|\bigl|\d{w}_t^n\bigr| \biggr]
\\
&&\qquad\leq \bigl[\bE \bigl(m_n(t_{k-2})|\Delta
w_{k-1}| \bigr)^2 \bigr]^{1/2}
\\
&&\qquad\quad{}\times \biggl[\bE \biggl(\max_{t\in[t_{k-1},t_k]}\int_{t_{k-2}}^t
\bigl(\sigma \bigl(Y^n_s\bigr)-\sigma
\bigl(Y^n_{\bar{s}_n}\bigr)\bigr) \bigl(\d w_s-\d
w^n_s\bigr) \biggr)^2
\biggr]^{1/2}
\\
&&\qquad\hspace*{-1pt}\quad{}+\bE \biggl[m_n(t_{k-2})\hspace*{-1pt}\int_{t_{k-1}}^{t_k}\biggl|
\biggl(\hspace*{-1pt}\int_{t_{k-2}}^{t_{k-1}}+\hspace*{-1pt}\int_{t_{k-1}}^{t}
\biggr)\sigma\bigl(Y^n_{\bar{s}_n}\bigr) \bigl(\d
w_s-\d w^n_s\bigr)\biggr|\hspace*{-1pt}|\Delta
w_{k-1}|\Delta^{-1}\,\d t \biggr]
\\
&&\qquad\quad{}+\bE \biggl[m_n(t_{k-2})\int_{t_{k-1}}^{t_k}\biggl|
\int_{\bar
{t}_n}^t\sigma \bigl(Y^n_s
\bigr)\dot{h}_s\,\d s+\int_{\bar{t}_n}^t
\tilde{b}\bigl(Y^n_s\bigr)\,\d s\biggr|\bigl|\d {w}_t^n\bigr|
\biggr]
\\
&&\qquad\quad{}+\bE \biggl[m_n(t_{k-2})\int_{t_{k-1}}^{t_k}
\max_{t\in
[t_{k-1},t_k]}\bigl|\phi ^n_{\bar{t}_n}-
\phi^n_t\bigr|\bigl|\d{w}_t^n\bigr| \biggr]
\\
&&\qquad\leq C\Delta^{3/2}+Ca(t_{k-2})\Delta+C\bE \biggl[|\Delta
w_{k-1}|\int_{t_{k-2}}^{t_k}|
\dot{h}_s|\,\d s \biggr]+C\bE G_k^1,
\end{eqnarray*}
where $G_k^1$ is defined in (\ref{Gk1}) and
\[
\sum_{k=2}^{2^n}C\bE \biggl[\int
_{t_{k-2}}^{t_k}|\dot{h}_s|\,\d s|\Delta
w_{k-1}| \biggr]\leq C\Delta^{\theta/2}\int_0^1|
\dot{h}_s|\,\d s\qquad \forall \theta\in(0,1).
\]
Similarly,
%
%e3.13 #&#
\begin{eqnarray}
\label{Bk} &&\biggl|\bE \biggl(m_n(t_{k-2})\int
_{t_{k-1}}^{t_k}{\bigl\langle}D\vph \bigl(Y^n_{\bar
{t}_n}
\bigr),\sigma\bigl(Y^n_t\bigr)-\sigma
\bigl(Y^n_{\bar{t}_n}\bigr){\bigr\rangle}\,\d {w}_t^n
\biggr)\biggr|
\nonumber
\\
&&\qquad\leq C a(t_{k-2})\Delta+C\Delta^{3/2}+C\bE
G_k^1,
\nonumber
\\[-8pt]
\\[-8pt]
\nonumber
&&\bE \biggl(m_n(t_{k-2})\int_{t_{k-1}}^{t_k}{
\bigl\langle}D\vph\bigl(Y^n_{\bar
{t}_n}\bigr),\sigma
\bigl(Y^n_{\bar{t}_n}\bigr){\bigr\rangle}\bigl(\d
w_t-\d{w}_t^n\bigr) \biggr)\\
&&\qquad=:\bE
G^2_k,\nonumber
\end{eqnarray}
while
%
%e3.14 #&#
\begin{eqnarray}
\label{sumgk2} \Biggl|\sum_{k=2}^{2^n}\bE
G^2_k\Biggr| &\leq& C \Bigl(\bE\max_{1\leq k\leq{2^n}}m_n^2(t_k)
\Bigr)^{1/2}
\nonumber
\\
&&{}\times \biggl(\bE\max_{1\leq k\leq{2^n}}\biggl|\int_{0}^{t_k}{
\bigl\langle }D\vph \bigl(Y^n_{\bar
{t}_n}\bigr),\sigma
\bigl(Y^n_{\bar{t}_n}\bigr){\bigr\rangle}\bigl(\d
w_t-\d{w}_t^n\bigr)\biggr|^2
\biggr)^{1/2}
\\
&\leq& C\Delta^{\theta}\qquad \forall\theta\in(0,1).\nonumber
\end{eqnarray}
Here the last inequality follows since
\begin{eqnarray*}
&&\bE\max_{1\leq k\leq{2^n}}\biggl|\int_{0}^{t_k}{
\bigl\langle}D\vph \bigl(Y^n_{\bar
{t}_n}\bigr),\sigma
\bigl(Y^n_{\bar{t}_n}\bigr){\bigr\rangle}\bigl(\d
w_t-\d{w}_t^n\bigr)\biggr|^2
\\
&&\qquad=\bE\max_{1\leq k\leq{2^n}}\Biggl|\sum_{i=1}^{k-1}{
\bigl\langle}D\vph \bigl(Y^n_{t_{i-1}}\bigr),\sigma
\bigl(Y^n_{t_{i-1}}\bigr){\bigr\rangle}(w_{t_{i+1}}-w_{t_i})
\\
&&\hspace*{42pt}\qquad\quad{}-\sum_{i=1}^{k-1}{\bigl\langle}D\vph
\bigl(Y^n_{t_{i-1}}\bigr),\sigma \bigl(Y^n_{t_{i-1}}
\bigr){\bigr\rangle} (w_{t_{i}}-w_{t_{i-1}})\Biggr|^2
\\
&&\qquad=\bE\max_{1\leq k\leq{2^n}}\Biggl|\sum_{i=1}^{k-1}
\bigl({\bigl\langle}D\vph \bigl(Y^n_{t_{i-1}}\bigr),\sigma
\bigl(Y^n_{t_{i-1}}\bigr){\bigr\rangle}-{\bigl\langle}D\vph
\bigl(Y^n_{t_{i}}\bigr),\sigma \bigl(Y^n_{t_{i}}
\bigr){\bigr\rangle} \bigr) (w_{t_{i+1}}-w_{t_i})
\\
&&\hspace*{45pt}\qquad\quad{}-{\bigl\langle}D\vph(x),\sigma(x){\bigr\rangle}w_{t_1}+{\bigl
\langle}D\vph \bigl(Y^n_{t_{k-1}}\bigr),\sigma
\bigl(Y^n_{t_{k-1}}\bigr){\bigr\rangle}(w_{t_{k}}-w_{t_{k-1}})\Biggr|^2
\\
&&\qquad\leq C\sum_{i=1}^3
J_{7,1}^{6,i},
\end{eqnarray*}
where by Lemma~\ref{estimatevarphi} and the conditions $\vph,  \sigma
\in
\cC_b^2$,
\begin{eqnarray*}
J_{7,1}^{6,1}&:=&\bE \biggl(\max_{1\leq k\leq{2^n}}\biggl|
\int_{0}^{t_k} \bigl({\bigl\langle} D\vph
\bigl(Y^n_{\bar{t}_n}\bigr),\sigma\bigl(Y^n_{\bar{t}_n}
\bigr){\bigr\rangle}-{\bigl\langle }D\vph\bigl(Y^n_{\hat
{t}_n}
\bigr),\sigma\bigl(Y^n_{\hat{t}_n}\bigr){\bigr\rangle} \bigr)\,\d
w_t\biggr|^2 \biggr)
\leq C\Delta,
\\
J_{7,1}^{6,2}&:=&\bE\bigl(w_{t_1}^2
\bigr)=\Delta,
\\
J_{7,1}^{6,3}&:=&\bE \Bigl(\max_{1\leq k\leq{2^n}}\bigl|{
\bigl\langle}D\vph \bigl(Y^n_{t_{k-1}}\bigr),\sigma
\bigl(Y^n_{t_{k-1}}\bigr){\bigr\rangle}\bigr|^2|\Delta
w_k|^2 \Bigr)\leq C\Delta^\theta \qquad\theta\in(0,1).
\end{eqnarray*}

We then consider $I_{7,2}$.
\begin{eqnarray*}
I_{7,2}&=& -\frac{4}{\gamma}\int_{t_{k-1}}^{t_k}
\mu_n(t)\int_{\bar
{t}_n}^t{\bigl\langle}
Y^n_s-Z_s-\bigl(Y^n_{\bar{t}_n}-Z_{\bar{t}_n}
\bigr),\sigma\bigl(Y^n_s\bigr) \bigl(\d
w_s-\d {w}_s^n\bigr){\bigr\rangle}
\\
&&\hspace*{80pt}{}\times{\bigl\langle}D\vph\bigl(Y^n_t\bigr),\sigma
\bigl(Y^n_t\bigr) \bigl(\d w_t-\d
{w}_t^n\bigr){\bigr\rangle}
\\
&&{}-\frac{4}{\gamma}\int_{t_{k-1}}^{t_k}
\mu_n(t)\int_{\bar
{t}_n}^t{\bigl\langle}
Y^n_{\bar{t}_n}-Z_{\bar{t}_n},\sigma\bigl(Y^n_s
\bigr) \bigl(\d w_s-\d {w}_s^n\bigr){\bigr
\rangle}
\\
&&\hspace*{85pt}{}\times{\bigl\langle}D\vph\bigl(Y^n_t\bigr),\sigma
\bigl(Y^n_t\bigr) \bigl(\d w_t-\d
{w}_t^n\bigr){\bigr\rangle}
\\
&=:&I_{7,2,1}+I_{7,2,2}
\end{eqnarray*}
and
\begin{eqnarray*}
I_{7,2,1}&=&-\frac{4}{\gamma}\int_{t_{k-1}}^{t_k}
\mu_n(t)\int_{\bar
{t}_n}^t{\bigl
\langle}Y^n_s-Z_s-\bigl(Y^n_{\bar{t}_n}-Z_{\bar{t}_n}
\bigr),\sigma \bigl(Y^n_s\bigr) \bigl(\d
w_s-\d{w}_s^n\bigr){\bigr\rangle}
\\
&&\hspace*{81pt}{}\times{\bigl\langle}D\vph\bigl(Y^n_t\bigr),\sigma
\bigl(Y^n_t\bigr)\,\d w_t{\bigr\rangle}
\\
&&{}+\frac{4}{\gamma}\int_{t_{k-1}}^{t_k}
\mu_n(t)\int_{\bar
{t}_n}^t{\bigl\langle}
Y^n_s-Z_s-\bigl(Y^n_{\bar{t}_n}-Z_{\bar{t}_n}
\bigr),\sigma\bigl(Y^n_s\bigr) \bigl(\d
w_s-\d {w}_s^n\bigr){\bigr\rangle}
\\
&&\hspace*{84pt}{}\times{\bigl\langle}D\vph\bigl(Y^n_t\bigr),\sigma
\bigl(Y^n_t\bigr)\,\d w^n_t{\bigr
\rangle}
\\
&=:&I_{7,2,1}^{(1)}+I_{7,2,1}^{(2)}.
\end{eqnarray*}
However, note that $\bE(I_{7,2,1}^{(1)})=0$, and by using the Schwartz
and BDG inequalities, Lemmas \ref{holderZ} and~\ref{estimatevarphi},
\begin{eqnarray*}
&&\bigl|\bE I_{7,2,1}^{(2)}\bigr|
\\
&&\qquad\leq C \bigl(\bE|\Delta w_{k-1}|^2 \bigr)^{1/2}
\\
&&\qquad\quad{}\times \biggl(\bE\max_{t\in[t_{k-1},t_k]}\biggl|\int_{\bar
{t}_n}^t{
\bigl\langle} Y^n_s-Z_s-
\bigl(Y^n_{\bar{t}_n}-Z_{\bar{t}_n}\bigr),\sigma
\bigl(Y^n_s\bigr) \bigl(\d w_s-\d
{w}_s^n\bigr){\bigr\rangle}\biggr|^2
\biggr)^{1/2}
\\
&&\qquad\leq C\Delta^{3/2}.
\end{eqnarray*}

Now according to Lemma~\ref{lemmaI5710},
\[
\bE\bigl(I_{7,2,2}(t_k)+I_{5,2,1}^5(t_k)+I_{10}(t_k)
\bigr)\leq C\Delta^{3/2}.
\]

On the other hand, by the assumptions on $f$ and $\sigma, b$, as well
as Lemma~\ref{holderZ} and Proposition~\ref{phiBV},
%
%e3.15 #&#
\begin{eqnarray}
\label{Vk} |\bE I_{7,3}|&\leq&C\bE\int_{t_{k-1}}^{t_k}
\mu_n(t)\int_{\bar
{t}_n}^t\bigl|Y^n_s-Z_s\bigr|^2|
\dot{h}_s|\,\d s\nonumber\\
&&\hspace*{77pt}{}\times\bigl|D\vph\bigl(Y^n_t\bigr)\bigr|
\bigl\|\sigma \bigl(Y^n_t\bigr)\bigr\|\bigl|\d w^n_t\bigr|
\nonumber
\\[-8pt]
\\[-8pt]
\nonumber
&\leq&C\bE \Bigl(\sup_{t\in[t_{k-2},t_k]}\bigl|Y^n_t-Z_t\bigr|^2\bigl|
\Delta w_{k-1}\bigr| \Bigr)\int_{t_{k-2}}^{t_k}|
\dot{h}_s|\,\d s
\\
&\leq&C\Delta^{1/2}\int_{t_{k-2}}^{t_k}|
\dot{h}_s|\,\d s\nonumber
\end{eqnarray}
and
\begin{eqnarray*}
|\bE I_{7,4}|&\leq&\frac{4}{\gamma}\biggl|\bE\int_{t_{k-1}}^{t_k}
\mu _n(t)\int_{\bar{t}_n}^t\bigl|Y^n_s-Z_s\bigr|\bigl|
\tilde{b}\bigl(Y^n_s\bigr)-b(Z_s)\bigr|\,\d s
\\
&&\hspace*{80pt}{}\times{\bigl\langle}D\vph\bigl(Y^n_t\bigr),\sigma
\bigl(Y^n_t\bigr) \bigl(\d w_t-\d
{w}_t^n\bigr){\bigr\rangle}\biggr|
\\
&\leq&C\biggl|\bE\int_{t_{k-1}}^{t_k}\mu_n(t)\int
_{\bar
{t}_n}^t\bigl(\bigl|Y^n_s-Z_s\bigr|^2+1
\bigr)\,\d s\,\d w^n_t\biggr|
\\
&\leq&C\bE \biggl(\int_{t_{k-1}}^{t_k}\int
_{\bar{t}_n}^tm_n(s)\,\d s\bigl|\d
w^n_t\bigr| \biggr)\\
&&{}+C\Delta\bE\biggl( \int
_{t_{k-2}}^{t_k}\bigl|\d w^n_t\bigr|\biggr)
\\
&\leq&C\Delta^{3/2},
\\
I_{7,5}&=& -\frac{4}{\gamma}\int_{t_{k-1}}^{t_k}
\mu_n(t)\int_{\bar
{t}_n}^t{\bigl\langle}
Y^n_s-Z_s-Y^n_{\bar{t}_n}+Z_{\bar{t}_n},
\d\phi_s^n-\d\psi _s{\bigr\rangle}
\\
&&\hspace*{83pt}{}\times{\bigl\langle}D\vph\bigl(Y^n_t\bigr),\sigma
\bigl(Y^n_t\bigr) \bigl(\d w_t-\d
{w}_t^n\bigr){\bigr\rangle}
\\
&&{}-\frac{4}{\gamma}\int_{t_{k-1}}^{t_k}
\mu_n(t){\bigl\langle}Y^n_{\bar
{t}_n}-Z_{\bar
{t}_n},
\phi_t^n-\phi_{\bar{t}_n}^n-(
\psi_t-\psi_{\bar
{t}_n}){\bigr\rangle}
\\
&&\hspace*{64pt}{}\times{\bigl\langle}D\vph\bigl(Y^n_t\bigr),\sigma
\bigl(Y^n_t\bigr) \bigl(\d w_t-\d
{w}_t^n\bigr){\bigr\rangle}
\\
&=:&I_{7,5,1}+I_{7,5,2}.
\end{eqnarray*}
By using (\ref{holderestimatevarphi}),
\begin{eqnarray}
&&|\bE I_{7,5,1}|
\nonumber\\
&&\qquad\leq C\bE\int_{t_{k-1}}^{t_k}\mu_n(t)\sup
_{s\in
[t_{k-2},t]}\bigl|Y^n_s-Z_s-Y^n_{t_{k-2}}+Z_{t_{k-2}}\bigr|
\nonumber\\
&&\hspace*{72pt}{}\times\bigl(\bigl|\phi ^n\bigr|_{t}^{t_{k-2}}+|\psi|_{t}^{t_{k-2}}
\bigr) \bigl|\d w^n_t\bigr|
\nonumber\\
&&\qquad\leq C\bE \biggl(\bigl(\bigl\|Y^n\bigr\|_{[t_{k-2},t_k]}+\|Z
\|_{[t_{k-2},t_k]}\bigr) \bigl(\bigl|\phi ^n\bigr|_{t_k}^{t_{k-2}}+|
\psi|_{t_k}^{t_{k-2}}\bigr)\int_{t_{k-1}}^{t_k}\bigl|
\d w^n_t\bigr| \biggr)
\nonumber\\
&&\qquad\leq C\Delta^{3/2},\nonumber
\\
\label{Gk}
&& |I_{7,5,2}|\leq C\bigl|Y^n_{t_{k-2}}-Z_{t_{k-2}}\bigr|
\nonumber\\
&&\hspace*{20pt}\qquad{}\times \biggl(\max_{t\in
[t_{k-1},t_k]}\biggl|\int_{t_{k-1}}^t{
\bigl\langle}D\vph\bigl(Y^n_t\bigr),\sigma
\bigl(Y^n_t\bigr)\,\d w_t{\bigr\rangle} \biggr|+|
\Delta w_{k-1}| \biggr)
\\
&&\hspace*{20pt}\qquad{}\times\bigl(\bigl|\phi^n\bigr|_{t_{k-2}}^{t_k}+|
\psi|_{t_{k-2}}^{t_k}\bigr)=:G_k^3,\nonumber
\end{eqnarray}
while according to Lemma~\ref{holderZ} and Proposition~\ref{phiBV},
\begin{eqnarray*}
&&\sum_{k=1}^{2^n}\bE G_k^3\\
&&\qquad\leq C\max_{1\leq k\leq2^n} \biggl(\bE \biggl[\bigl|Y^n_{t_{k-2}}-Z_{t_{k-2}}\bigr|^2\\
&&\qquad\hspace*{58pt}{}\times\biggl(\max_{t\in[t_{k-1},t_k]}\biggl|\int_{t_{k-1}}^t{
\bigl\langle}D\vph\bigl(Y^n_t\bigr),\sigma
\bigl(Y^n_t\bigr)\,\d w_t{\bigr\rangle
}\biggr|^2+|\Delta w_{k-1}|^2\biggr) \biggr]
\biggr)^{1/2}
\\
&&\qquad\quad{}\times \bigl(\bE\bigl(\bigl|\phi^n\bigr|_{1}^{0}+|
\psi|_{1}^{0}\bigr)^2 \bigr)^{1/2}
\\
&&\qquad\leq C\Delta^{1/2}.
\end{eqnarray*}
Also by the boundedness of $\sigma$ we have
\begin{eqnarray*}
|\bE I_{7,6}|&\leq& \frac{\gamma}{2}\biggl|\bE\int_{t_{k-1}}^{t_k}
\mu _n(t)\int_{\bar{t}_n}^t\operatorname{tr}
\bigl(\sigma\sigma^*\bigr) \bigl(Y^n_s\bigr)\,\d s\times{
\bigl\langle }D\vph \bigl(Y^n_t\bigr),\sigma
\bigl(Y^n_t\bigr) \bigl(\d w_t-
\d{w}_t^n\bigr){\bigr\rangle}\biggr|
\\
&\leq& C\Delta^{3/2}.
\end{eqnarray*}

Hence by applying all the above estimates to (\ref{akminus}),
\begin{eqnarray*}
a(t_k)-a(t_{k-1})&\leq&Ca(t_{k-2})\Delta+C
\Delta^{3/2}+b_k+\bE\int_{t_{k-1}}^{t_k}m_n(t)|
\dot{h}_t|\,\d t
\\
&\leq&C\bE \bigl(m_n(t_{k-2})-m_n(t_{k-1})+m_n(t_{k-1})
\bigr)\Delta +C\Delta ^{3/2}
\\
&&{}+Cb_k+\bE\int_{t_{k-1}}^{t_k}m_n(t)|
\dot{h}_t|\,\d t
\\
&\leq&Ca(t_{k-1})\Delta+C\Delta^{3/2}+Cb_k+\bE
\int_{t_{k-1}}^{t_k}m_n(t)|
\dot{h}_t|\,\d t,
\end{eqnarray*}
where
\begin{eqnarray*}
b_k&:=&G_k+G_k^1+G_k^2+G^3_k+A_k^n+C
\Delta^{1/2}\int_{t_{k-2}}^{t_k}|
\dot{h}_t|\,\d t,
\end{eqnarray*}
and hence
\begin{eqnarray*}
\sum_{k=1}^{2^n} b_k&=& \bE
\Biggl[\sum_{k=1}^{2^n}\biggl(A_k^n+G_k+G_k^1+G_k^2+G^3_k+C
\Delta^{1/2}\int_{t_{k-2}}^{t_k}|
\dot{h}_t|\,\d t\biggr) \Biggr]
\\
&\leq& C\biggl[\Delta^{\theta/2}+\sup_{2\leq k\leq2^n}\biggl(\int
_{t_{k-2}}^{t_k}|\dot{h}_s|^2\,\d s
\biggr)^{1/2}\biggr],\qquad \theta\in(0,1),
\end{eqnarray*}
and $A_k^n, G_k, G_k^1, G_k^2, G^3_k$, are defined in (\ref{Ak}),
(\ref
{defGk}), (\ref{Gk1}), (\ref{Bk}) and (\ref{Gk}), respectively.

Therefore according to Bihari's inequality, by denoting $h_k:=e^{C\int
_{t_{k-1}}^{t_k}|\dot{h}_s|\,\d s}$, we have
\begin{eqnarray*}
a(t_k)&\leq& \bigl[a(t_{k-1}) (1+C\Delta)+C
\Delta^{3/2}+b_k \bigr]h_k
\\
&\leq&a(t_{k-1})h_k(1+C\Delta)+C\Delta^{3/2}h_k+b_kh_k
\\
&\leq&\bigl[\bigl(a(t_{k-2}) (1+C\Delta)+C\Delta ^{3/2}+b_{k-1}
\bigr)h_{k-1}\bigr]h_k(1+C\Delta)\\
&&{}+C\Delta^{3/2}h_k+b_kh_k
\\
&\leq&\cdots
\\
&\leq&(1+C\Delta)^kh_k\cdots h_1
a(t_0)+\sum_{i=0}^{k-1}(1+C
\Delta)^i h_k\cdots h_{k-i}\bigl(C
\Delta^{3/2}+b_{k-i}\bigr)
\\
&\leq& e^{CT}e^{C\int_0^T|\dot{h}_t|\,\d t}C\biggl[\Delta^{\theta/2}+\sup
_{1\leq k\leq2^n}\biggl(\int_{t_{k-2}}^{t_k}|
\dot{h}_s|^2\,\d s\biggr)^{1/2}\biggr], \qquad\theta
\in(0,1),
\end{eqnarray*}
and we obtain the desired result.
\end{pf}

%le3.10 #&#
\begin{lemma}\label{Aklemma}
\[
\Biggl|\bE \Biggl(\sum_{i=1}^{2^n}A_i^n
\Biggr)\Biggr|\leq C\biggl[\Delta^{1/2}+\sup_{2\leq
k\leq2^n}\biggl(
\int_{t_{k-2}}^{t_k}|\dot{h}_s|^2
\,\d s\biggr)^{1/2}\biggr].
\]
\end{lemma}

\begin{pf}
Set
\begin{eqnarray*}
\zeta^n_t:=\int_{\bar{t}_n}^t
\bigl[\bigl(\sigma\bigl(Y^n_s\bigr)-\sigma
\bigl(Y^n_{\bar
{t}_n}\bigr)\bigr) \bigl(\d w_s-\d
w^n_s\bigr)+\sigma\bigl(Y^n_s
\bigr)\dot{h}_s\,\d s +\tilde{b}\bigl(Y^n_s
\bigr)\,\d s \bigr].
\end{eqnarray*}
We have for any $2\leq k\leq2^n$,
\begin{eqnarray*}
&&\bigl|\bE \bigl(A_{k}^n+A_{k+1}^n
\bigr)\bigr|
\\
&&\qquad=\biggl|\bE \biggl(\int_{t_{k-1}}^{t_{k+1}}\mu_n(
\bar{t}_n){\biggl\langle }Y^n_{\bar
{t}_n}-Z_{\bar{t}_n},
\bigl(\sigma\bigl(Y^n_t\bigr)-\sigma\bigl(Y^n_{\bar{t}_n}
\bigr)\bigr)\dot {w}^n_t-\frac{1}{2}(\nabla\sigma)
\sigma\bigl(Y^n_{\bar{t}_n}\bigr){\biggr\rangle }\,\d t
\biggr)\biggr|
\\
&&\qquad=\biggl|\bE \biggl(\int_{t_{k-1}}^{t_{k}}
\mu_n({t}_{k-2}){\biggl\langle }Y^n_{\bar
{t}_n}-Z_{\bar{t}_n},
\bigl(\sigma\bigl(Y^n_t\bigr)-\sigma\bigl(Y^n_{\bar{t}_n}
\bigr)\bigr)\dot {w}^n_t-\frac{1}{2}(\nabla\sigma)
\sigma\bigl(Y^n_{\bar{t}_n}\bigr){\biggr\rangle }\,\d t
\biggr)
\\
&&\qquad\quad{}+\bE \biggl(\int_{t_{k}}^{t_{k+1}} \bigl(
\mu_n({t}_{k-1})-\mu _n({t}_{k-2})
\bigr){\biggl\langle}Y^n_{\bar{t}_n}-Z_{\bar{t}_n},\bigl(
\sigma \bigl(Y^n_t\bigr)-\sigma \bigl(Y^n_{\bar{t}_n}
\bigr)\bigr)\dot{w}^n_t\\
&&\hspace*{253pt}{}-\frac{1}{2}(\nabla\sigma)
\sigma \bigl(Y^n_{\bar
{t}_n}\bigr){\biggr\rangle}\,\d t
\biggr)
\\
&&\qquad\quad{}+\bE \biggl(\int_{t_{k}}^{t_{k+1}}
\mu_n({t}_{k-2}){\biggl\langle }Y^n_{\bar
{t}_n}-Z_{\bar{t}_n},
\bigl(\sigma\bigl(Y^n_t\bigr)-\sigma\bigl(Y^n_{\bar{t}_n}
\bigr)\bigr)\dot {w}^n_t\\
&&\hspace*{195pt}{}-\frac{1}{2}(\nabla\sigma)
\sigma\bigl(Y^n_{\bar{t}_n}\bigr){\biggr\rangle }\,\d t
\biggr)\biggr|
\\
&&\qquad\leq \biggl|\bE \biggl(\int_{t_{k-1}}^{t_{k+1}}
\mu_n({t}_{k-2}){\biggl\langle} Y^n_{\bar
{t}_n}-Z_{\bar{t}_n},
\bigl(\sigma\bigl(Y^n_t\bigr)-\sigma\bigl(Y^n_{\bar{t}_n}
\bigr)\bigr)\dot {w}^n_t\\
&&\hspace*{181pt}{}-\frac{1}{2}(\nabla\sigma)
\sigma\bigl(Y^n_{\bar{t}_n}\bigr){\biggr\rangle }\,\d t
\biggr)\biggr|+C\Delta^{3/2}.
\end{eqnarray*}
Thus by continuing this procedure we get
\begin{eqnarray*}
&&\Biggl|\bE \Biggl(\sum_{i=1}^{2^n}A_i^n
\Biggr)\Biggr|
\\
&&\qquad\leq C2^n\Delta^{3/2}+\biggl|\bE \biggl(\int
_0^{1}\mu_n({t}_1){
\biggl\langle }Y^n_{\bar
{t}_n}-Z_{\bar{t}_n},\bigl(\sigma
\bigl(Y^n_t\bigr)-\sigma\bigl(Y^n_{\bar{t}_n}
\bigr)\bigr)\dot {w}^n_t\\
&&\hspace*{214pt}{}-\frac{1}{2}(\nabla\sigma)
\sigma\bigl(Y^n_{\bar{t}_n}\bigr){\biggr\rangle }\,\d t
\biggr)\biggr|
\\
&&\qquad\leq C\Delta^{1/2}+ \Bigl[\bE \Bigl(\mu_n^2(t_1)
\sup_{0\leq t\leq
1}\bigl|Y^n_{t}-Z_{t}\bigr|^2
\Bigr) \Bigr]^{1/2}
\\
&&\qquad\quad\times{} \biggl[\bE \biggl(\int_0^{1} \biggl(
\bigl(\sigma\bigl(Y^n_t\bigr)-\sigma
\bigl(Y^n_{\bar
{t}_n}\bigr)\bigr)\dot{w}^n_t-
\frac{1}{2}(\nabla\sigma)\sigma\bigl(Y^n_{\bar
{t}_n}\bigr)
\biggr)\,\d t \biggr)^2 \biggr]^{1/2}.
\end{eqnarray*}
Note that $\sigma\in\mathcal{C}_b^2$,
\begin{eqnarray*}
&&\bigl|\sigma\bigl(Y_t^n\bigr)-\sigma\bigl(Y_{\bar{t}_n}^n
\bigr)-(\nabla\sigma) \bigl(Y^n_{\bar
{t}_n}\bigr)
\bigl(Y^n_t-Y^n_{\bar{t}_n}\bigr)\bigr| \leq
C\bigl|Y^n_t-Y^n_{\bar{t}_n}\bigr|^2,
\\
&& Y^n_t-Y^n_{\bar{t}_n}=\int
_{\bar{t}_n}^t\sigma\bigl(Y^n_s
\bigr) \bigl(\d w_s-\d w^n_s\bigr)+\int
_{\bar{t}_n}^t\sigma\bigl(Y^n_s
\bigr)\dot{h}_s\,\d s\\
&&\hspace*{28pt}\qquad{} +\int_{\bar{t}_n}^t
\tilde{b}\bigl(Y^n_s\bigr)\,\d s+\phi^n_t-
\phi^n_{\bar{t}_n}.
\end{eqnarray*}
Then
\begin{eqnarray*}
&&\bE \biggl[\sup_{2\leq k\leq2^n}\biggl|\int_0^{t_k}
\biggl(\bigl(\sigma \bigl(Y^n_t\bigr)-\sigma
\bigl(Y^n_{\bar{t}_n}\bigr)\bigr)\dot{w}^n_t-
\frac{1}{2}(\nabla\sigma)\sigma \bigl(Y^n_{\bar
{t}_n}\bigr)
\biggr)\,\d t\biggr|^2 \biggr]
\\
&&\qquad\leq C\bE \biggl[ \biggl(\int_0^{1}\bigl|Y^n_t-Y^n_{\bar{t}_n}\bigr|^2\bigl|
\dot {w}^n_t\bigr|\,\d t \biggr)^2 \biggr]
\\
&&\qquad\quad{}+C\bE \biggl[\sup_{1\leq k\leq2^n} \biggl(\int_0^{t_k}
\nabla\sigma \bigl(Y^n_{\bar{t}_n}\bigr) \bigl(
\zeta_t+\phi^n_t-\phi^n_{\bar{t}_n}
\bigr)\dot {w}^n_t\,\d t \biggr)^2 \biggr]
\\
&&\qquad\quad{}+C\bE \biggl[\sup_{1\leq k\leq2^n} \biggl(\int_0^{t_k}
\biggl((\nabla \sigma )\sigma\bigl(Y^n_{\bar{t}_n}\bigr)\int
_{\bar{t}_n}^t\,\d w_s\dot
{w}^n_t-(\nabla \sigma)\sigma\bigl(Y^n_{\bar{t}_n}
\bigr) \biggr)\,\d t \biggr)^2 \biggr]
\\
&&\qquad\quad{}+C\bE \biggl[\sup_{1\leq k\leq2^n} \biggl(\int_0^{t_k}
\biggl((\nabla \sigma )\sigma\bigl(Y^n_{\bar{t}_n}\bigr)\int
_{\bar{t}_n}^t\,\d w^n_s\dot
{w}^n_t-\frac
{1}{2}(\nabla\sigma)\sigma
\bigl(Y^n_{\bar{t}_n}\bigr) \biggr)\,\d t \biggr)^2
\biggr]
\\
&&\qquad=: \sum_{\alpha=1}^4 T_\alpha.
\end{eqnarray*}

Note that by (\ref{holderestimatevarphi}),
\begin{eqnarray*}
T_1&=&C\bE \biggl[ \biggl(\int_0^{1}\bigl|Y^n_t-Y^n_{\bar{t}_n}\bigr|^2\bigl|
\dot {w}^n_t\bigr|\,\d t \biggr)^2 \biggr]
\\
&\leq&C\bE \Biggl[ \Biggl(\sum_{i=1}^{2^n}
\int_{t_{i-1}}^{t_i}\bigl|Y^n_t-Y^n_{\bar
{t}_n}\bigr|^2
\bigl|\dot{w}^n_t\bigr|\,\d t \Biggr)^2 \Biggr]
\\
&\leq&C2^{2n}\max_{1\leq i\leq2^n}\bE \Bigl(\sup
_{t\in
[t_{i-1},t_i]}\bigl|Y^n_t-Y^n_{t_{i-2}\vee0}\bigr|^4|
\Delta w_{i-1}|^2 \Bigr) \leq C\Delta.
\end{eqnarray*}
By Lemma~\ref{estimatevarphi} and Proposition~\ref{phiBV},
\begin{eqnarray*}
T_2&=&C\bE \biggl[\sup_{1\leq k\leq2^n} \biggl(\int
_0^{t_k}\nabla \sigma \bigl(Y^n_{\bar{t}_n}
\bigr) \bigl(\zeta_t+\phi^n_t-
\phi^n_{\bar{t}_n}\bigr)\dot {w}^n_t\,\d t
\biggr)^2 \biggr]
\\
&\leq& C\bE \biggl[ \biggl(\int_0^{1}\bigl|\nabla
\sigma\bigl(Y^n_{\bar{t}_n}\bigr)\zeta _t\dot
{w}^n_t\bigr|\,\d t \biggr)^2 \biggr]
\\
&&{}+C\bE \biggl[\sup_{1\leq k\leq2^n} \biggl(\int_0^{t_k}
\nabla\sigma \bigl(Y^n_{\bar{t}_n}\bigr) \bigl(
\phi^n_t-\phi^n_{\bar{t}_n}\bigr)
\dot{w}^n_t\,\d t \biggr)^2 \biggr]
\\
&\leq&C \biggl\{\bE\int_0^1|
\zeta_t|^{4}\,\d t \biggr\}^{1/2} \biggl\{\bE \int
_0^1\bigl|\dot{w}^n_t\bigr|^{4}
\,\d t \biggr\}^{1/2}
\\
&&{}+C\bE \Biggl[\sup_{1\leq k\leq2^n}\Biggl(\sum
_{i=1}^k\int_{t_{i-1}}^{t_i}
\nabla\sigma\bigl(Y^n_{\bar{t}_n}\bigr) \bigl(
\phi^n_t-\phi ^n_{\bar
{t}_n}\bigr)
\dot{w}^n_t\,\d t\Biggr)^2 \Biggr]
\\
&\leq&C\Delta \biggl(\Delta+\sup_{2\leq k\leq2^n}\biggl(\int
_{t_{k-2}}^{t_k}|\dot{h}_s|^2\,\d s
\biggr) \biggr)\Delta^{-1}+C\bE \Bigl[\Bigl(\bigl|\phi ^n\bigr|_{1}
\sup_{1\leq i\leq2^n}|\Delta w_{i}|\Bigr)^2 \Bigr]
\\
&\leq&C\biggl(\Delta+\sup_{1\leq k\leq2^n}\int_{t_{k-2}}^{t_k}|
\dot {h}_s|^2\,\d s\biggr)\\
&&{}+ \Bigl[\bE\Bigl(\sup
_{1\leq k\leq2^n}|\Delta w_{k}|^{2p}\Bigr)
\Bigr]^{1/p} \bigl[\bE\bigl(\bigl|\phi^n\bigr|_{1}
\bigr)^{2q} \bigr]^{1/q}
\\
&\leq&C\biggl(\Delta^{1-1/p}+\sup_{1\leq k\leq2^n}\int
_{t_{k-2}}^{t_k}|\dot {h}_s|^2\,\d
s\biggr) \qquad\forall p, q>1, 1/p+1/q=1.
\end{eqnarray*}
Note that $T_3=T_3^1+T_3^2$, where
\begin{eqnarray*}
T_3^1&:=&C\bE\sup_{1\leq k\leq2^n} \biggl(\int
_0^{t_k} \biggl((\nabla \sigma )\sigma
\bigl(Y^n_{\bar{t}_n}\bigr)\int_{\bar{t}_n}^{\hat{t}_n}
\,\d w_s\dot {w}^n_t-(\nabla\sigma)\sigma
\bigl(Y^n_{\bar{t}_n}\bigr) \biggr)\,\d t \biggr)^2
\\
&\leq&C\bE \Biggl[\sum_{i=1}^{2^n} \bigl((
\nabla\sigma)\sigma \bigl(Y^n_{{t}_{i-1}\vee0}\bigr)
\bigr)^2\bigl(|\Delta w_i|^2-\Delta
\bigr)^2
\\
&&\hspace*{20pt}{}+2\sum_{i<j} (\nabla\sigma)\sigma
\bigl(Y^n_{{t}_{i-1}\vee0}\bigr) (\nabla \sigma )\sigma
\bigl(Y^n_{{t}_{j-1}\vee0}\bigr) \bigl(|\Delta w_i|^2-
\Delta\bigr) \bigl(|\Delta w_j|^2-\Delta\bigr) \Biggr]
\\
&\leq&C\bE \Biggl[\sum_{i=1}^{2^n} \bigl((
\nabla\sigma)\sigma \bigl(Y^n_{{t}_{i-1}\vee0}\bigr)
\bigr)^2\bigl(|\Delta w_i|^2-\Delta
\bigr)^2 \Biggr]
\\
&\leq&C2^n\Delta^2\leq C\Delta,
\\
T_3^2&:=&C\bE\sup_{1\leq k\leq2^n} \biggl(\int
_0^{t_k}(\nabla\sigma )\sigma
\bigl(Y^n_{\bar{t}_n}\bigr)\int_{\hat{t}_n}^t
\,\d w_s\dot{w}^n_t\,\d t \biggr)^2
\\
&\leq&C\bE \biggl(\int_0^1(\nabla\sigma)
\sigma\bigl(Y^n_{\bar
{t}_n}\bigr)\frac{\hat
{t}_n+\Delta-t}{\Delta}(w_{\hat{t}_n}-w_{\bar{t}_n})
\,\d w_t \biggr)^2
\\
&\leq&C\bE \biggl(\int_0^1\bigl|(\nabla\sigma)
\sigma\bigl(Y^n_{\bar
{t}_n}\bigr)\bigr|^2|w_{\hat
{t}_n}-w_{\bar{t}_n}|^2
\,\d t \biggr)
\\
&\leq& C\Delta.
\end{eqnarray*}
Also, $T_4=T_4^1+T_4^2$, where
\begin{eqnarray*}
T_4^1&:=&C\bE \biggl[\sup_{1\leq k\leq2^n}
\biggl(\int_0^{t_k}(\nabla \sigma )\sigma
\bigl(Y^n_{\bar{t}_n}\bigr)\int_{\bar{t}_n}^{\hat{t}_n}
\,\d w^n_s\dot {w}_t^n\,\d t
\biggr)^2 \biggr]
\\
&\leq&C\bE \biggl(\int_0^1(\nabla\sigma)
\sigma\bigl(Y^n_{\bar
{t}_n}\bigr) \biggl(\int
_{\bar
{t}_n}^{\hat{t}_n}\,\d w^n_s
\biggr)\,\d w_t \biggr)^2
\\
&\leq&C\bE \biggl(\int_0^1\biggl|\int
_{\bar{t}_n}^{\hat{t}_n}\,\d w^n_s\biggr|^2
\,\d t \biggr)
\\
&=&C\sum_{i=0}^{2^n-1}\bE \biggl(\int
_{t_i}^{t_{i+1}}\biggl|\int_{\bar
{t}_n}^{\hat{t}_n}
\,\d w^n_s\biggr|^2\,\d t \biggr)
\\
&\leq& C\Delta,
\\
T_4^2&:=&C\bE\sup_{1\leq k\leq2^n} \biggl(\int
_0^{t_k} \biggl((\nabla \sigma )\sigma
\bigl(Y^n_{\bar{t}_n}\bigr)\int_{\hat{t}_n}^t
\,\d w^n_s\dot {w}^n_t-
\frac
{1}{2}(\nabla\sigma)\sigma\bigl(Y^n_{\bar{t}_n}\bigr)
\biggr)\,\d t \biggr)^2
\\
&\leq&C\bE\sup_{1\leq k\leq2^n} \Biggl(\sum
_{i=0}^{k-2}(\nabla \sigma )\sigma
\bigl(Y^n_{t_i}\bigr) \biggl(\biggl(\Delta^{-2}
\int_{t_{i+1}}^{t_{i+2}}\int_{t_{i+1}}^t
\,\d s\,\d t\biggr) (w_{t_{i+1}}-w_{t_{i}})^2\\
&&\hspace*{265pt}{}-
\frac{1}{2}\Delta \biggr) \Biggr)^2
\\
&\leq&C\Delta.
\end{eqnarray*}
Summing these estimates we get
\begin{eqnarray*}
\Biggl|\bE\sup_{1\leq k\leq2^n} \Biggl(\sum_{i=1}^kA_i^n
\Biggr)\Biggr|&\leq &C\biggl[\Delta ^{1/2}+\sup_{1\leq k\leq2^n}
\biggl(\int_{t_{k-2}}^{t_k}|\dot{h}_s|^2
\,\d s\biggr)^{1/2}\biggr].
\end{eqnarray*}
\upqed\end{pf}

%le3.11 #&#
\begin{lemma}\label{lemmaI51}
$ |\bE (\sum_{i=0}^{2^n}{I}_{5,1}(t_i) )|\leq C\Delta
^{\theta
/2}, \forall\theta\in(0,1)$.
\end{lemma}

\begin{pf}
Since
\begin{eqnarray*}
I_{5,1}=2\int_{t_{k-1}}^{t_k}
\mu_n(\bar{t}_n){\bigl\langle }Y_t^n-Z_t,
\sigma \bigl(Y^n_{\bar{t}_n}\bigr){\bigr\rangle}\bigl(\d
w_t-\d w^n_t\bigr)+\int_{t_{k-1}}^{t_k}
\mu _n(\bar {t}_n)\operatorname{tr}\bigl(\sigma\sigma^*\bigr)
\bigl(Y^n_{\bar{t}_n}\bigr)\,\d t,
\end{eqnarray*}
it is trivial to see that $\bE(I_{5,1}-\tilde{I}_{5,1})=0$, where
\begin{eqnarray*}
\tilde{I}_{5,1}(t_k)&:=&2\int_{t_{k-1}}^{t_k}
\mu_n(\bar {t}_n){\bigl\langle} Y_t^n-Z_t-
\bigl(Y^n_{\bar{t}_n}-Z_{\bar{t}_n}\bigr),\sigma
\bigl(Y^n_{\bar
{t}_n}\bigr) \bigl(\d w_t-\d
w^n_t\bigr){\bigr\rangle}
\\
&&{}+\int_{t_{k-1}}^{t_k}\mu_n(
\bar{t}_n)\operatorname{tr}\bigl(\sigma\sigma ^*\bigr)
\bigl(Y^n_{\bar{t}_n}\bigr)\,\d t.
\end{eqnarray*}

Note that
\begin{eqnarray*}
&&\bigl|\bE \bigl(\tilde{I}_{5,1}(t_i)+\tilde{I}_{5,1}(t_{i+1})
\bigr)\bigr|
\\
&&\qquad\leq C\Delta^{3/2}
\\
&&\qquad\quad{}+\biggl|\bE\int_{t_{i}}^{t_{i+1}}\mu_n({t}_{i-2})
\bigl(2{\bigl\langle} Y_t^n-Z_t-
\bigl(Y^n_{\bar{t}_n}-Z_{\bar{t}_n}\bigr),\sigma
\bigl(Y^n_{\bar
{t}_n}\bigr) \bigl(\d w_t-\d
w^n_t\bigr){\bigr\rangle}\\
&&\hspace*{243pt}{}+\operatorname{tr}\bigl(\sigma
\sigma^*\bigr) \bigl(Y^n_{\bar
{t}_n}\bigr)\,\d t \bigr)\biggr|.
\end{eqnarray*}
Continuing this process and using arguments similar to those used in
(\ref{sumgk2}) above, we get
\begin{eqnarray*}
&&\Biggl|\bE \Biggl(\sum_{i=0}^{2^n}
\tilde{I}_{5,1}(t_i) \Biggr)\Biggr|
\\
&&\qquad\leq C\Delta^{1/2}
\\
&&\qquad\quad{}+\biggl|\bE\int_{0}^{1}\mu_n(0)
\bigl(2{\bigl\langle}Y_t^n-Z_t-
\bigl(Y^n_{\bar
{t}_n}-Z_{\bar
{t}_n}\bigr),\sigma
\bigl(Y^n_{\bar{t}_n}\bigr) \bigl(\d w_t-\d
w^n_t\bigr){\bigr\rangle}\\
&&\hspace*{221pt}{}+\operatorname {tr}\bigl(\sigma
\sigma^*\bigr) \bigl(Y^n_{\bar{t}_n}\bigr)\,\d t \bigr)\biggr|
\\
&&\qquad\leq C\Delta^{1/2}
\\
&&\qquad\quad{}+\biggl|\bE\int_{0}^{1}\mu_n(0)
\biggl(2{\biggl\langle}\int_{0}^{t}\sigma
\bigl(Y^n_{\bar
{t}_n}\bigr) \bigl(\d w_s-\d
w^n_s\bigr),\sigma\bigl(Y^n_{\bar{t}_n}
\bigr) \bigl(\d w_t-\d w^n_t\bigr){\biggr
\rangle}\\
&&\hspace*{234pt}{} +\operatorname{tr}\bigl(\sigma\sigma^*\bigr) \bigl(Y^n_{\bar{t}_n}
\bigr)\,\d t \biggr)\biggr|
\\
&&\qquad\quad{}+\biggl|\bE\int_{0}^{1}\mu_n(0)
\biggl(2{\biggl\langle}\int_0^{\bar
{t}_n}\sigma
\bigl(Y^n_{\bar
{t}_n}\bigr) \bigl(\d w_s-\d
w^n_s\bigr),\sigma\bigl(Y^n_{\bar{t}_n}
\bigr) \bigl(\d w_t-\d w^n_t\bigr){\biggr
\rangle} \biggr)\biggr|
\\
&&\qquad\quad{}+\biggl|\bE\int_{0}^{1}\mu_n(0)
\biggl(2{\biggl\langle}\int_{\bar
{t}_n}^t\bigl(\sigma
\bigl(Y^n_s\bigr)-\sigma\bigl(Y^n_{\bar{t}_n}
\bigr)\bigr) \bigl(\d w_s-\d w^n_s\bigr),\\
&&\hspace*{186pt}\sigma\bigl(Y^n_{\bar
{t}_n}\bigr) \bigl(\d w_t-\d
w^n_t\bigr){\biggr\rangle} \biggr)\biggr|
\\
&&\qquad\quad{}+\biggl|\bE\int_{0}^{1}\mu_n(0)
\biggl(2{\biggl\langle}\int_{\bar
{t}_n}^t\bigl(\sigma
\bigl(Y^n_s\bigr)-\sigma(Z_s)\bigr)
\dot{h}_s\,\d s,\sigma\bigl(Y^n_{\bar{t}_n}\bigr)
\bigl(\d w_t-\d w^n_t\bigr){\biggr\rangle}
\biggr)\biggr|
\\
&&\qquad\quad{}+\biggl|\bE\int_{0}^{1}\mu_n(0)
\biggl(2{\biggl\langle}\int_{\bar
{t}_n}^t\bigl(\tilde
{b}\bigl(Y^n_s\bigr)-b(Z_s)\bigr)\,\d s,
\sigma\bigl(Y^n_{\bar{t}_n}\bigr) \bigl(\d w_t-\d
w^n_t\bigr){\biggr\rangle} \biggr)\biggr|
\\
&&\qquad\quad{}+\biggl|\bE\int_{0}^{1}\mu_n(0)
\biggl(2{\biggl\langle}\int_{\bar{t}_n}^t\bigl(\d \phi
^n_s-\d \psi_s\bigr),\sigma
\bigl(Y^n_{\bar{t}_n}\bigr) \bigl(\d w_t-\d
w^n_t\bigr){\biggr\rangle} \biggr)\biggr|
\\
&&\qquad\leq C\Delta^{1/2}+C\bE\sup_{1\leq k\leq2^n} \biggl(\int
_{0}^{t_{k}}\sigma\bigl(Y^n_{\bar{t}_n}
\bigr) \bigl(\d w_t-\d w^n_t\bigr)
\biggr)^2+C\Delta ^{\theta/2}
\\
&&\qquad\leq C\Delta^{\theta/2} \qquad\forall\theta\in(0,1).
\end{eqnarray*}
\upqed\end{pf}

%le3.12 #&#
\begin{lemma}\label{lemmaI5710}
      $ \bE(I_{7,2,2}(t_k)+I_{5,2,1}^5(t_k)+I_{10}(t_k))\leq C\Delta^{3/2}$.
\end{lemma}

\begin{pf}
Since $I_{7,2,2}=I_{7,2,2}^{(1)}+I_{7,2,2}^{(2)}$, where
\begin{eqnarray*}
\bE\bigl[I_{7,2,2}^{(1)}\bigr]&=&-\frac{4}{\gamma}\bE
\biggl(\int_{t_{k-1}}^{t_k}\mu _n(t)\int
_{t_{k-2}}^t{\bigl\langle}Y^n_{\bar{t}_n}-Z_{\bar{t}_n},
\sigma \bigl(Y^n_s\bigr) \bigl(\d w_s-
\d{w}_s^n\bigr){\bigr\rangle}
\\
&&\hspace*{135pt}{}\times{\bigl\langle}D\vph\bigl(Y^n_t\bigr),\sigma
\bigl(Y^n_t\bigr)\,\d w_t{\bigr\rangle} \biggr)
\\
&=&0,
\\
I_{7,2,2}^{(2)}&=& \frac{4}{\gamma}\int_{t_{k-1}}^{t_k}
\bigl(\mu _n(t)-\mu _n(\bar{t}_n)\bigr)\int
_{\bar{t}_n}^t{\bigl\langle}Y^n_{\bar{t}_n}-Z_{\bar
{t}_n},
\sigma \bigl(Y^n_s\bigr) \bigl(\d w_s-
\d{w}_s^n\bigr){\bigr\rangle}
\\
&&\hspace*{117pt}{}\times{\bigl\langle}D\vph\bigl(Y^n_t\bigr),\sigma
\bigl(Y^n_t\bigr)\,\d{w}_t^n{\bigr
\rangle}
\\
&&{}+\frac{4}{\gamma}\int_{t_{k-1}}^{t_k}
\mu_n(\bar{t}_n)\int_{\bar
{t}_n}^t{
\bigl\langle}Y^n_{\bar{t}_n}-Z_{\bar{t}_n},\bigl(\sigma
\bigl(Y^n_s\bigr)-\sigma \bigl(Y^n_{\bar
{t}_n}
\bigr)\bigr) \bigl(\d w_s-\d{w}_s^n\bigr){
\bigr\rangle}
\\
&&\hspace*{89pt}{}\times{\bigl\langle}D\vph\bigl(Y^n_t\bigr),\sigma
\bigl(Y^n_t\bigr)\,\d{w}_t^n{\bigr
\rangle}
\\
&&{}+\frac{4}{\gamma}\int_{t_{k-1}}^{t_k}
\mu_n(\bar{t}_n)\int_{\bar
{t}_n}^t{
\bigl\langle}Y^n_{\bar{t}_n}-Z_{\bar{t}_n},\sigma
\bigl(Y^n_{\bar
{t}_n}\bigr) \bigl(\d w_s-
\d{w}_s^n\bigr){\bigr\rangle}
\\
&&\hspace*{85pt}{}\times{\bigl\langle}D\vph\bigl(Y^n_t\bigr)-D\vph
\bigl(Y^n_{\bar{t}_n}\bigr),\sigma \bigl(Y^n_t
\bigr)\,\d {w}_t^n{\bigr\rangle}
\\
&&{}+\frac{4}{\gamma}\int_{t_{k-1}}^{t_k}
\mu_n(\bar{t}_n)\int_{\bar
{t}_n}^t{
\bigl\langle}Y^n_{\bar{t}_n}-Z_{\bar{t}_n},\sigma
\bigl(Y^n_{\bar
{t}_n}\bigr) \bigl(\d w_s-
\d{w}_s^n\bigr){\bigr\rangle}
\\
&&\hspace*{89pt}{}\times{\bigl\langle}D\vph\bigl(Y^n_{t_{k-2}}\bigr),\sigma
\bigl(Y^n_t\bigr)-\sigma \bigl(Y^n_{\bar
{t}_n}
\bigr){\bigr\rangle} \,\d{w}_t^n
\\
&&{}+\frac{4}{\gamma}\int_{t_{k-1}}^{t_k}
\mu_n(t_{k-2})\int_{\bar
{t}_n}^t{
\bigl\langle}Y^n_{t_{k-2}}-Z_{t_{k-2}},\sigma
\bigl(Y^n_{t_{k-2}}\bigr) \bigl(\d w_s-\d
{w}_s^n\bigr){\bigr\rangle}
\\
&&\hspace*{98pt}{}\times{\bigl\langle}D\vph\bigl(Y^n_{t_{k-2}}\bigr),\sigma
\bigl(Y^n_{t_{k-2}}\bigr){\bigr\rangle }\,\d{w}_t^n
\\
&=:&\sum_{i=1}^5 J_{7,2,2}^i.
\end{eqnarray*}

By applying Lemmas \ref{holderZ}, \ref{estimatevarphi}, \ref{mulemma}
and Proposition~\ref{phiBV}, we can get
\[
\bigl|\bE\bigl( J_{7,2,2}^i\bigr)\bigr|\leq C\Delta^{3/2},\qquad
i=1,2,3,4
\]
and
%
%e3.17 #&#
\begin{eqnarray}\qquad
\label{J7225} &&\bE\bigl( J_{7,2,2}^5|\sF_{t_{k-2}}
\bigr)
\nonumber
\\[-8pt]
\\[-8pt]
\nonumber
&&\qquad=\frac{2\Delta}{\gamma}\mu_n({t_{k-2}}) \sum
_i{\bigl\langle}D\vph\bigl(Y^n_{t_{k-2}}
\bigr),\sigma \bigl(Y^n_{t_{k-2}}\bigr)e_i{\bigr
\rangle} {\bigl\langle} Y_{t_{k-2}}^n-Z_{t_{k-2}},\sigma
\bigl(Y^n_{t_{k-2}}\bigr)e_i{\bigr\rangle}.
\end{eqnarray}

Moreover, since
\begin{eqnarray*}
I_{10}&=&-\frac{4}{\gamma}\int_{t_{k-1}}^{t_k}
\mu_n(t)\sum_i{\bigl\langle} D\vph
\bigl(Y^n_t\bigr),\sigma\bigl(Y^n_t
\bigr)e_i{\bigr\rangle} {\bigl\langle}Y_t^n-Z_t,
\sigma \bigl(Y^n_t\bigr)e_i{\bigr\rangle}\,\d t
\\
&=&-\frac{4}{\gamma}\int_{t_{k-1}}^{t_k}\bigl(
\mu_n(t)-\mu_n(\bar {t}_n)\bigr)\sum
_i{\bigl\langle}D\vph\bigl(Y^n_t
\bigr),\sigma\bigl(Y^n_t\bigr)e_i{\bigr
\rangle} {\bigl\langle }Y_t^n-Z_t,\sigma
\bigl(Y^n_t\bigr)e_i{\bigr\rangle}\,\d t
\\
&&{}-\frac{4}{\gamma}\int_{t_{k-1}}^{t_k}
\mu_n(\bar{t}_n)\sum_i{
\bigl\langle} D\vph \bigl(Y^n_t\bigr)-D\vph
\bigl(Y^n_{\bar{t}_n}\bigr),\sigma\bigl(Y^n_t
\bigr)e_i{\bigr\rangle} {\bigl\langle }Y_t^n-Z_t,
\sigma \bigl(Y^n_t\bigr)e_i{\bigr\rangle}\,\d t
\\
&&{}-\frac{4}{\gamma}\int_{t_{k-1}}^{t_k}
\mu_n(\bar{t}_n)\sum_i{
\bigl\langle} D\vph \bigl(Y^n_{\bar{t}_n}\bigr),\bigl(\sigma
\bigl(Y^n_t\bigr)-\sigma\bigl(Y^n_{\bar
{t}_n}
\bigr)\bigr)e_i{\bigr\rangle} {\bigl\langle} Y_t^n-Z_t,
\sigma\bigl(Y^n_t\bigr)e_i{\bigr\rangle}\,\d t
\\
&&{}-\frac{4}{\gamma}\int_{t_{k-1}}^{t_k}\hspace*{-1pt}
\mu_n(\bar{t}_n)\sum_i{
\bigl\langle} D\vph \bigl(Y^n_{\bar{t}_n}\bigr),\sigma
\bigl(Y^n_{\bar{t}_n}\bigr)e_i{\bigr\rangle} {\bigl
\langle }Y_t^n-Z_t-\bigl(Y_{\bar
{t}_n}^n-Z_{\bar{t}_n}
\bigr),\sigma\bigl(Y^n_t\bigr)e_i{\bigr
\rangle}\,\d t
\\
&&{}-\frac{4}{\gamma}\int_{t_{k-1}}^{t_k}
\mu_n(\bar{t}_n)\sum_i{
\bigl\langle} D\vph \bigl(Y^n_{\bar{t}_n}\bigr),\sigma
\bigl(Y^n_{\bar{t}_n}\bigr)e_i{\bigr\rangle} {\bigl
\langle }Y_{\bar
{t}_n}^n-Z_{\bar
{t}_n},\bigl(\sigma
\bigl(Y^n_t\bigr)-\sigma\bigl(Y^n_{\bar{t}_n}
\bigr)\bigr)e_i{\bigr\rangle}\,\d t
\\
&&{}-\frac{4\Delta}{\gamma}\mu_n(t_{k-2})\sum
_i{\bigl\langle}D\vph \bigl(Y^n_{t_{k-2}}
\bigr),\sigma\bigl(Y^n_{t_{k-2}}\bigr)e_i{\bigr
\rangle} {\bigl\langle} Y_{t_{k-2}}^n-Z_{t_{k-2}},\sigma
\bigl(Y^n_{t_{k-2}}\bigr)e_i{\bigr\rangle}
\\
&=:&\sum_{i=1}^6I_{10,i}.
\end{eqnarray*}
By Lemmas \ref{holderZ}, \ref{estimatevarphi} and \ref{mulemma} it is
easy to get
\[
\bigl|\bE(I_{10,i})\bigr|\leq C\Delta^{3/2},\qquad i=1,\ldots,5,
\]
and by using (\ref{I5215}), (\ref{J7225}),
\begin{eqnarray*}
\bE\bigl(I_{10,6}+J_{7,2,2}^5+I_{5,2,1}^5
\bigr)=\bE \bigl[\bE \bigl(I_{10,6}+J_{7,2,2}^5+I_{5,2,1}^5|
\sF_{t_{k-2}}\bigr) \bigr]=0.
\end{eqnarray*}
Thus summing all the estimates above we have
\begin{eqnarray*}
\bigl|\bE\bigl(I_{7,2,2}+I_{5,2,1}^5+I_{10}
\bigr)\bigr|\leq C\Delta^{3/2}.
\end{eqnarray*}
\upqed\end{pf}

%pr3.6 #&#
\begin{proposition}\label{mainconvergence} $
\lim_n\bE\sup_{t\in[0,T]}|Y^n_t-Z_t|^2=0$.
\end{proposition}

\begin{pf}
Take two arbitrary integers $n>n_0$, and
set $s_k=k2^{-n_0}$, $k=0,1,\ldots,  2^{n_0}-1$. Then $\{s_k\}
_{k=1}^{2^{n_0}}\subset\{t_k^n\}_{k=1}^{2^n}$. Since for every $t\in
[s_k, s_{k+1}]$
\[
\bigl|Y_t^n-Z_t\bigr|\leq\bigl|Y^n_t-Y^n_{s_k}\bigr|+
\bigl|Y^n_{s_k}-Z_{s_k}\bigr|+|Z_{s_k}-Z_t|,
\]
we have
\[
\bigl|Y_t^n-Z_t\bigr|\leq\sup_{|s-t|\leq
2^{-n_0}}
\bigl(\bigl|Y^n_t-Y^n_s\bigr|+|Z_t-Z_s|
\bigr)+\sup_{k=0,\ldots, 2^{n_0}-1}\bigl|Y^n_{s_k}-Z_{s_k}\bigr|.
\]
Consequently
\begin{eqnarray*}
\sup_{t\in[0,1]}\bigl|Y_t^n-Z_t\bigr|^2&
\leq& 3\sup_{t\in[0,1],|s-t|\leq2^{-n_0}}\bigl(\bigl|Y^n_t-Y^n_s\bigr|^2+|Z_t-Z_s|^2
\bigr)
\\
&&{}+3\sup_{k=0,\ldots, 2^{n_0}-1}\bigl|Y^n_{s_k}-Z_{s_k}\bigr|^2.
\end{eqnarray*}
Thus by Lemma~\ref{holderZ}, Proposition~\ref{convergence3} and (\ref
{holderestimatevarphi}), for any $\theta\in(0,1)$,
\begin{eqnarray*}
&&\bE\Bigl[\sup_{t\in[0,1]}\bigl|Y_t^n-Z_t\bigr|^2
\Bigr]\\
&&\qquad\leq C\bE\Bigl[\sup_{t\in[0,1],|s-t|\leq
2^{-n_0}}\bigl(\bigl|Y^n_t-Y^n_s\bigr|^2+|Z_t-Z_s|^2
\bigr)\Bigr]
\\
&&\qquad\quad{}+C\sum_{k=0}^{2^{n_0}-1}\bE
\bigl[\bigl|Y^n_{s_k}-Z_{s_k}\bigr|^2\bigr]
\\
&&\qquad\leq C\bE\Bigl[\sup_{t\in[0,1],|s-t|\leq
2^{-n_0}}\bigl(\bigl|Y^n_t-Y^n_s\bigr|^2+|Z_t-Z_s|^2
\bigr)\Bigr]
\\
&&\qquad\quad{}+C2^{n_0}\sup_{0\leq k\leq2^{n_0}}
\bE\bigl[\bigl|Y^n_{s_k}-Z_{s_k}\bigr|^2
\bigr]
\\
&&\qquad\leq C2^{-n_0\theta/2}+C2^{n_0}\biggl[2^{-n\theta/2}+\sup
_{2\leq k\leq
2^n}\biggl(\int_{(k-2)2^{-n}}^{k2^{-n}}|
\dot{h}_s|^2\,\d s\biggr)^{1/2}\biggr].
\end{eqnarray*}
Letting first $n\to\infty$ and then $n_0\to\infty$ gives the result.
\end{pf}

\begin{pf*}{Proof of Theorem~\ref{support}}
As we have already noticed, we only need to prove~(\ref{convergence2}),
and this follows from Proposition~\ref{millet2}, Lemma~\ref{holderZ},
Propositions \ref{holderC} and~\ref{mainconvergence}. We have completed the proof of
Theorem~\ref{support}.
\end{pf*}

%s4 #&#
\section{An elementary application}\label{sec4}
In this section we give an elementary application of the support
theorem to maximum principle. Recall that a typical application of the
support theorem for ordinary diffusions is the maximum principle for
degenerate elliptic and parabolic operators of second order, so it
comes with no surprise that the support theorem for reflected
diffusions is applicable to obtain a boundary-interior maximum
principle for the same operators.

%de4.1 #&#
\begin{definition}\label{subharmonic}
Let $A$ be the generator of the Markov family $\{X_t(x)\}$, the
solution to equation (\ref{msde1}). A function $u$ defined in $\bar{D}$
is said to be $A$-subharmonic in $\bar{D}$ if:
\begin{longlist}[(ii)]
\item[(i)] it is locally bounded and upper semicontinuous;

\item[(ii)] $t\rightarrow u(X_t(x))$ is a local submartingale for every $x\in
\bar{D}$.
\end{longlist}
\end{definition}

For $x\in\bar{D}$, set
\begin{eqnarray*}
\sP&:=&\bigl\{h\in\cC\bigl([0,+\infty);\mR^d\bigr);
h_0=0, t\rightarrow h_t \mbox{ is smooth}\bigr\},
\\
D(x)&:=&\overline{\bigl\{y\in\bar{D},\exists h\in\sP, t_0>0, \mbox{s.t. }
y=Z_{t_0}(x,h)\bigr\}},
\end{eqnarray*}
where $Z_t(x,h)$ solves the following deterministic Skorohod problem in $D$:
\begin{eqnarray*}
Z_t(x,h)&=&x+\int_0^t\sigma
\bigl(Z_s(x,h)\bigr)\dot{h}_s\,\d s+\int
_0^tb\bigl(Z_s(x,h)\bigr)\,\d s+
\psi_t, \\
 Z_0(x,h)&=&x\in\bar{D}.
\end{eqnarray*}

We have the following theorem, the proof of which is a modification of
\cite{iw}, Theorem~6.8.3.

%th4.1 #&#
\begin{theorem}\label{maximum}
Let $u$ be an $A$-harmonic function on $D$ and $x\in\bar{D}$. If
$u(x)=\max_{y\in D(x)} u(y)$, then $u\equiv u(x)$ on $D(x)$.
\end{theorem}

\begin{pf}
Let $y\in D(x)$ and $\{\tau_n\}$ be a localization sequence of
stopping times for $\{u(X_t(x))\}$. For every $m$ we set
\[
\varsigma_m:=\inf\bigl\{t\dvtx \bigl|X_t(x)-y\bigr|
\leq2^{-m}\bigr\}.
\]
Then for every $m$ there exists an $N_m$ such that $\forall n>N_m$,
%
%e4.1 #&#
\begin{equation}
\label{varsigmatau} \bP(\varsigma_m<\tau_n)>0.
\end{equation}

By the submartingale property we have
\[
u(x)\leq\bE \bigl[u\bigl(X_{t\wedge\varsigma_m\wedge\tau_n}(x)\bigr) \bigr].
\]
Note that since
\[
\bP \bigl(X_{t\wedge\varsigma_m\wedge\tau_n}(x)\in D(x)\bigr) =1,
\]
we have
\[
\bP \bigl(u\bigl(X_{t\wedge\varsigma_m\wedge\tau_n}(x)\bigr)=u(x) \bigr)=1.
\]
Therefore for $n>N_m$,
\[
u\bigl(X_{t\wedge\varsigma_m}(x)\bigr)\mathbh{1}_{\{\varsigma_m<\tau_n\}
}=u(x)
\mathbh{1}_{\{\varsigma_m<\tau_n\}}.
\]
This together with (\ref{varsigmatau}) implies that for every $m$ there
exists $y_m$ such that
\begin{eqnarray*}
|y_m-y|=2^{-m},\qquad u(y_m)=u(x).
\end{eqnarray*}
Hence by the upper semicontinuity property of $u$,
\begin{eqnarray*}
u(y)\geq\limsup_{m\rightarrow\infty}u(y_m)=u(x),
\end{eqnarray*}
implying that $u(y)=u(x)$.
\end{pf}

\begin{example*}
Suppose $u\in\cC^2(D)\cap\cC(\bar{D})$, and let
\begin{eqnarray*}
Lu=\sum_{i,j}a^{ij}(x)\frac{\partial^2 u}{\partial x^i\,\partial
x^j}+
\sum_ib^i(x)\frac{\partial u}{\partial x^i},\qquad  x\in D.
\end{eqnarray*}
If
\begin{eqnarray*}
Lu(x)&\geq&0 \qquad\forall x\in D,\qquad
\\
\frac{\partial u(x)}{\partial n}&\geq&0,\qquad  x\in\partial D
\end{eqnarray*}
holds in the viscosity sense, then $u$ is $A$-subharmonic.
\end{example*}

\section*{Acknowledgments}
A part of this work was done when the first author
visited the Institute of Mathematics, Academia Sinica and National Central
University. He would like to express his hearty thanks to the two
institutions, in particular to Professors
Tzuu-Shuh Chiang, Shuenn-Jyi Sheu and Yunshyong Chow for their warm
hospitality and helpful discussions. Both authors are very grateful to
the referees for their careful reading of the manuscript and valuable
suggestions.

% imsref loaded by akundreckaite, 2015-04-07 14:15:53

%
%\begin{appendix}
%\section{}
%\end{appendix}

%\section*{Acknowledgments}

%\begin{supplement}[id=suppA]
%\slink[doi]{10.1214/00-AOPXXXXSUPP} %[doi,text={...}] - jei reikia
%suskaldyti doi
%\sfilename{aopXXXX\_supp.pdf}
%\end{supplement}

%\begin{thebibliography}{99}
%\bibitem[\protect\citeauthoryear{}{}]{r1}
%\bibitem{r1}
%\end{thebibliography}

\printaddresses
\end{document}